\pdfoutput=1
\RequirePackage{ifpdf}
\ifpdf 
\documentclass[pdftex]{sigma}
\else
\documentclass{sigma}
\fi

\usepackage{bbm}
\numberwithin{equation}{section}

\newtheorem{Theorem}{Theorem}[section]
\newtheorem{Corollary}[Theorem]{Corollary}
\newtheorem{Lemma}[Theorem]{Lemma}
\newtheorem{Proposition}[Theorem]{Proposition}

 { \theoremstyle{definition}
\newtheorem{Definition}[Theorem]{Definition}
\newtheorem{Convention}[Theorem]{Convention}

\newtheorem{Remark}[Theorem]{Remark} }

\begin{document}
\allowdisplaybreaks

\newcommand{\arXivNumber}{2405.10609}

\renewcommand{\PaperNumber}{086}

\FirstPageHeading

\ShortArticleName{Quasi-Polynomial Extensions of Nonsymmetric Macdonald--Koornwinder Polynomials}

\ArticleName{Quasi-Polynomial Extensions of Nonsymmetric\\ Macdonald--Koornwinder Polynomials}

\Author{Jasper STOKMAN}

\AuthorNameForHeading{J.~Stokman}

\Address{KdV Institute for Mathematics, University of Amsterdam, \\ Science Park 105-107, 1098 XG Amsterdam, The Netherlands}
\Email{\href{mailto:j.v.stokman@uva.nl}{j.v.stokman@uva.nl}}
\URLaddress{\url{https://staff.fnwi.uva.nl/j.v.stokman/}}

\ArticleDates{Received March 26, 2025, in final form October 06, 2025; Published online October 14, 2025}

\Abstract{In a recent joint paper with S.~Sahi and V.~Venkateswaran (2025), families of actions of the double affine Hecke algebra on spaces of quasi-polynomials were introduced. These so-called quasi-polynomial representations led to the introduction of quasi-polynomial extensions of the nonsymmetric Macdonald polynomials, which reduce to metaplectic Iwahori--Whittaker functions in the $\mathfrak{p}$-adic limit. In this paper, these quasi-polynomial representations are extended to Sahi's $5$-parameter double affine Hecke algebra, and the quasi-polynomial extensions of the nonsymmetric Koornwinder polynomials are introduced.}


\Classification{33D80; 20C08}

{
\renewcommand{\theTheorem}{\Alph{Theorem}}

\section{Introduction}

\newcounter{paranum}[section]      \renewcommand{\theparanum}{\thesection.\arabic{paranum}}      \newcommand{\mypar}{\vspace{10pt}\noindent\refstepcounter{paranum}\textbf{\theparanum.}\ }

\mypar
The double affine Hecke algebra $\mathbb{H}$ with adjoint root data depends on a deformation parameter~$q$ and on a number of Hecke parameters. The Hecke parameters are most conveniently encoded by a multiplicity function, which is an affine Weyl group invariant function on the associated reduced affine root system. Cherednik's polynomial representation is a faithful representation of $\mathbb{H}$ on Laurent polynomials in several variables, given explicitly in terms of Demazure--Lusztig operators. Up to a multiplicative scalar, the nonsymmetric Macdonald polynomials can be characterised as the simultaneous eigenfunctions for the action of Bernstein's~\cite{Lu} commuting elements $Y^\lambda$ within the (double) affine Hecke algebra. See the monographs by Cherednik~\cite{ChBook} and Macdonald~\cite{Ma} for details and further references.

If the underlying finite root system is of type ${\rm C}_r$, then a nonreduced extension of the Che\-red\-nik--Macdonald theory was developed in~\cite{No,Sa,St}. It depends on five Hecke parameters (four when $r=1$), which are encoded by a multiplicity function
on Macdonald's~\cite{MacAff} nonreduced affine root system of type \smash{${\rm C}^\vee{\rm C}_r$}.
The resulting simultaneous polynomial eigenfunctions of the $Y^\lambda$ are Sahi's~\cite{Sa} nonsymmetric Koornwinder polynomials.
Their symmetric versions are the celebrated Koornwinder polynomials~\cite{Ko}, which reduce to Askey--Wilson polynomials~\cite{AW} for $r=1$.

In a joint paper~\cite{SSV2} with Sahi and Venkateswaran, a {\it quasi-polynomial} extension of the Cherednik--Macdonald theory was developed when the multiplicity function on the reduced affine root system is invariant for the action of the extended affine Weyl group (so it depends on one Hecke parameter if the underlying finite root system has a single Weyl group orbit, and two otherwise).
The role of the polynomial representation is then replaced by explicit families of $\mathbb{H}$-representations on spaces of quasi-polynomials, which are linear combinations of monomials with possibly non-integral exponents, with the action given in terms of truncated versions of Demazure--Lusztig operators.
The resulting quasi-polynomial extensions of the nonsymmetric Macdonald polynomials are $q$-analogs
of Iwahori--Whittaker functions on metaplectic covers of reductive groups over non-Archimedean local fields. Metaplectic Iwahori--Whittaker functions have been studied from the perspective of Hecke algebras in~\cite{BBBF,CGP,PP,PP2,SSV2}.
In this context, the truncated Demazure--Lusztig operators reduce to metaplectic Demazure operators (see~\cite{SSV2}).

In this paper, we will construct the nonreduced extensions of the quasi-polynomial representations and introduce the quasi-polynomial analogs of the nonsymmetric Koornwinder polynomials.
We will proceed by extending the framework for the quasi-polynomial theory in such a~way that it gives the nonreduced quasi-polynomial theory when the underlying finite root system is of type ${\rm C}_r$.
For other types, it will simply reduce to the quasi-polynomial theory from~\cite{SSV2}.
 The setup of the extended framework is modelled by the treatment of the twisted polynomial theory with adjoint root datum from~\cite{StCh}.

In the remainder
of the introduction, we will explicitly state the main results when the underlying root system is of type ${\rm C}_r$. In Section~\ref{prel}, we will introduce affine root systems and the double affine Hecke algebra $\mathbb{H}$ in the general, extended framework.
Following~\cite{SSV1}, we start Section~\ref{qpsect} by
introducing an $H^X$-action in terms of truncated Demazure--Lusztig type operators on the space of all quasi-polynomials, where $H^X$ is the copy of the affine Hecke algebra inside~$\mathbb{H}$ that contains the monomials.
This $H^X$-representation is reducible, with subrepresentations being naturally parametrised by affine Weyl group orbits in the ambient Euclidean space of the root system.
Following~\cite{SSV2}, we then give for each subrepresentation a multiparameter extension of the $H^X$-action to an action of the double affine Hecke algebra.
It gives the quasi-polynomial representations from~\cite[Section~4]{SSV2} as well as the new, nonreduced extensions when the underlying root system is of
type~${\rm C}_r$.
In Section~\ref{qKsection}, we introduce the quasi-polynomial extensions of the nonsymmetric Macdonald--Koornwinder polynomials. Finally, in Section~\ref{InducedModuleSection}, we~identify the quasi-polynomial representations with $Y$-parabolically induced $\mathbb{H}$-modules.

 In~\cite[Section~6]{SSV2}, various additional properties of the quasi-polynomial extensions of the nonsymmetric Macdonald polynomials were derived, such as creation formulas, (anti)sym\-me\-tric versions of the quasi-polynomials, and orthogonality relations. It is straightforward to derive the analogous results for the nonreduced extension of the type ${\rm C}_r$ quasi-polynomials using the intertwiners of Sahi's~\cite{Sa} nonreduced extension of the double affine Hecke algebra, but we will not discuss the details in this paper.
We also do not discuss the theory for extended lattices, which follows quite easily from the theory for adjoint root data, cf.~\cite[Section~7]{SSV2}.

\mypar\label{1.2}
Let $\mathbf{F}$ be a field of characteristic zero. The algebra of quasi-polynomials~\cite{SSV2} in $r$ variables over~$\mathbf{F}$ is the group algebra $\mathbf{F}[\mathbb{R}^r]$ of the Euclidean space $\mathbb{R}^r$, viewed as additive group. We write~$x^y$ for the standard basis element in $\mathbf{F}[\mathbb{R}^r]$ associated to
the vector $y=(y_1,\dots,y_r)\in\mathbb{R}^r$, so that
\[
\mathbf{F}[\mathbb{R}^r]=\bigoplus_{y\in\mathbb{R}^r}\mathbf{F}x^y,\qquad x^yx^{y^\prime}=x^{y+y^\prime},\qquad x^0=1.
\]
Denote by $\{\epsilon_i\}_{i=1}^r$ the standard orthonormal basis of $\mathbb{R}^r$. Then $x^y=x_1^{y_1}\cdots x_r^{y_r}$ with
\[
x_i^{\xi}:=x^{\xi\epsilon_i}\in\mathbf{F}[\mathbb{R}^r]\qquad \text{for}\, \xi\in\mathbb{R}.
\]
We call $x^y$ the quasi-monomial with quasi-exponent $y\in\mathbb{R}^r$.

Consider the hyperoctahedral group $S_r\ltimes (\pm 1)^r$. It is a Coxeter group with Coxeter system $\{s_1,\dots,s_r\}$ given by the simple neighboring transpositions
$s_i=(i,i+1)$ for $1\leq i<r$ and $s_r=(1,\dots,1,-1)$ (we identify $S_r$ and $(\pm 1)^r$ with the corresponding subgroups in $S_r\ltimes (\pm 1)^r$).
The Coxeter generators satisfy the type ${\rm C}_r$ braid relations
\begin{gather}
s_is_{i+1}s_i=s_{i+1}s_is_{i+1},\qquad 1\leq i<r-1,\nonumber\\
s_{r-1}s_rs_{r-1}s_r=s_rs_{r-1}s_rs_{r-1},\qquad
s_is_{i^\prime}=s_{i^\prime}s_i\qquad \text{if}\ |i-i^\prime|>1.\label{braidC}
\end{gather}
The formulas
\begin{align}
&s_iy:=(y_1,\dots,y_{i-1},y_{i+1},y_i,y_{i+2},\dots,y_r),\qquad 1\leq i<r,\nonumber\\
&s_ry:=(y_1,\dots,y_{r-1},-y_r)\label{trivact}
\end{align}
define a linear action of the hyperoctahedral group $S_r\ltimes (\pm 1)^r$ on $\mathbb{R}^r$.
The action~\eqref{trivact} naturally gives rise to a $S_r\ltimes (\pm 1)^r$-action
by algebra automorphisms on $\mathbf{F}[\mathbb{R}^r]$
by letting
$s_j$ act on the quasi-exponents of the quasi-monomials $x^y$ according to~\eqref{trivact}.
The subalgebra $\mathbf{F}\smash{\bigl[x^{\pm 1}\bigr]}$ of $\mathbf{F}[\mathbb{R}^r]$ spanned $x^\mu$, $\mu\in\mathbb{Z}^r$, is the algebra of Laurent polynomials in the variables $x_i:=x_i^1=x^{\epsilon_i}$, $1\leq i\leq r$, which inherits a $S_r\ltimes (\pm 1)^r$-action from $\mathbf{F}[\mathbb{R}^r]$.

For $\xi\in\mathbb{R}$, write
\begin{enumerate}\itemsep=0pt
\item[] $\lfloor\xi\rfloor$ for the largest integer $\leq\xi$,
\item[] $\lfloor\xi\rfloor_e$ for the largest even integer $\leq\xi$,
\item[]
$\lfloor\xi\rfloor_o$ for the largest odd integer $\leq\xi$.
\end{enumerate}

\begin{Definition}
For $1\leq i<r$, let $\nabla_i$, $\nabla_r^e$, $\nabla_r^o$ be the linear operators on $\mathbf{F}[\mathbb{R}^r]$ defined by
\begin{align*}
&\nabla_i(x^y):=\biggl(\frac{1-(x_{i+1}/x_i)^{\lfloor y_i-y_{i+1}\rfloor}}{1-x_i/x_{i+1}}\biggr)x^y,\\
&\nabla_r^e(x^y):=\biggl(\frac{1-x_r^{-\lfloor 2y_r\rfloor_e}}{1-x_r^2}\biggr)x^y,\qquad
\nabla_r^o(x^y):=\biggl(\frac{x_r-x_r^{-\lfloor 2y_r\rfloor_o}}{1-x_r^2}\biggr)x^y
\end{align*}
for $y\in\mathbb{R}^r$.
\end{Definition}

Note that the $\nabla_i$ are well defined because of the truncation by floor functions of the exponents in the numerator.
For $y=\mu\in\mathbb{Z}^r$, the $\nabla_i$ reduce to divided difference operators
\begin{equation}\label{nablapol}
\nabla_i(x^\mu)=\frac{x^\mu-x^{s_i\mu}}{1-x_i/x_{i+1}},\qquad \nabla_r^e(x^\mu)=\frac{x^\mu-x^{s_r\mu}}{1-x_r^2}=x_r^{-1}\nabla_r^o(x^\mu).
\end{equation}

For a subset $B$ in a set $X$, we denote by
\[
\chi_B\colon \ X\rightarrow \{0,1\}
\]
the indicator function of $B$.
We use the shorthand notations $\chi_e$ and $\chi_o$ for the indicator function of the even and odd integers inside $\mathbb{R}$.

\begin{Theorem}\label{TheoremA}
For $k,k_r,u_r\in\mathbf{F}^\times$, the linear operators $\mathcal{T}_1,\dots,\mathcal{T}_r$ on $\mathbf{F}[\mathbb{R}^r]$ defined by
\begin{align*}
&\mathcal{T}_i(x^y):=k^{\chi_{\mathbb{Z}}(y_i-y_{i+1})}x^{s_iy}+\bigl(k-k^{-1}\bigr)\nabla_i(x^y)\qquad \text{for}\ 1\leq i<r,\\
&\mathcal{T}_r(x^y):=k_r^{\chi_e(2y_r)}u_r^{\chi_o(2y_r)}x^{s_ry}+\bigl(k_r-k_r^{-1}\bigr)\nabla_r^e(x^y)+\bigl(u_r-u_r^{-1}\bigr)\nabla_r^o(x^y)
\end{align*}
satisfy the type ${\rm C}_r$ braid relations~\eqref{braidC} and the quadratic Hecke relations
\begin{align*}
&(\mathcal{T}_i-k)\bigl(\mathcal{T}_i+k^{-1}\bigr)=0\qquad\text{for}\ 1\leq i<r,\nonumber\\
&(\mathcal{T}_r-k_r)\bigl(\mathcal{T}_r+k_r^{-1}\bigr)=0.
\end{align*}
\end{Theorem}

Theorem~\ref{TheoremA} will be proven in Section~\ref{qpHX} (it is the special case of Theorem~\ref{theoremH} when the underlying finite root system is of type ${\rm C}_r$).
When $k_r=u_r$, Theorem~\ref{TheoremA} was obtained before in~\cite{SSV2, SSV1}. Following~\cite{SSV2}, we call the operators $\mathcal{T}_i$, $1\leq i\leq r$, truncated Demazure--Lusztig operators.

Theorem~\ref{TheoremA} provides a representation of the $2$-parameter Hecke algebra $H_0=H_0(k,k_r)$ of type ${\rm C}_r$ on $\mathbf{F}[\mathbb{R}^r]$.
Using the natural \smash{$\mathbf{F}\bigl[x^{\pm 1}\bigr]$}-module structure on $\mathbf{F}[\mathbb{R}^r]$, the $H_0$-action on $\mathbf{F}[\mathbb{R}^r]$ extends to a representation of the $3$-parameter affine Hecke algebra $\widetilde{H}=H(u_r,k,k_r)$ of type~${\rm C}_r$ \big(this uses the Bernstein presentation
of \smash{$\widetilde{H}$}\big). Note that the $\mathcal{T}_j$ preserve \smash{$\mathbb{C}\big[x^{\pm 1}\big]$}, in which case they reduce to the operators
\begin{gather}
\mathcal{T}_i(x^\mu)=kx^{s_i\mu}+(k-k^{-1})\biggl(\frac{x^\mu-x^{s_i\mu}}{1-x_i/x_{i+1}}\biggr)\nonumber\\
\hphantom{\mathcal{T}_i(x^\mu)}{}
=kx^\mu+k^{-1}\frac{1-k^2x_i/x_{i+1}}{1-x_i/x_{i+1}}\bigl(x^{s_i\mu}-x^\mu\bigr),\nonumber\\
\mathcal{T}_r(x^\mu)=k_rx^{s_r\mu}+\bigl(\bigl(k_r-k_r^{-1}\bigr)+\bigl(u_r-u_r^{-1}\bigr)x_r\bigr)\biggl(\frac{x^\mu-x^{s_r\mu}}{1-x_r^2}\biggr)\nonumber\\
\hphantom{\mathcal{T}_r(x^\mu)}{}
=k_rx^r+k_r^{-1}\frac{(1-ax_r)(1-bx_r)}{1-x_r^2}\bigl(x^{s_r\mu}-x^\mu\bigr)\label{tauipol}
\end{gather}
for $\mu\in\mathbb{Z}^r$ and $1\leq i<r$ by~\eqref{nablapol}, with $\{a,b\}=\bigl\{k_ru_r,-k_ru_r^{-1}\bigr\}$. These Demazure--Lusztig type operators on \smash{$\mathbb{C}\bigl[x^{\pm 1}\bigr]$} arise in the definition of the polynomial representation of the double affine Hecke algebra
of type ${\rm C}^\vee{\rm C}_r$, see~\cite{No,Sa}.

\mypar\label{1.3}
The affine extension of Theorem~\ref{TheoremA}
depends on various additional parameters. First of all, it depends on a dilation parameter \smash{$q^{\frac{1}{2}}\in\mathbf{F}^\times$}, which naturally appears in the following affine extension of the $S_r\ltimes (\pm 1)^r$-action on $\mathbf{F}\smash{\bigl[x^{\pm 1}\bigr]}$.

The affine Weyl group of type ${\rm C}_r$ is
$W:=(S_r\ltimes (\pm 1)^r)\ltimes\mathbb{Z}^r$, with the rightmost semidirect product defined in terms of the action~\eqref{trivact} restricted to $\mathbb{Z}^r$.
It acts by algebra automorphisms on $\mathbf{F}\smash{\bigl[x^{\pm 1}\bigr]}$ by
\begin{align}
&(s_ip)(x)=p(x_1,\dots,x_{i-1},x_{i+1},x_i,x_{i+2},\dots,x_r),\qquad 1\leq i<r,\nonumber\\
&(s_rp)(x)=p\bigl(x_1,\dots,x_{r-1},x_r^{-1}\bigr),\nonumber\\
&(\tau(\lambda)p)(x)=p(q^{-\lambda_1}x_1,\dots,q^{-\lambda_r}x_r),\qquad \lambda\in\mathbb{Z}^r,\label{Wintroaction}
\end{align}
for $p(x)=p(x_1,\dots,x_r)\in\mathbf{F}\smash{\bigl[x^{\pm 1}\bigr]}$, with $\tau(\lambda)$ the affine Weyl group element corresponding to~${\lambda\in\mathbb{Z}^r}$. Regarding $\mathbf{F}\smash{\bigl[x^{\pm 1}\bigr]}$ as the algebra of regular functions on
the $\mathbf{F}$-torus
\begin{equation}\label{F_torus_CC}
\mathbf{T}:=(\mathbf{F}^\times)^r,
\end{equation}
we may view~\eqref{Wintroaction} as the $W$-action on $\mathbf{F}\smash{\bigl[x^{\pm 1}\bigr]}$ contragredient to the left $W$-action
\begin{align*}
&s_it:=(t_1,\dots,t_{i-1},t_{i+1},t_i,t_{i+2},\dots,t_r),\qquad 1\leq i<r,\\
&s_rt:=\bigl(t_1,\dots,t_{r-1},t_r^{-1}\bigr),\\
&\tau(\lambda)t:=\bigl(q^{\lambda_1}t_1,\dots,q^{\lambda_r}t_r\bigr)
\end{align*}
on $t=(t_1,\dots,t_r)\in\mathbf{T}$.

The affine Weyl group $W$ is a Coxeter group with Coxeter system
$\{s_0,s_1,\dots,s_r\}$ containing the extra simple reflection
\[
s_0:=\tau(\epsilon_1)s_{\epsilon_1},\qquad s_{\epsilon_1}:=s_1\cdots s_{r-1}s_rs_{r-1}\cdots s_1.
\]
The simple reflection $s_0$ acts on $\mathbf{F}\smash{\bigl[x^{\pm 1}\bigr]}$ by
\[
(s_0p)(x)=p\bigl(qx_1^{-1},x_2,\dots,x_r\bigr).
\]
The braid relations involving $s_0$ are
\begin{gather}
s_0s_1s_0s_1=s_1s_0s_1s_0,\qquad
s_0s_i=s_is_0,\qquad 1<i\leq r.\label{braidCaff}
\end{gather}

The linear $S_r\ltimes (\pm 1)^r$-action~\eqref{trivact} on $\mathbb{R}^r$ extends to an affine linear $W$-action with $\tau(\lambda)$,
$\lambda\in\mathbb{Z}^r$, acting on $\mathbb{R}^r$
as translation operators,
\[
\tau(\lambda)y:=y+\lambda,\qquad y\in\mathbb{R}^r.
\]
The alcove
\[
C_+:=\bigl\{y\in\mathbb{R}^r \mid 0<y_r<y_{r-1}<\cdots<y_1<\tfrac{1}{2}\bigr\}
\]
is the intersection of the half-spaces
$\{y\in\mathbb{R}^r\mid \alpha_j(y)>0\}$, $0\leq j\leq r$, where the affine linear functionals $\alpha_j\colon \mathbb{R}^r\rightarrow\mathbb{R}$ are defined by
\[
\alpha_0(y):=1-2y_1,\qquad \alpha_i(y):=y_i-y_{i+1},\quad 1\leq i<r,\qquad \alpha_r(y):=2y_r.
\]
The closure $\overline{C_+}$ of $C_+$ in $\mathbb{R}^r$ is a fundamental domain for the $W$-action on $\mathbb{R}^r$. For a $W$-orbit~$\mathcal{O}$ in $\mathbb{R}^r$ we denote by
$c^\mathcal{O}$ the unique vector in $\mathcal{O}\cap \overline{C_+}$. Note that $\mathbb{Z}^r$ is a $W$-orbit in $\mathbb{R}^r$, and~$c^{\mathbb{Z}^r}=0$.

For a $W$-orbit $\mathcal{O}$ in $\mathbb{R}^r$, consider the free $\mathbf{F}\smash{\bigl[x^{\pm 1}\bigr]}$-submodule
\[
\mathbf{F}[\mathcal{O}]:=\bigoplus_{y\in\mathcal{O}}\mathbf{F}x^y
\]
of $\mathbf{F}[\mathbb{R}^r]$ of finite rank. The operators $\mathcal{T}_i$, $1\leq i\leq r$, preserve $\mathbf{F}[\mathcal{O}]$, and we write
\[
\mathcal{T}_i^{\mathcal{O}}:=\mathcal{T}_i\vert_{\mathbf{F}[\mathcal{O}]}
\]
for the resulting linear operators on $\mathbf{F}[\mathcal{O}]$. We now define a linear operator $\mathcal{T}_0^\mathcal{O}$ on $\mathbf{F}[\mathcal{O}]$ that will provide the
local affine extension of Theorem~\ref{TheoremA} (`local' in the sense that the operators should be restricted to $\mathbf{F}[\mathcal{O}]$). The operator $\mathcal{T}_0^\mathcal{O}$ will depend on additional parameters that lie in an~$\mathcal{O}$-dependent affine subtorus $\mathbf{T}_\mathcal{O}$ of $\mathbf{T}$. We define the affine subtorus $\mathbf{T}_\mathcal{O}$ now first.

Consider the simple co-roots \smash{$\alpha_j^\vee\in\mathbb{Z}^r\times\frac{1}{2}\mathbb{Z}$}, defined by
\[
\alpha_0^\vee=\bigl(-\epsilon_1,\tfrac{1}{2}\bigr),\qquad \alpha_i^\vee=(\epsilon_i-\epsilon_{i+1},0),\quad 1\leq i<r,\qquad \alpha_r^\vee:=(\epsilon_r,0).
\]
The (affine) subtorus $\mathbf{T}_\mathcal{O}$ is then given by
\begin{equation}\label{F_torus_CCC}
\mathbf{T}_{\mathcal{O}}:=\bigl\{\,t\in\mathbf{T} \mid t^{\alpha_j^\vee}=1\ \text{for}\ j\in\{0,\dots,r\}\ \hbox{satisfying}\ \alpha_j\bigl(c^{\mathcal{O}}\bigr)=0\bigr\},
\end{equation}
where
\[
t^{(\mu,\ell)}:=q^\ell t^\mu=q^\ell t_1^{\mu_1}\cdots t_r^{\mu_r}\qquad \textup{for}\ (\mu,\ell)\in\mathbb{Z}^r\times\tfrac{1}{2}\mathbb{Z}.
\]
Note that $\mathbf{T}_\mathcal{O}=\mathbf{T}$
if $\mathcal{O}$ is a regular $W$-orbit, while $\mathbf{T}_{\mathbb{Z}^r}=\{1_{\mathbf{T}}\}$.

For $y\in\mathbb{R}^r$, let $\mathbf{g}_y\in W$ be the unique element of minimal length such that $\mathbf{g}_y^{-1}y\in\overline{C_+}$.

\begin{Definition}
Let $\mathcal{O}$ be a $W$-orbit in $\mathbb{R}^r$. For $k_0, u_0\in\mathbf{F}^\times$ and $t\in\mathbf{T}_\mathcal{O}$, let $\mathcal{T}_0^\mathcal{O}$ be the linear operator on $\mathbf{F}[\mathcal{O}]$ defined by
\begin{equation*}
\mathcal{T}_0^\mathcal{O}(x^y):=k_0^{\chi_e(2y_1)}u_0^{\chi_o(2y_1)}(\mathbf{g}_yt)^{\epsilon_1}s_{\epsilon_1}(x^y)+\bigl(k_0-k_0^{-1}\bigr)\nabla_0^e(x^y)+\bigl(u_0-u_0^{-1}\bigr)\nabla_0^o(x^y)
\end{equation*}
for $y\in \mathcal{O}$, where $\nabla_0^e$, $\nabla_0^o$ are the linear operators on $\mathbf{F}[\mathbb{R}^r]$ defined by
\begin{equation*}
\nabla_0^e(x^y):=\biggl(\frac{1-\bigl(q^{-\frac{1}{2}}x_1\bigr)^{\lfloor -2y_1\rfloor_e}}{1-qx_1^{-2}}\biggr)x^y,\qquad
\nabla_0^o(x^y):=\biggl(\frac{q^{\frac{1}{2}}x_1^{-1}-\bigl(q^{-\frac{1}{2}}x_1\bigr)^{\lfloor -2y_1\rfloor_o}}{1-qx_1^{-2}}\biggr)x^y
\end{equation*}
for $y\in\mathbb{R}^r$.
\end{Definition}

It is straightforward to check that $\mathcal{T}_0^{\mathcal{O}}$ is a well-defined linear operator on $\mathbf{F}[\mathcal{O}]$. Furthermore, in case of the $W$-orbit $\mathcal{O}=\mathbb{Z}^r$ we have $\mathbf{g}_\mu 1_T=(q^{\mu_1},\dots,q^{\mu_r})$ for $\mu\in\mathbb{Z}^r$, and hence \smash{$\mathcal{T}_0^{\mathbb{Z}^r}$} reduces~to
\begin{align}
\mathcal{T}_0^{\mathbb{Z}^r}(x^\mu)&=k_0s_0(x^\mu)+\bigl(\bigl(k_0-k_0^{-1}\bigr)+\bigl(u_0-u_0^{-1}\bigr)q^{\frac{1}{2}}x_1^{-1}\bigr)\biggl(\frac{x^\mu-s_0(x^\mu)}{1-qx_1^{-2}}\biggr)\nonumber\\
&=k_0x^\mu+k_0^{-1}\frac{\bigl(1-cx_1^{-1}\bigr)\bigl(1-dx_1^{-1}\bigr)}{1-qx_1^{-2}}(s_0(x^\mu)-x^\mu)\label{T0Zr}
\end{align}
for $\mu\in\mathbb{Z}^r$ with $\{c,d\}=\bigl\{q^{\frac{1}{2}}k_0u_0,-q^{\frac{1}{2}}k_0u_0^{-1}\bigr\}$, which is the Demazure--Lusztig operator associated to the affine simple reflection $s_0$ appearing in the polynomial representation of the double affine Hecke algebra of type ${\rm C}^\vee{\rm C}_r$, see~\cite{No,Sa}.

\begin{Theorem}\label{TheoremB}
For $q^{\frac{1}{2}},k_0,u_0,k,k_r,u_r\in\mathbf{F}^\times$ and $t\in\mathbf{T}_{\mathcal{O}}$,
the operators $\mathcal{T}_0^{\mathcal{O}},\dots,\mathcal{T}_r^{\mathcal{O}}$
satisfy the affine type ${\rm C}_r$ braid relations~\eqref{braidC} and \eqref{braidCaff} and the Hecke relations
\begin{gather*}
\bigl(\mathcal{T}_0^{\mathcal{O}}-k_0\bigr)\bigl(\mathcal{T}_0^{\mathcal{O}}+k_0^{-1}\bigr)=0,\\
\bigl(\mathcal{T}_i^{\mathcal{O}}-k\bigr)\bigl(\mathcal{T}_i^{\mathcal{O}}+k^{-1}\bigr)=0,\\
\bigl(\mathcal{T}_r^{\mathcal{O}}-k_r\bigr)\bigl(\mathcal{T}_r^{\mathcal{O}}+k_r^{-1}\bigr)=0
\end{gather*}
for $1\leq i<r$.
\end{Theorem}

Theorem~\ref{TheoremB} will be proven in Section~\ref{wtsection} (it is the special case of
Theorem~\ref{mthm} when the underlying finite root system is of type ${\rm C}_r$).
In Section~\ref{InducedModuleSection}, we will also show that the quasi-polynomial representation is isomorphic to a $Y$-parabolically induced $\mathbb{H}$-module. These results were obtained before in~\cite{SSV2} when $k_0=u_0=k_r=u_r$.

Theorem~\ref{TheoremB} gives rise to a representation of the $3$-parameter Hecke algebra $H:=H(k_0,k,k_r)$ of type ${\rm C}_r$ on $\mathbf{F}[\mathcal{O}]$ (using now the Coxeter presentation of the affine Hecke algebra $H$). Adding the action of $\mathbf{F}\smash{\bigl[x^{\pm 1}\bigr]}$ by multiplication operators yields a representation of Sahi's~\cite{Sa} double affine Hecke algebra \smash{$\mathbb{H}=\mathbb{H}\bigl(k_0,u_0,k,k_r,u_r;q^{\frac{1}{2}}\bigr)$} of type ${\rm C}^\vee{\rm C}_r$ on $\mathbf{F}[\mathcal{O}]$, which depends on $t\in\mathbf{T}_\mathcal{O}$. We call it the quasi-polynomial representation of $\mathbb{H}$. By~\eqref{tauipol} and~\eqref{T0Zr}, the quasi-polynomial representation for $\mathcal{O}=\mathbb{Z}^r$ is the polynomial representation of $\mathbb{H}$ from~\cite{No,Sa}, which governs the Koornwinder polynomials.

When $k_0=u_0=k_r=u_r$, a particular reparametrisation of the extra parameters $t\in\mathbf{T}_{\mathcal{O}}$ in terms of so-called $g$-parameters allows to glue the quasi-polynomial representations from Theorem~\ref{TheoremB} into a family of $\mathbb{H}$-representations on $\mathbf{F}[\mathbb{R}^r]$ with the Coxeter type generators of $H$ acting by global $g$-dependent truncated Demazure--Lusztig type operators. This is an important intermediate step in establishing the link to representation theory of metaplectic covers of symplectic groups over non-Archimedean local fields when taking the $\mathfrak{p}$-adic limit $q\rightarrow\infty$ (the $g$-parameters are then given in terms of Gauss sums). In this metaplectic context the global truncated Demazure--Lusztig operators
reduce to the type ${\rm C}_r$ metaplectic Demazure--Lusztig operators from~\cite{BBBF,CGP,PP,PP2}. See~\cite{SSV2} for details. It is unknown how Theorem~\ref{TheoremB} for arbitrary parameters $k_0$, $u_0$, $k_r$, $u_r$ relates to
representation theory of metaplectic covers of symplectic groups over non-Archimedean local fields.

\mypar\label{1.4}
The quasi-polynomial extensions of the monic nonsymmetric Koornwinder polynomials are defined as follows.
For $y^\prime,y\in\mathbb{R}^r$, write
\[
y^\prime\leq y
\]
if $y^\prime$ and $y$ lie in the same $W$-orbit and $\mathbf{g}_{y^\prime}\leq_B\mathbf{g}_{y}$, where $\leq_B$ is the Bruhat order
on $W$. For a $W$-orbit $\mathcal{O}$ and an element $y\in\mathcal{O}$, the subset $\{y^\prime\in\mathbb{R}^r\mid y^\prime\leq y\}$ of $\mathcal{O}$ is finite, and $c^\mathcal{O}$ is its unique minimal element.

We say that $p(x)\in\mathbf{F}[\mathbb{R}^r]$ is a quasi-polynomial of degree $y\in\mathbb{R}^r$ if
\[
p(x)-dx^y\in\bigoplus_{y^\prime<y}\mathbf{F}x^{y^\prime}
\]
for some $d\in\mathbf{F}^\times$. We then say that $d$ is the leading term of $p(x)$, and $p(x)$ is said to be monic if $d=1$.
If $p(x)$ is of degree $y$ with $y$ lying in the $W$-orbit $\mathcal{O}$, then
$p(x)\in\mathbf{F}[\mathcal{O}]$.

Fix a $W$-orbit $\mathcal{O}$ in $\mathbb{R}^r$ and fix $t\in\mathbf{T}_{\mathcal{O}}$. Consider the invertible linear operators
\begin{equation*}
\mathcal{Y}_i^{\mathcal{O}}:=\bigl(\mathcal{T}_{i-1}^{\mathcal{O}}\bigr)^{-1}\cdots \bigl(\mathcal{T}_1^{\mathcal{O}}\bigr)^{-1}\mathcal{T}_0^{\mathcal{O}}\mathcal{T}_1^{\mathcal{O}}\cdots\mathcal{T}_{r-1}^{\mathcal{O}}\mathcal{T}_r^{\mathcal{O}}\mathcal{T}_{r-1}^{\mathcal{O}}\cdots
\mathcal{T}_i^{\mathcal{O}},\qquad 1\leq i\leq r,
\end{equation*}
on $\mathbf{F}[\mathcal{O}]$. These operators are the images under $\pi^{\mathcal{O}}$ of the commuting elements $Y^{\epsilon_i}\in H$ in the Bernstein presentation of $H$
(cf.\ Section~\ref{AHAsection}). In particular,
\[
\bigl[ \mathcal{Y}_i^{\mathcal{O}},\mathcal{Y}_j^{\mathcal{O}}\bigr]=0\qquad \text{for}\ 1\leq i,j\leq r.
\]

\begin{Theorem}\label{TheoremC}
Fix a $W$-orbit $\mathcal{O}$ and fix generic parameters $q^{\frac{1}{2}},k_0,u_0,k,k_r,u_r\in\mathbf{F}^\times$ and $t\in\mathbf{T}_{\mathcal{O}}$. For each $y\in\mathcal{O}$, there exists a unique quasi-polynomial
\[
E_y^{\mathcal{O}}(x)=E_y^{\mathcal{O}}\bigl(x;k_0,u_0,k,k_r,u_r,t;q^{\frac{1}{2}}\bigr)\in\mathbf{F}[\mathcal{O}]
\]
satisfying the following two properties:
\begin{enumerate}\itemsep=0pt
\item[$(1)$] $E_y^{\mathcal{O}}(x)$ is a monic quasi-polynomial of degree $y$.
\item[$(2)$] $E_y^{\mathcal{O}}(x)$ is a joint eigenfunction of the commuting operators $\mathcal{Y}_i^{\mathcal{O}}$, $1\leq i\leq r$.
\end{enumerate}
\end{Theorem}

We will prove Theorem~\ref{TheoremC} in Section~\ref{qKsection} (it is the special case of Theorem~\ref{EJycor} when the underlying finite root system is of type ${\rm C}_r$). It was derived before in~\cite{SSV2} when $k_0=u_0=k_r=u_r$.
The main step in proving Theorem~\ref{TheoremC} is showing that the $\mathcal{Y}_i^{\mathcal{O}}$, $1\leq i\leq r$, are triangular operators relative to
the partially ordered quasi-monomial basis $\{x^y\}_{y\in\mathcal{O}}$ of $\mathbf{F}[\mathcal{O}]$, with the partial order on $\{x^y\}_{y\in\mathcal{O}}$ induced from the partial order $\leq$ on the corresponding set $\mathcal{O}$ of quasi-exponents. Concretely, we will show that $\mathcal{Y}_i^{\mathcal{O}}(x^y)\in\mathbf{F}[\mathcal{O}]$ is a quasi-polynomial of degree $y$ with leading term \smash{$\gamma_i^{\mathcal{O}}(y)\in\mathbf{F}^\times$} given explicitly by
\begin{align*}
&\gamma_i^{\mathcal{O}}(y)=\bigl(\mathbf{g}_y\bigl(\mathfrak{s}^{\mathcal{O}}t\bigr)\bigr)_i^{-1},\\
&\mathfrak{s}^{\mathcal{O}}:=\bigl((k_0k_r)^{-\chi_e(2c_1^{\mathcal{O}})}(u_0u_r)^{\chi_o(2c_1^{\mathcal{O}})}k^{n_1^{\mathcal{O}}},
\dots,(k_0k_r)^{-\chi_e(2c_r^{\mathcal{O}})}(u_0u_r)^{\chi_o(2c_r^{\mathcal{O}})}k^{n_r^{\mathcal{O}}}\bigr),
\end{align*}
where
\[
n_i^{\mathcal{O}}:=\sum_{j=i+1}^r\bigl(\eta\bigl(c_i^{\mathcal{O}}-c_j^{\mathcal{O}}\bigr)+\eta\bigl(c_i^{\mathcal{O}}+c_j^{\mathcal{O}}\bigr)\bigr)+
\sum_{j=1}^{i-1}\bigl(\eta\bigl(c_j^{\mathcal{O}}+c_i^{\mathcal{O}}\bigr)-\eta\bigl(c_j^{\mathcal{O}}-c_i^{\mathcal{O}}\bigr)\bigr)
\]
and $\eta:=\chi_{\mathbb{Z}_{>0}}-\chi_{\mathbb{Z}_{\leq 0}}$. Then
\[
\mathcal{Y}_i^\mathcal{O}(E_y(x))=\gamma_i^\mathcal{O}(y)E_y(x)\qquad \text{for}\ 1\leq i\leq r\ \text{and}\ y\in\mathcal{O},
\]
and the generic conditions on the parameters in Theorem~\ref{TheoremC} boil down to the requirement that the map
\[
\mathcal{O}\rightarrow\mathbf{T},\qquad y\mapsto \bigl(\gamma_1^{\mathcal{O}}(y),\dots,\gamma_r^{\mathcal{O}}(y)\bigr)
\]
is an embedding.

The $E_\lambda^{\mathbb{Z}^r}(x)\in\mathbf{F}\smash{\bigl[x^{\pm 1}\bigr]}$, $\lambda\in\mathbb{Z}^r$, are Sahi's~\cite{Sa} monic nonsymmetric Koornwinder polynomials (recall that in this case necessarily $t=1_{\mathbf{T}}$). Hence the $E_y^\mathcal{O}(x)$ ($y\in\mathcal{O}$) may be viewed as quasi-polynomial generalisations of the nonsymmetric Koornwinder polynomials, depending on the extra parameters $t\in\mathbf{T}_\mathcal{O}$.
If $\mathcal{O}$ and $\mathcal{O}^\prime$ are two $W$-orbits in $\mathbb{R}^r$ intersecting $\overline{C_+}$ in the same face, then the corresponding families of quasi-polynomials are essentially the same (cf.~\cite[Theorem~6.2\,(4)]{SSV2}).

}

\section{Preliminaries}\label{prel}

\subsection{Reduced affine root systems}

Let $(E,\langle\cdot,\cdot\rangle)$ be an Euclidean space of dimension $r$. Transferring the inner product $\langle\cdot,\cdot\rangle$ on~$E$ to~$E^*$ through the linear isomorphism \smash{$E\overset{\sim}{\longrightarrow} E^*$}, $y\mapsto \langle y,\cdot\rangle$, is turning $E^*$ into an Euclidean space. We denote its inner product again by $\langle\cdot,\cdot\rangle$, and its norm by $\|\cdot\|$.

Let $\Phi_0$ be an irreducible reduced root system in $E^*$ with Weyl group $W_0$. Its dual root system $\Phi_0^\vee=\{\alpha^\vee\}_{\alpha\in\Phi_0}$ in $E$ consists of the co-roots $\alpha^\vee\in E$ ($\alpha\in\Phi_0$), which are the vectors in~$E$ satisfying
\begin{equation}\label{corootdef}
\bigl\langle y,\alpha^\vee\bigr\rangle=\frac{2\alpha(y)}{\|\alpha\|^2}
\end{equation}
for all $y\in E$.

Consider the corresponding reduced affine root system
\[
\Phi:=\Phi_0\times \mathbb{Z}\subset E^*\times\mathbb{R}.
\]
We will view an element $(\phi,\xi)\in E^*\times\mathbb{R}$ in the ambient space as an affine linear functional on~$E$ by $y\mapsto\phi(y)+\xi$ for $y\in E$.

The projection $E^*\times\mathbb{R}\rightarrow E^*$ on the first component will be denoted by
\[
f\mapsto\overline{f}.
\]
It restricts to a surjective map $\Phi\twoheadrightarrow\Phi_0$. Furthermore, we have $a=(\overline{a},a(0))$ for $a\in\Phi$.
Throughout the paper, we will identify a root $\alpha\in\Phi_0$ with $(\alpha,0)\in\Phi$.

For $a\in\Phi$, denote by $s_a\colon E\rightarrow E$ the orthogonal reflection in the affine root hyperplane $a^{-1}(0)\subset E$. Then
\begin{equation*}
s_a(y)=y-a(y)\overline{a}^\vee
\end{equation*}
for $y\in E$.
The affine Weyl group $W$ of $\Phi$ is the subgroup of affine linear transformations of $E$ generated by the orthogonal reflections $s_a$, $a\in\Phi$. The finite Weyl group $W_0$ is the subgroup generated by $s_\alpha$, $\alpha\in\Phi_0$.

For $y\in E$, let $\tau(y)\colon E\rightarrow E$ be the translation map $z\mapsto z+y$. Then
\begin{equation}\label{ssplit}
s_a=s_{\overline{a}}\,\tau\bigl(a(0)\overline{a}^\vee\bigr)
\end{equation}
for $a\in\Phi$. Consequently,
$W\simeq W_0\ltimes Q^\vee$ with $Q^\vee=\mathbb{Z}\Phi_0^\vee$ the co-root lattice of $\Phi_0$.

The linear, contragredient $W$-action
on the space $E^*\times\mathbb{R}$ of affine linear functionals on $E$ restricts to a $W$-action on $\Phi$. It satisfies
\begin{align}
&s_a(b)=b-\overline{b}\bigl(\overline{a}^\vee\bigr)a=\bigl(s_{\overline{a}}(\overline{b}),b(0)-a(0)\overline{b}\bigl(\overline{a}^\vee\bigr)\bigr),\nonumber\\
&\tau(\lambda)b=\bigl(\overline{b},b(0)-\overline{b}(\lambda)\bigr)\label{WPhiaction}
\end{align}
for $a,b\in\Phi$ and $\lambda\in Q^\vee$.

We fix an ordered basis $\Delta_0=\{\alpha_1,\dots,\alpha_r\}$ of the root system $\Phi_0$ once and for all. We will choose the ordering such that the following convention
holds true.

\begin{Convention}\label{root_length_convention}
The simple root $\alpha_r$ is a long root.
\end{Convention}

If all the roots in $\Phi_0$ have the same root length, then all roots are considered to be long as well as short.
We denote by $\Phi_0^{+}$ the set of positive roots in $\Phi_0$ relative to $\Delta_0$. The corresponding set of negative roots is denoted by $\Phi_0^-:=-\Phi_0^+$.

The Weyl group $W_0$ is a Coxeter group with Coxeter generators $\{s_1,\dots,s_r\}$ given by the simple reflections $s_i:=s_{\alpha_i}$, $1\leq i\leq r$. The closure of the positive Weyl chamber
\begin{equation*}
E_+:=\bigl\{y\in E \mid \alpha(y)>0\ \forall \alpha\in\Phi_0^+ \bigr\}
\end{equation*}
is a fundamental domain for the $W_0$-action on $E$.

The ordered basis $\Delta_0$ of $\Phi_0$ extends to an ordered basis
$\Delta=\{\alpha_0,\dots,\alpha_r\}$ of $\Phi$ with the additional affine simple root
\[
\alpha_0=(-\varphi,1),
\]
where $\varphi$ is the highest root of $\Phi_0$ relative to $\Delta_0$. The corresponding sets of positive and negative roots are denoted by $\Phi^+$ and $\Phi^-$, respectively.
The affine Weyl group $W$ is a Coxeter group with Coxeter generators $\{s_0,\dots,s_r\}$ the simple reflections $s_j:=s_{\alpha_j}$, $0\leq j\leq r$. By~\eqref{ssplit}, we~have
\[
s_0=s_\varphi\tau\bigl(-\varphi^\vee\bigr)=\tau\bigl(\varphi^\vee\bigr)s_\varphi.
\]

The closure $\overline{C_+}$ of the fundamental alcove
\begin{equation*}
C_+:=\{y\in E_+ \mid \alpha_0(y)>0\}
\end{equation*}
is a fundamental domain for the action of $W\simeq W_0\ltimes Q^\vee$ on $E$ by reflections and translations.
For a $W$-orbit $\mathcal{O}$ in $E$, we denote by
$c^\mathcal{O}$ the unique vector in $\mathcal{O}\cap \overline{C_+}$.

Since $\alpha_0=(-\varphi,1)$, we have the following alternative description:
\[
C_+=\bigl\{y\in E \mid 0<\alpha(y)<1 \ \forall \alpha\in\Phi_0^+\bigr\}
\]
of the fundamental alcove.
Note furthermore that
\begin{equation}\label{W0C}
\bigcup_{w\in W_0}w\bigl(\overline{C_+}\bigr)=\{y\in E \mid |\alpha(y)|\leq 1\ \forall \alpha\in\Phi_0\}.
\end{equation}

\subsection{Nonreduced extensions and multiplicity functions}

For $a\in\Phi$ such that $\overline{a}\smash{\bigl(Q^\vee\bigr)}=\mathbb{Z}$, we have
\[
Wa=W_0 \overline{a}\times\mathbb{Z}
\]
by~\eqref{WPhiaction}. A case by case inspection of the Dynkin diagrams shows that $\alpha\smash{\bigl(Q^\vee\bigr)}=\mathbb{Z}$ for $\alpha\in\Phi_0$ unless $\alpha\in\Phi_0$ is long and
 $\Phi_0$ is of type ${\rm C}_r$, $r\geq 1$, in which case $\alpha\smash{\bigl(Q^\vee\bigr)}=2\mathbb{Z}$ (note that ${\rm C}_1=\textup{A}_1$ and ${\rm C}_2=\textup{B}_2$).
If $\Phi_0$ is of type ${\rm C}_r$, $r\geq 1$, then $\alpha_r$ is the only long simple root in $\Delta_0$ in view of convention~\ref{root_length_convention}.

The $W$-orbits in $\Phi$ can now be described as follows:
\begin{enumerate}\itemsep=0pt
\item[$(1)$] If all the roots in $\Phi_0$ have the same root length but $\Phi_0$ is not of type $\textup{A}_1$, then $W$ acts transitively on $\Phi$.
\item[$(2)$] If $\Phi_0$ is of type ${\rm C}_1=\textup{A}_1$, then $\Phi=W\alpha_0\sqcup W\alpha_1$ and
\[
W\alpha_0=\Phi_0\times\mathbb{Z}_o,\qquad W\alpha_1=\Phi_0\times\mathbb{Z}_e
\]
with $\mathbb{Z}_o$ (resp.\ $\mathbb{Z}_e$) the set of odd (resp.\ even) integers.
\item[$(3)$] If $\Phi_0$ is of type ${\rm C}_r$, $r\geq 2$, then $\Phi=W\alpha_0\sqcup W\alpha_1\sqcup W\alpha_r$ and $\alpha_i\in W\alpha_1$ for all $1\leq i<r$. Furthermore,
\[
W\alpha_0=\Phi_0^\ell\times\mathbb{Z}_o,\qquad W\alpha_1=\Phi_0^s\times\mathbb{Z},\qquad W\alpha_r=\Phi_0^\ell\times\mathbb{Z}_e
\]
with $\Phi_0^\ell$ (resp.\ $\Phi_0^s$) the long (resp.\ short) roots in $\Phi_0$.
\item[$(4)$] If $\Phi_0$ is of type $\textup{B}_r$, $r\geq 3$, $F_4$ or $G_2$ and if $\alpha_i\in\Delta_0$, $1\leq i<r$, is a short simple root, then $\Phi=W\alpha_0\sqcup W\alpha_i$ and
\[
W\alpha_0=\Phi_0^\ell\times\mathbb{Z}=W\alpha_r,\qquad W\alpha_i=\Phi_0^s\times\mathbb{Z}.
\]
\end{enumerate}
The set
\[
\Phi^{\textup{nr}}:=\Phi\sqcup\bigl\{a/2\mid a\in\Phi\ \textup{such that}\ \overline{a}\bigl(Q^\vee\bigr)=\mathbb{Z}_e\big\}
\]
forms an affine root system in the affine space $E^*\times\mathbb{R}$ (see~\cite{MacAff}). If $\Phi_0$ is of type $C_r$, $r\geq 1$, then $\Phi^{\textup{nr}}$ is the nonreduced irreducible affine root system of type ${\rm C}^\vee{\rm C}_r$ (see~\cite{MacAff}). In this case, $\Phi^{\textup{nr}}$~has five $W$-orbits $W\alpha_0$, $W\frac{\alpha_0}{2}$, $W\alpha_1$, $W\alpha_r$, $W\frac{\alpha_r}{2}$ when $r\geq 2$, and four $W$-orbits
$W\alpha_0$, $W\frac{\alpha_0}{2}$, $W\alpha_1$, $W\frac{\alpha_1}{2}$ when $r=1$. If $\Phi_0$ is not of type ${\rm C}_r$, $r\geq 1$, then $\Phi^{\textup{nr}}=\Phi$.

Let $\mathbf{F}$ be a field of characteristic zero. We call a $W$-invariant function
\[
\mathbf{k}\colon\ \Phi^{\textup{nr}}\rightarrow\mathbf{F}^\times,\qquad a\mapsto \mathbf{k}_a
\]
a multiplicity function. Denote by $\mathcal{K}$ the set of multiplicity functions.
In order to obtain uniform notations, we extend a multiplicity function
$\mathbf{k}\colon \Phi^{\textup{nr}}\rightarrow\mathbf{F}^\times$ to a $W$-invariant function
\[
\mathbf{k}\colon\ \Phi\sqcup\tfrac{1}{2}\Phi\rightarrow\mathbf{F}^\times
\]
by declaring \smash{$\mathbf{k}_{\frac{a}{2}}:=\mathbf{k}_{a}$} when \smash{$\frac{a}{2}\not\in\Phi^{\textup{nr}}$}. Note that for any root $\alpha\in\Phi_0$,
\begin{equation}
\label{equalparameter}
\mathbf{k}_\alpha=\mathbf{k}_{(\alpha,1)}=\mathbf{k}_{\frac{\alpha}{2}}=\mathbf{k}_{(\frac{\alpha}{2},\frac{1}{2})}\qquad \textup{if $\Phi_0$ is not of type ${\rm C}_r$, $r\geq 1$}.
\end{equation}
For the value of a multiplicity function $\mathbf{k}$ at a simple root $\alpha_j$ and at \smash{$\frac{\alpha_j}{2}$}, we will use the shorthand notations
\[
k_j:=\mathbf{k}_{\alpha_j},\qquad u_j:=\mathbf{k}_{\frac{\alpha_j}{2}}.
\]

If $\Phi_0$ has rank $r=1$, then
the multiplicity function $\mathbf{k}$ is determined by the four parameters $k_0$, $u_0$, $k_1$, $u_1$, which can be chosen arbitrarily.
If $\Phi_0$ has rank $r>1$, then $\mathbf{k}$ is determined by the five parameters $k_0$, $u_0$, $k:=k_i$, $k_r$, $u_r$, where $1\leq i<r$ is such that $\alpha_i$ is a short root. These parameters can be chosen arbitrarily
when $\Phi_0$ is of type~${\rm C}_r$, $r\geq 2$. If $\Phi_0$ is not of type~${\rm C}_r$, $r\geq 1$,
then $k_0=u_0=k_r=u_r$ by~\eqref{equalparameter}, hence $\mathcal{K}\simeq\mathbf{F}^\times$ if all the roots in $\Phi_0$ have the same length and $\mathcal{K}\simeq (\mathbf{F}^\times)^2$ otherwise (i.e., if $\Phi_0$ is of type $\textup{B}_r$, $r\geq 3$, $F_4$ or $G_2$).

The extended affine Weyl group is the subgroup
\[
W^{\textup{ext}}:=W_0\tau\bigl(P^\vee\bigr)
\]
of the group of affine linear transformations of $E$,
with $P^\vee$ the co-weight lattice of $\Phi_0$. The linear, contragredient $W^{\textup{ext}}$-action on the space $E^*\times\mathbb{R}$ of affine linear functionals on $E$ restricts to a $W^{\textup{ext}}$-action on $\Phi$. The explicit formulas for this action are again given by~\eqref{WPhiaction}, now with $\lambda\in P^\vee$ in the second formula.

Write $\mathcal{K}^{{\rm res}}\subseteq\mathcal{K}$ for the subset of multiplicity functions $\mathbf{k}$ satisfying the following
two additional conditions:
\begin{enumerate}\itemsep=0pt
\item[$ (1)$] $\mathbf{k}$ is $W^{\textup{ext}}$-invariant,
\item[$(2)$] $\mathbf{k}_{\frac{a}{2}}=\mathbf{k}_a$ for all $a\in\Phi$.
\end{enumerate}
By $(2)$, a restricted multiplicity function $\mathbf{k}\in\mathcal{K}^{{\rm res}}$ is uniquely determined by its values on $\Phi$, and $\mathbf{k}_a=\mathbf{k}_{\overline{a}}$ for $a\in\Phi$. Its value $\mathbf{k}_\alpha$ at $\alpha\in\Phi_0$ only depends on the length of $\alpha$. Hence $\mathcal{K}^{{\rm res}}\simeq\mathbf{F}^\times$ if all roots in $\Phi_0$ have the same length and
$\mathcal{K}^{{\rm res}}\simeq (\mathbf{F}^\times)^2$ otherwise, which implies that
\[
\mathcal{K}^{{\rm res}}=\mathcal{K}\qquad \textup{if }  \Phi_0 \ \textup{is not of type}\ {\rm C}_r,\ r\geq 1.
\]
In particular, for $\mathbf{k}\in\mathcal{K}^{{\rm res}}$ formula~\eqref{equalparameter} holds true for root systems $\Phi_0$ of any type,
\[
k_0=u_0=k_r=u_r\qquad \text{if}\ \mathbf{k}\in\mathcal{K}^{{\rm res}}.
\]

Define an involution
\begin{equation*}
\mathcal{K}\overset{\sim}{\longrightarrow}\mathcal{K},\qquad \mathbf{k}\mapsto \widetilde{\mathbf{k}}
\end{equation*}
by interchanging the values $k_0$ and $u_r$ of $\mathbf{k}$ on the $W$-orbits $W\alpha_0$ and \smash{$W\frac{\alpha_r}{2}$}. We call it the duality involution.
Note that it is the identity unless $\Phi_0$ is of type ${\rm C}_r$, $r\geq 1$. Furthermore, its restriction to $\mathcal{K}^{{\rm res}}$ is the identity for all types.

The standard Cherednik--Macdonald theory for parameters $\mathbf{k}\in\mathcal{K}^{{\rm res}}$ admits an extension to
parameters $\mathbf{k}\in\mathcal{K}$~\cite{No,Sa}. It only gives new results when $\Phi_0$ is of type ${\rm C}_r$, $r\geq 1$, since otherwise $\mathcal{K}=\mathcal{K}^{{\rm res}}$.
The resulting theory is sometimes referred to as the Koornwinder case \big(since the associated analogs of the symmetric Macdonald polynomials are the Koornwinder polynomials~\cite{Ko}\big), or as the ${\rm C}^\vee{\rm C}_r$ case (since $\mathcal{K}$ is the natural set of multiplicity functions on the nonreduced root system of type
${\rm C}^\vee{\rm C}_r$).
It is common in the literature on Koornwinder~$\big({\rm C}^\vee{\rm C}_r\big)$ extensions of the Cherednik--Macdonald theory
to develop the theory directly using the following explicit realisation of the root system $\Phi_0$ of type ${\rm C}_r$, $r\geq 1$, see, e.g.,~\cite{Ko,No,Sa,St}:
\begin{itemize}\itemsep=0pt
\item $E=\mathbb{R}^r$ with orthonormal basis $\{\epsilon_i\}_{i=1}^r$,
\item $\Phi_0^s=\{\pm(\epsilon_i\pm\epsilon_j)\}_{1\leq i<j\leq r}$ ($=\varnothing$ when $r=1$) and $\Phi_0^\ell=\{\pm 2\epsilon_i\}_{i=1}^r$, where we identify~${E\simeq E^*}$
via the scalar product \big(in particular, $Q^\vee=\bigoplus_{i=1}^r\mathbb{Z}\epsilon_i$\big),
\item $\alpha_i=\epsilon_i-\epsilon_{i+1}$, $1\leq i<r$, and $\alpha_r=2\epsilon_r$.
\end{itemize}
Note that
$\varphi=2\epsilon_1$ is the highest root, hence $\alpha_0=(-2\epsilon_1,1)$. With these choices, the type~${\rm C}_r$ results presented in Sections~\ref{1.2}--\ref{1.4} follow immediately from the general results as discussed below. In the remainder of the paper, $\mathbf{k}$ will be a multiplicity function in $\mathcal{K}$ unless stated explicitly otherwise.

\begin{Remark}
The setup in this subsection follows~\cite{StAM, StCh}. In terms of the initial data $D$ from~\cite[Section~1.1]{StAM}, we are considering the case of twisted adjoint root data $D=(R_0,t,\Lambda,\Lambda)$ with $\Lambda$ the root lattice of $R_0$. Then
$(\Phi_0,\Phi^{\textup{nr}})$ corresponds to $\bigl(R_0^\vee,R(D)^\vee\bigr)$, with $R(D)^\vee$ the dual of the (possibly non-reduced) affine root system $R(D)$ from~\cite[Section~1.1]{StAM}. All nonreduced cases of the Cherednik--Macdonald theory as described in Macdonald's book~\cite{Ma} can be recovered from the case that $\Phi_0$ is of type ${\rm C}_r$ by appropriate specialisations of the multiplicity parameters, see, e.g.,~\cite[Section~9.2.3]{StCh} for details.
\end{Remark}

\subsection{Quasi-polynomials}
The space of quasi-polynomials~\cite{SSV2} is the group algebra
\[
\mathbf{F}[E]=\bigoplus_{y\in E}\mathbf{F}x^y
\]
of $E$, viewed as abelian additive group. Here we denote the canonical basis elements $x^y$, ${y\in E}$, of~$\mathbf{F}[E]$ multiplicatively, so \smash{$x^y x^{y^\prime}=x^{y+y^\prime}$} and $x^0$ is the unit element. We call $x^y$ the quasi-monomial with quasi-exponent $y\in E$. The Weyl group $W_0$ acts on $\mathbf{F}[E]$
by $\mathbf{F}$-algebra automorphisms by
\[
w(x^y):=x^{wy}
\]
for $w\in W_0$ and $y\in E$.

For any subset $Z\subseteq E$, we write
\[
\mathbf{F}[Z]:=\bigoplus_{y\in Z}\mathbf{F}x^y
\]
for the subspace of $\mathbf{F}[E]$ spanned by the quasi-monomials $x^y$, $y\in Z$.

The subspace $\mathbf{F}[P^\vee]$ is a $W_0$-stable subalgebra of $\mathbf{F}[E]$, which we call the subalgebra of Laurent polynomials in $\mathbf{F}[E]$. For $\mu\in P^\vee$, we say that $x^\mu$ is a monomial with exponent $\mu\in P^\vee$.
For generic multiplicity parameter $\mathbf{k}\in\mathcal{K}^{{\rm res}}$, the nonsymmetric Macdonald polynomials form a~basis of $\mathbf{F}[P^\vee]$.
In this paper, we are focussing on the theory for the extended set $\mathcal{K}$ of multiplicity parameters, in which case the corresponding nonsymmetric Macdonald--Koornwinder polynomials form a basis of the $W_0$-stable subalgebra $\mathbf{F}\smash{\bigl[Q^\vee\bigr]}$ of $\mathbf{F}[P^\vee]$.

Let $E/W$ be the set of $W$-orbits in $E$. For a $W$-orbit $\mathcal{O}\in E/W$, we denote by $c^{\mathcal{O}}$ the unique point in the intersection of $\mathcal{O}$
and the closure $\overline{C_+}$ of the fundamental alcove. Note that the $W$-orbit in $E$ containing the origin is $Q^\vee$.
A~$W$-orbit $\mathcal{O}$ in $E$ has a finite number of $\tau\smash{\bigl(Q^\vee\bigr)}$-orbits $\tau\smash{\bigl(Q^\vee\bigr)}y_i$, $1\leq i\leq N$, hence the corresponding space $\mathbf{F}[\mathcal{O}]$ of quasi-polynomials with quasi-exponents in $\mathcal{O}$ is a free $\mathbf{F}\smash{\bigl[Q^\vee\bigr]}$-module of finite rank,
\begin{equation}\label{FOpoldec}
\mathbf{F}[\mathcal{O}]=\bigoplus_{i=1}^N\mathbf{F}\smash{\bigl[Q^\vee\bigr]}x^{y_i}.
\end{equation}
Furthermore, the
space $\mathbf{F}[E]$ of quasi-polynomials decomposes as
\begin{equation}\label{Porbit}
\mathbf{F}[E]=\bigoplus_{\mathcal{O}\in E/W}\mathbf{F}[\mathcal{O}].
\end{equation}

For generic multiplicity functions $\mathbf{k}\in\mathcal{K}^{{\rm res}}$ and an arbitrary $W$-orbit $\mathcal{O}$, quasi-polynomial extensions of the nonsymmetric Macdonald polynomials were introduced in~\cite{SSV2}. They depend on a generic dilation parameter $q_\varphi\in\mathbf{F}^\times$ and additional $\mathcal{O}$-dependent representation parameters, and they form a basis of $\mathbf{F}[\mathcal{O}]$. For $\mathcal{O}=Q^\vee$, they are the nonsymmetric Macdonald polynomials.
The goal of this paper is to extend this result to multiplicity functions $\mathbf{k}$ in $\mathcal{K}$.

This boils down to introducing the Koornwinder $\bigl({\rm C}^\vee{\rm C}_r\bigr)$ analogs of the quasi-polynomial extensions of the type ${\rm C}_r$ nonsymmetric Macdonald polynomials from~\cite{SSV2}. They will now depend on five (four in case of $r=1$) multiplicity parameters instead of two (one in case of~${r=1}$). To stay close to the notations and results from~\cite{SSV2}, we will give a uniform treatment
of the theory for multiplicity functions $\mathbf{k}\in\mathcal{K}$ when the root system $\Phi_0$ is of arbitrary type. For type ${\rm C}_r$, the results are made more concrete in Sections~\ref{1.2}--\ref{1.4}.

\subsection{The affine Hecke algebra}\label{AHAsection}

The affine Hecke algebra $H=H(\mathbf{k})$ is the unital associative $\mathbf{F}$-algebra with generators $T_j$, $0\leq j\leq r$, and relations
\begin{enumerate}\itemsep=0pt\samepage
\item[(a)] The $(W,\{s_0,\dots,s_r\})$-braid relations for $T_0,\dots,T_r$,
\item[(b)] $(T_j-k_j)\bigl(T_j+k_j^{-1}\bigr)=0$ for $j=0,\dots,r$.
\end{enumerate}\pagebreak

\noindent
Here (a) means{\samepage
\[
T_iT_jT_i\cdots=T_jT_iT_j\cdots,\qquad 0\leq i\not=j\leq r,
\]
with on each side $m_{ij}$ terms, where $m_{ij}$ is the order of $s_is_j$ in $W$.}

A reduced expression of $g\in W$ is an expression $g=s_{j_1}\cdots s_{j_{\ell(g)}}$ of $g$ as product of simple reflections with $\ell(g)$ minimal.
The length function $\ell\colon W\rightarrow\mathbb{Z}_{\geq 0}$ satisfies $\ell(g)=\#\Pi(g)$, with
\[
\Pi(g):=\Phi^+\cap g^{-1}\Phi^-.
\]
The braid relations ensure that the element
\[
T_g:=T_{j_1}\cdots T_{j_{\ell(g)}}
\]
in $H$ does not depend on the choice of reduced expression $g=s_{j_1}\cdots s_{j_{\ell(g)}}$, and $\{T_g\}_{g\in W}$ is a~basis of $H$.

The finite Hecke algebra $H_0=H_0(\mathbf{k})$ is the subalgebra of $H$ generated by $T_1,\dots,T_r$. The defining relations of $H_0$ in terms of $T_1,\dots,T_r$
are the $(W_0,\{s_1,\dots,s_r\})$-braid relations and the quadratic relations $(T_i-k_i)\bigl(T_i+k_i^{-1}\bigr)=0$ for $i=1,\dots,r$.
For $w\in W_0\subset W$, a~reduced expression $w=s_{i_1}\cdots s_{i_\ell}$ in $W$ can be chosen with simple reflections from $W_0$, i.e., with $1\leq i_j\leq r$. Furthermore,
$\Pi(w)=\Phi_0^+\cap w^{-1}\Phi_0^-$ for $w\in W_0$, and $\{T_w\}_{w\in W_0}$ is a basis of $H_0$.

Define $\chi\colon \Phi_0\rightarrow \{\pm 1\}$ by
\begin{equation}\label{chi}
\chi(\alpha):=
\begin{cases}
\hphantom{-}1&\text{if}\ \alpha\in\Phi_0^{+},\\
-1&\text{if}\ \alpha\in\Phi_0^-.
\end{cases}
\end{equation}
In other words, $\chi=\chi_+-\chi_-$ with
\begin{equation*}
\chi_{\pm}:=\chi_{\Phi_0^{\pm}}\colon\ \Phi_0\rightarrow \{0,1\}
\end{equation*}
the indicator function of $\Phi_0^{\pm}$ in $\Phi_0$.
For $j=0,\dots,r$ and $g\in W$, we have
\begin{equation*}
\ell(s_jg)=\ell(g)+\chi\bigl(g^{-1}\alpha_j\bigr),
\end{equation*}
and hence
\begin{equation}\label{TjTw}
T_jT_g=\chi_-\bigl(g^{-1}\alpha_j\bigr)\bigl(k_j-k_j^{-1}\bigr)T_g+T_{s_jg}
\end{equation}
in the affine Hecke algebra $H$.

We now describe the Bernstein decomposition of $H$ (see~\cite{Lu} for details).
Let $H^\times$ be the group of units in $H$. There exists a unique group homomorphism
\[
Q^\vee\rightarrow H^\times,\qquad \lambda\mapsto Y^\lambda
\]
such that $Y^\lambda=T_{\tau(\lambda)}$ for \smash{$\lambda\in\overline{E_+}\cap Q^\vee$}.
The resulting algebra map
\begin{equation}\label{injY}
\mathbf{F}\smash{\bigl[Q^\vee\bigr]}\hookrightarrow H,\qquad p\mapsto p(Y),
\end{equation}
mapping $x^\lambda$ to $Y^\lambda$ for all $\lambda\in Q^\vee$, is injective. The image of the embedding is denoted by~\smash{$\mathbf{F}_Y\smash{\bigl[Q^\vee\bigr]}$}.
The multiplication map of $H$ restricts to a linear isomorphism
\[
H_0\otimes\mathbf{F}_Y\smash{\bigl[Q^\vee\bigr]}\overset{\sim}{\longrightarrow} H.
\]
The commutation relations between elements in $H_0$ and \smash{$\mathbf{F}_Y\smash{\bigl[Q^\vee\bigr]}$} are described (and determined) by the {Bernstein--Lusztig} cross relations
\begin{equation}\label{crossY}
Y^\lambda T_i-T_i Y^{s_i\lambda}=\Biggl(\frac{\widetilde{\mathbf{k}}_{\alpha_i}-\widetilde{\mathbf{k}}_{\alpha_i}^{-1}+
\Bigl(\widetilde{\mathbf{k}}_{\frac{\alpha_i}{2}}-\widetilde{\mathbf{k}}_{\frac{\alpha_i}{2}}^{-1}\Bigr)Y^{-\alpha_i^\vee}}{1-Y^{-2\alpha_i^\vee}}\Biggr)\bigl(Y^\lambda-Y^{s_i\lambda}\bigr)
\end{equation}
for $i=1,\dots,r$ and $\lambda\in Q^\vee$.

The right-hand side of~\eqref{crossY} appears to be in the quotient field $\mathbf{F}_Y\smash{\bigl(Q^\vee\bigr)}$ of $\mathbf{F}_Y\smash{\bigl[Q^\vee\bigr]}$, but it lies
in $\mathbf{F}_Y\smash{\bigl[Q^\vee\bigr]}$. Indeed, if $\alpha_i\smash{\bigl(Q^\vee\bigr)}=\mathbb{Z}$, then
\smash{$\widetilde{\mathbf{k}}_{\frac{\alpha_i}{2}}=\widetilde{\mathbf{k}}_{\alpha_i}=k_i$}, hence the right-hand side of~\eqref{crossY} reduces to
\[
\bigl(k_i-k_i^{-1}\bigr)\biggl(\frac{Y^\lambda-Y^{s_i\lambda}}{1-Y^{-\alpha_i^\vee}}\biggr)=
\bigl(k_i-k_i^{-1}\bigr)Y^\lambda\biggl(\frac{1-Y^{-\alpha_i(\lambda)\alpha^\vee}}{1-Y^{-\alpha_i^\vee}}\biggr),
\]
which lies in $\mathbf{F}_Y\smash{\bigl[Q^\vee\bigr]}$ since $\alpha_i(\lambda)\in\mathbb{Z}$. If $\alpha_i\smash{\bigl(Q^\vee\bigr)}=\mathbb{Z}_e$, then $\Phi_0$ is of type ${\rm C}_r$, $r\geq 1$,
and $i=r$, in view of our convention that $\alpha_r$ is a long root.
The right-hand side of~\eqref{crossY} then reads
\[
\bigl(k_r-k_r^{-1}+\bigl(k_0-k_0^{-1}\bigr)Y^{-\alpha_r^\vee}\bigr)Y^\lambda\biggl(\frac{1-Y^{-\alpha_r(\lambda)\alpha_r^\vee}}{1-Y^{-2\alpha_r^\vee}}\biggr),
\]
which lies in $\mathbf{F}_Y\smash{\bigl[Q^\vee\bigr]}$ since $\alpha_r(\lambda)\in \mathbb{Z}_e$.

Denote by
\[
\mathbf{F}\smash{\bigl[Q^\vee\bigr]}^{W_0}\subset\mathbf{F}\smash{\bigl[Q^\vee\bigr]}
\]
the subalgebra of $W_0$-invariant elements in $\mathbf{F}\smash{\bigl[Q^\vee\bigr]}$, and \smash{$\mathbf{F}_Y\bigl[Q^\vee\bigr]^{W_0}$} for its image in $H$ under the embedding~\eqref{injY}. Then
\[
Z(H)=\mathbf{F}_Y\smash{\bigl[Q^\vee\bigr]}^{W_0},
\]
where $Z(H)$ denotes the center of $H$.

\subsection{The double affine Hecke algebra}\label{2.5}
Fix a parameter $q_\varphi\in\mathbf{F}^\times$ and set
\[
q_\alpha:=q_\varphi^{\|\varphi\|^2/\|\alpha\|^2}\qquad \text{for}\ \alpha\in\Phi_0.
\]
It equals either $q_\varphi$, $q_\varphi^2$ or $q_\varphi^3$ (see~\cite[Section~9.4, Table~1]{Hu}), and $q_{w\alpha}=q_\alpha$ for all $\alpha\in\Phi_0$.
We like to think of $q_\varphi$ as being equal to \smash{$q^{2/\|\varphi\|^2}$} for $q$ in some field extension of $\mathbf{F}$ (this is done in~Sections~\ref{1.2}--\ref{1.4}, where we described the results explicitly for $\Phi_0$ of type ${\rm C}_r$).
To circumvent field extensions, we will introduce instead a group homomorphism $q\colon \smash{\frac{2}{\|\varphi\|^2}}\mathbb{Z}\rightarrow\mathbf{F}^\times$, $m\mapsto q^m$, with
\begin{equation*}
q^{m}:=q_\varphi^{m\|\varphi\|^2/2}\qquad \text{for}\ m\in\tfrac{2}{\|\varphi\|^2}\mathbb{Z}.
\end{equation*}
Note that $q^m\in\mathbf{F}$ is meaningful for $m\in \bigl\langle Q^\vee, Q^\vee\bigr\rangle$ since
$\bigl\langle Q^\vee, Q^\vee\bigr\rangle\subseteq\frac{2}{\|\varphi\|^2}\mathbb{Z}$.

\begin{Definition}\label{T_def}
We denote by $\mathbf{T}$ the $\mathbf{F}$-torus of rank $r$ consisting of the group homomorphisms~$Q^\vee\rightarrow\mathbf{F}^\times$.
\end{Definition}
The abelian group structure on $\mathbf{T}$ is by pointwise multiplication,
\[ (st)^\mu:=s^\mu t^\mu,\qquad s,t\in\mathbf{T},\ \mu\in Q^\vee,
\]
where $t^\mu\in\mathbf{F}^\times$ denotes the value of $t\in\mathbf{T}$ at $\mu\in Q^\vee$.
\begin{Remark}
Consider the root system $\Phi_0=\{\pm(\epsilon_i\pm\epsilon_j)\}_{1\leq i<j\leq r}\cup\{\pm 2\epsilon_i\}_{i=1}^r$ of type ${\rm C}_r$, with~$\{\epsilon_1,\dots,\epsilon_r\}$ the standard orthonormal basis of $\mathbb{R}^r$. Its co-root lattice $Q^\vee$
equals $\mathbb{Z}^r$. In the introduction, we identified the corresponding $\mathbf{F}$-torus $\mathbf{T}$ with $(\mathbf{F}^\times)^r$ by the isomorphism
\[
\mathbf{T}\overset{\sim}{\longrightarrow}(\mathbf{F}^\times)^r,\qquad t\mapsto (t^{\epsilon_1},\dots,t^{\epsilon_r}),
\]
cf.~\eqref{F_torus_CC} and~\eqref{F_torus_CCC}.
\end{Remark}

For $\lambda\in Q^\vee$, we define the torus element
\[
q^\lambda\in\mathbf{T}
\]
by $\mu\mapsto q^{\langle \lambda,\mu\rangle}$, $\mu\in Q^\vee$. Then $\mathbf{T}$ admits a left $W$-action by
\begin{align}
&(wt)^\mu:=t^{w^{-1}\mu},\nonumber\\
&(\tau(\lambda)t)^\mu:=(q^\lambda t)^\mu=q^{\langle\lambda,\mu\rangle}t^\mu\label{Waction}
\end{align}
for $t\in\mathbf{T}$, $w\in W_0$ and $\lambda,\mu\in Q^\vee$.

We will view a polynomial $p=\sum_\mu d_\mu x^\mu\in\mathbf{F}\smash{\bigl[Q^\vee\bigr]}$ as regular function on $\mathbf{T}$ by
\[
p(t):=\sum_\mu d_\mu t^\mu\qquad \text{for}\ t\in\mathbf{T}.
\]
The formula
\[
(gp)(t):=p\bigl(g^{-1}t\bigr)
\]
for $p\in\mathbf{F}\smash{\bigl[Q^\vee\bigr]}$, $g\in W$ and $t\in\mathbf{T}$ then turns $\mathbf{F}\smash{\bigl[Q^\vee\bigr]}$ into a $W$-module, with $W$ acting by algebra automorphisms.
Concretely, the action on the basis of monomials is given by
\begin{equation}\label{Wa}
w(x^\mu)=x^{w\mu},\qquad
\tau(\lambda)(x^\mu)=q^{-\langle\mu,\lambda\rangle}x^\mu
\end{equation}
for $w\in W_0$ and $\lambda,\mu\in Q^\vee$.

Note that by~\eqref{ssplit},
\begin{equation}\label{ssplitaction}
s_a(x^\mu)=q_{\overline{a}}^{-a(0)\overline{a}(\mu)}x^{s_{\overline{a}}\mu}
\end{equation}
for $a\in\Phi$ and $\mu\in Q^\vee$. In particular,
\begin{equation*}
s_0(x^\mu)=q_\varphi^{\varphi(\mu)}x^{s_\varphi\mu}.
\end{equation*}

In various computations, it is convenient to use co-roots of affine roots and incorporate $q$-powers in the exponents of the monomials. The co-root $b^\vee$ of $b\in\Phi$ is defined by
\[
b^\vee:=\biggl(\overline{b}^\vee,\frac{2b(0)}{\|\overline{b}\|^2}\biggr)\in E\times\mathbb{R}.
\]
The resulting set $\Phi^\vee:=\{b^\vee\}_{b\in\Phi}$ of co-roots is an affine root system in $E\times\mathbb{R}$~\cite{MacAff}. Note that
\begin{equation}\label{dualactionW}
(s_ab)^\vee=b^\vee-\overline{a}\bigl(\overline{b}^\vee\bigr)a^\vee\qquad \text{for}\, a,b\in\Phi
\end{equation}
in view of
\eqref{WPhiaction} and the fact that \smash{$\frac{2}{\|\beta\|^2}\beta(\alpha^\vee)=\frac{2}{\|\alpha\|^2}\alpha(\beta^\vee)$} for $\alpha,\beta\in\Phi_0$ (indeed, both sides are equal to $\langle\alpha^\vee,\beta^\vee\rangle$ by~\eqref{corootdef}).

We set
\[
x^{\widehat{\mu}}:=q^mx^\mu\in\mathbf{F}\smash{\bigl[Q^\vee\bigr]}\qquad\textup{for}\ \widehat{\mu}=(\mu,m)\in Q^\vee\times \tfrac{2}{\|\varphi\|^2}\mathbb{Z}
\]
and we write $t^{\widehat{\mu}}$ for the evaluation of $x^{\widehat{\mu}}\in\mathbf{F}\smash{\bigl[Q^\vee\bigr]}$ at $t\in\mathbf{T}$ \big(so in particular, $t^{\widehat{\mu}}=q^mt^\mu$\big).
Note that \smash{$\Phi^\vee\subset Q^\vee\times \frac{2}{\|\varphi\|^2}\mathbb{Z}$}, hence \smash{$x^{b^\vee}$} makes sense for all $b\in\Phi$. Concretely,
\begin{equation}\label{arootreduction}
x^{b^\vee}:=q_{\overline{b}}^{b(0)}x^{\overline{b}^\vee}\in\mathbf{F}\smash{\bigl[Q^\vee\bigr]},
\end{equation}
which reduces to the monomial \smash{$x^{\beta^\vee}$} in case $b=(\beta,0)$.

Formula~\eqref{ssplitaction} can then be rewritten as
\begin{equation}\label{ssplitaction2}
s_a(x^\mu)=x^{\mu-\overline{a}(\mu)a^\vee}.
\end{equation}
We furthermore have
\begin{equation}\label{awcomp}
g(x^{b^\vee})=x^{(gb)^\vee}\qquad \textup{for}\ g\in W\ \text{and}\ b\in\Phi.
\end{equation}
Indeed, it suffices to check~\eqref{awcomp} for $w=s_a$, $a\in\Phi$. By~\eqref{arootreduction} and~\eqref{ssplitaction2}, we have
\[
s_a\bigl(x^{b^\vee}\bigr)=x^{b^\vee-\overline{a}(\overline{b}^\vee)a^\vee}
\]
which equals \smash{$x^{(s_ab)^\vee}$} by~\eqref{dualactionW}.

\begin{Definition}[\cite{ChKZ2,Sa}]
The double affine Hecke algebra $\mathbb{H}=\mathbb{H}(\mathbf{k};q_\varphi)$ is the unital associative $\mathbf{F}$-algebra generated by
$T_j$, $0\leq j\leq r$, and $x^\mu$, $\mu\in Q^\vee$, subject to the following relations:
\begin{enumerate}\itemsep=0pt
\item[(a)] the $(W,\{s_0,\dots,s_r\})$-braid relations for $T_0,\dots,T_r$,
\item[(b)] the quadratic relations $(T_j-k_j)\bigl(T_j+k_j^{-1}\bigr)=0$ for $j=0,\dots,r$,
\item[(c)] $x^\mu x^\nu=x^{\mu+\nu}$, $\mu,\nu\in Q^\vee$, and $x^0$ is the unit element of $\mathbb{H}$,
\item[(d)] the cross relations
\begin{equation}\label{crossX}
T_jx^\mu-s_j(x^\mu)T_j=\Biggl(\frac{k_j-k_j^{-1}+\bigl(u_j-u_j^{-1}\bigr)x^{\alpha_j^\vee}}{1-x^{2\alpha_j^\vee}}\Biggr)(x^\mu-s_j(x^\mu))
\end{equation}
for $j=0,\dots,r$ and $\mu\in Q^\vee$.
\end{enumerate}
\end{Definition}
The Poincar{\'e}--Birkhoff--Witt (PBW) theorem for $\mathbb{H}$ states that the canonical algebra maps $H\rightarrow \mathbb{H}$ and $\mathbf{F}\smash{\bigl[Q^\vee\bigr]}\rightarrow
\mathbb{H}$
are embeddings, and that
the multiplication map of $\mathbb{H}$ restricts to a~linear isomorphism
\[
\mathbf{F}\smash{\bigl[Q^\vee\bigr]}\otimes H\overset{\sim}{\longrightarrow} \mathbb{H}.
\]
By the Bernstein presentation of the affine Hecke algebra (see Section~\ref{AHAsection}), the subalgebra
\[
H^X\subset\mathbb{H}
\]
generated by $\mathbf{F}\smash{\bigl[Q^\vee\bigr]}$ and $H_0$ is isomorphic to the affine Hecke algebra $\widetilde{H}:=H(\widetilde{\mathbf{k}})$.

The double affine Hecke algebra with dual multiplicity parameters will be denoted by
\[
\widetilde{\mathbb{H}}:=\mathbb{H}\bigl(\widetilde{\mathbf{k}},q_\varphi\bigr).
\]
To keep the notations manageable, we will use the same notations $T_j$, $T_g$, $Y^\mu$, $x^\mu$ in both $\mathbb{H}$ and~$\widetilde{\mathbb{H}}$.

The duality anti-isomorphism~\cite{ChBook,Sa} is the unique anti-algebra isomorphism
\[
\delta=\delta_{\mathbf{k}}\colon\ \mathbb{H}\overset{\sim}{\longrightarrow}\widetilde{\mathbb{H}}
\]
satisfying
\[
\delta(T_i)=T_i,\qquad \delta\bigl(Y^\lambda\bigr)=x^{-\lambda},\qquad \delta(x^\mu)=Y^{-\mu}
\]
for $i=1,\dots,r$ and $\lambda,\mu\in Q^\vee$. Its inverse is $\widetilde{\delta}:=\delta_{\widetilde{\mathbf{k}}}$.

\section{Quasi-polynomial representations}\label{qpsect}

\subsection[The quasi-polynomial representation of H\^{}X]{The quasi-polynomial representation of $\boldsymbol{H^X}$}\label{qpHX}

In this subsection, we introduce the quasi-polynomial representation of the dual affine Hecke algebra $H^X$.
In case $\mathbf{k}\in\mathcal{K}^{{\rm res}}$,
the representation we obtain was derived before in~\cite{SSV2, SSV1}.

Let
\[
\lfloor\cdot\rfloor\colon \ \mathbb{R}\rightarrow \mathbb{Z}
\]
be the floor function, so $\lfloor s\rfloor$ is the largest integer $\leq s$. We denote by
\[
\lfloor\cdot\rfloor_e\colon \ \mathbb{R}\rightarrow\mathbb{Z}_e,\qquad
\lfloor\cdot\rfloor_o\colon \ \mathbb{R}\rightarrow\mathbb{Z}_o
\]
the functions which map $s\in\mathbb{R}$ to
the smallest even and odd integer $\leq s$, respectively. Note~that
\[
\lfloor s\rfloor_e:=2\lfloor s/2\rfloor,\qquad
\lfloor s\rfloor_o=\lfloor s+1\rfloor_e-1
\]
for $s\in\mathbb{R}$.

\begin{Definition}\label{truncation}
For $a\in\Phi$, define the even and odd truncated divided difference operator $\nabla_a=\nabla_a(\mathbf{k})\in\textup{End}(\mathbf{F}[E])$ by
\begin{align}
&\nabla_a^e(x^y):=\Biggl(\frac{1-x^{-\lfloor\overline{a}(y)\rfloor_ea^\vee}}{1-x^{2a^\vee}}\Biggr)x^y,\nonumber\\
&\nabla_a^o(x^y):=\Biggl(\frac{x^{a^\vee}-x^{-\lfloor\overline{a}(y)\rfloor_oa^\vee}}{1-x^{2a^\vee}}\Biggr)x^y\label{Ai}
\end{align}
for $y\in E$. We furthermore write $\nabla_j^e:=\nabla_{\alpha_j}^e$ and $\nabla_j^o:=\nabla_{\alpha_j}^o$ for $j=0,\dots,r$.
\end{Definition}

Truncation in Definition~\ref{truncation} refers to the fact that the real numbers $\overline{a}(y)$ in formula~\eqref{Ai} are truncated using the even and odd floor operations.
These truncations are necessary to turn the two quotients in~\eqref{Ai} into well defined elements in $\mathbf{F}\smash{\bigl[Q^\vee\bigr]}$. Note that $\nabla_a^e$ and $\nabla_a^o$ depend on~$q_\varphi$ when $a\in\Phi\setminus\Phi_0$, in view of~\eqref{arootreduction}. In this subsection, we only need the truncated divided difference operators for $a\in\Phi_0$, and there will be no dependence on $q_\varphi$.

We write
\begin{equation}\label{Aisum}
\nabla_a:=\nabla_a^e+\nabla_a^o
\end{equation}
for the sum of the even and odd truncated divided difference operator, and $\nabla_j:=\nabla_{\alpha_j}$ for $j=0,\dots,r$. Then
\begin{equation}\label{Aires}
\nabla_a(x^y)=\Biggl(\frac{1-x^{-\lfloor\overline{a}(y)\rfloor a^\vee}}{1-x^{a^\vee}}\Biggr)x^y
\end{equation}
for $y\in E$ since
\begin{equation*}
\{\lfloor s\rfloor_e,\lfloor s\rfloor_o\}=\{\lfloor s\rfloor,\lfloor s\rfloor-1\}
\end{equation*}
as unordered $2$-sets for any $s\in\mathbb{R}$.
The truncated difference operator $\nabla_a\in\textup{End}(\mathbf{F}[E])$ was introduced before in~\cite[Section~4.2]{SSV2}.

The link of the various truncated divided difference operators to the usual divided difference operator is as follows (the first part of the lemma was observed before in~\cite[Lemma~4.4]{SSV2}).
\begin{Lemma}\label{polredNabla}
\hfill
\begin{enumerate}\itemsep=0pt
\item[$(1)$] If $a\in\Phi$, then
\[
\nabla_a(x^\mu)=\frac{x^\mu-s_a(x^{\mu})}{1-x^{a^\vee}}
\]
for $\mu\in Q^\vee$.
\item[$(2)$] If $\Phi_0$ is of type ${\rm C}_r$, $r\geq 1$, and if $a\in\Phi$ is an affine root such that $\overline{a}\in\Phi_0^\ell$ \textup{(}in other words,
$a\in W\alpha_0\sqcup W\alpha_r$\textup{)},
then
\[
\nabla_a^e(x^\mu)=\frac{x^\mu-s_a(x^{\mu})}{1-x^{2a^\vee}}=x^{-a^\vee}\nabla_a^o(x^\mu)
\]
for $\mu\in Q^\vee$.
\end{enumerate}
\end{Lemma}
\begin{proof}
(1) This is immediate from~\eqref{ssplitaction2} and the fact that $\overline{a}\smash{\bigl(Q^\vee\bigr)}\subseteq\mathbb{Z}$.\\
(2)
Under these assumptions, we have $\overline{a}\smash{\bigl(Q^\vee\bigr)}=\mathbb{Z}_e$, hence
\[
\lfloor \overline{a}(\mu)\rfloor_e=\overline{a}(\mu)=\lfloor \overline{a}(\mu)\rfloor_o+1
\]
and the result follows again from~\eqref{ssplitaction2}.
\end{proof}

We will often make use of the indicator functions $\chi_{\mathbb{Z}_e}, \chi_{\mathbb{Z}_o}\colon \mathbb{R}\rightarrow\{0,1\}$ of the
subsets of even and odd integers, respectively. We will denote them by $\chi_e$ and $\chi_o$, respectively.

\begin{Theorem}\label{theoremH}
The formulas
\begin{align}
&\pi(x^\mu)x^y:=x^{y+\mu},\nonumber\\
&\pi(T_i)x^y:=k_i^{\chi_e(\alpha_i(y))}u_i^{\chi_o(\alpha_i(y))}x^{s_iy}+\bigl(k_i-k_i^{-1}\bigr)\nabla_i^e(x^y)+\bigl(u_i-u_i^{-1}\bigr)\nabla_i^o(x^y)\label{qprepHX}
\end{align}
for $\mu\in Q^\vee$, $i\in\{1,\dots,r\}$ and $y\in E$
define a representation $\pi\colon H^X\rightarrow\textup{End}(\mathbf{F}[E])$.
\end{Theorem}
\begin{proof}
Formula~\eqref{qprepHX} uniquely defines linear operators $\pi(x^\mu)$ and $\pi(T_i)$ on $\mathbf{F}[E]$, which in turn restrict to linear operators
\begin{equation}\label{piO}
\pi^{\mathcal{O}}(x^\mu):=\pi(x^\mu)\vert_{\mathbf{F}[\mathcal{O}]},\qquad \pi^{\mathcal{O}}(T_i):=\pi(T_i)\vert_{\mathbf{F}[\mathcal{O}]}
\end{equation}
on $\mathbf{F}[\mathcal{O}]$ for every $W$-orbit $\mathcal{O}$ in $E$. In view of~\eqref{Porbit} it thus suffices to show that the linear operators~\eqref{piO} define a representation $\pi^{\mathcal{O}}=\pi(\cdot)\vert_{\mathbf{F}[\mathcal{O}]}\colon H^X\rightarrow\textup{End}(\mathbf{F}[\mathcal{O}])$.

Consider $H^X$ as left regular $H^X$-module. By~\eqref{TjTw} and~\eqref{crossX} the $H^X$-action can be written down explicitly relative to the basis
\[
\bigl\{x^\lambda T_w\mid \lambda\in Q^\vee,\, w\in W_0\bigr\}
\]
of $H^X$. The resulting formulas are
\begin{gather}
x^\mu x^\lambda T_w=x^{\lambda+\mu}T_w,\nonumber\\
T_ix^\lambda T_w=x^{s_i\lambda}T_{s_iw}+\bigl(k_i-k_i^{-1}\bigr)\Biggl(\frac{1-x^{(2\chi_-\bigl(w^{-1}\alpha_i\bigr)-\alpha_i(\lambda))\alpha_i^\vee}}{1-x^{2\alpha_i^\vee}}\Biggr)x^\lambda T_w\nonumber\\
\hphantom{T_ix^\lambda T_w=}{}
+(u_i-u_i^{-1})\Biggl(\frac{x^{\alpha_i^\vee}-x^{(1-\alpha_i(\lambda))\alpha_i^\vee}}{1-x^{2\alpha_i^\vee}}\Biggr)x^{\lambda}T_w
\label{regularHX}
\end{gather}
for $\lambda,\mu\in Q^\vee$, $w\in W_0$ and $i=1,\dots,r$.

Let $\kappa_w^\mathcal{O}\in\mathbf{F}^\times$, $w\in W_0$, be a collection of nonzero scalars.
Consider the surjective linear map
\begin{equation}\label{psicmap}
\psi^\mathcal{O}\colon\ H^X\twoheadrightarrow \mathbf{F}[\mathcal{O}],\qquad \psi^\mathcal{O}\bigl(x^\lambda T_w\bigr):=\kappa_w^\mathcal{O}x^{\lambda+wc^{\mathcal{O}}},\quad \lambda\in Q^\vee,\ w\in W_0.
\end{equation}
We will fine tune the scalars $\kappa_w^\mathcal{O}\in\mathbf{F}^\times$, $w\in W_0$, in such a way that the kernel of $\psi^\mathcal{O}$
is a left ideal in $H^X$. This will allow us to push the left regular $H^X$-action through $\psi^\mathcal{O}$, giving rise to a
$H^X$-action on $\mathbf{F}[\mathcal{O}]$. We then show that the resulting action of $x^\mu$ and $T_i$ on $\mathbf{F}[\mathcal{O}]$ is by the linear operators $\pi^\mathcal{O}(x^\mu)$ and $\pi^\mathcal{O}(T_i)$, which completes the proof of the theorem.

For any choice of scalars $\kappa_w^\mathcal{O}\in\mathbf{F}^\times$, $w\in W_0$, the kernel of $\psi^\mathcal{O}$ is invariant under left multiplication by $\mathbf{F}\smash{\bigl[Q^\vee\bigr]}$.
Note that the kernel $\psi^\mathcal{O}$ is invariant under left multiplication by $T_i$ if there exists a linear operator $\mathcal{D}_i^\mathcal{O}\in\textup{End}(\mathbf{F}[\mathcal{O}])$ such that
\begin{equation*}
\mathcal{D}_i^\mathcal{O}\bigl(\psi^\mathcal{O}(h)\bigr)=\psi^\mathcal{O}(T_ih)
\end{equation*}
for all $h\in H^X$.
By~\eqref{regularHX}, we have
\begin{align}
\psi^\mathcal{O}\bigl(T_ix^\lambda T_w\bigr)={}&\kappa_{s_iw}^\mathcal{O}s_i\bigl(x^{\lambda+wc^{\mathcal{O}}}\bigr)
+ \kappa_w^\mathcal{O}\bigl(k_i-k_i^{-1}\bigr)\Biggl(\frac{1-x^{(2\chi_-\bigl(w^{-1}\alpha_i\bigr)-\alpha_i(\lambda))\alpha_i^\vee}}{1-x^{2\alpha_i^\vee}}\Biggr)x^{\lambda+wc^{\mathcal{O}}}\nonumber\\
&{}+
\kappa_w^\mathcal{O}\bigl(u_i-u_i^{-1}\bigr)\Biggl(\frac{x^{\alpha_i^\vee}-x^{(1-\alpha_i(\lambda))\alpha_i^\vee}}{1-x^{2\alpha_i^\vee}}\Biggr)x^{\lambda+wc^{\mathcal{O}}},
\label{spoint}
\end{align}
so we look for conditions on the scalars $\kappa_w^\mathcal{O}$ such that~\eqref{spoint} can be expressed as
a linear operator $\mathcal{D}_i^\mathcal{O}\in\textup{End}(\mathbf{F}[\mathcal{O}])$ acting on $\psi^\mathcal{O}\bigl(x^\lambda T_w\bigr)=\kappa_w^\mathcal{O}x^{\lambda+wc^\mathcal{O}}$ for all $\lambda\in Q^\vee$ and $w\in W_0$.

To simplify notations, we write
\[
y=\lambda+wc^\mathcal{O}
\]
with $w\in W_0$ and $\lambda\in Q^\vee$.
We first prove that
\begin{gather}
\frac{1}{\kappa_w^\mathcal{O}} \psi^\mathcal{O}\bigl(T_ix^\lambda T_w\bigr)
=\biggl(\frac{\kappa_{s_iw}^\mathcal{O}}{\kappa_w^\mathcal{O}}
+\bigl(k_i-k_i^{-1}\bigr)\chi_-\bigl(w^{-1}\alpha_i\bigr)\chi_e\bigl(\alpha_i\bigl(wc^\mathcal{O}\bigr)\bigr)\nonumber\\
\hphantom{\frac{1}{\kappa_w^\mathcal{O}} \psi^\mathcal{O}\bigl(T_ix^\lambda T_w\bigr)=\biggl(}{}
 +\bigl(u_i-u_i^{-1}\bigr)\chi_+\bigl(w^{-1}\alpha_i\bigr)\chi_o(\alpha_i\bigl(wc^\mathcal{O}\bigr))\biggr)x^{s_iy}\nonumber\\
 \hphantom{\frac{1}{\kappa_w^\mathcal{O}} \psi^\mathcal{O}\bigl(T_ix^\lambda T_w\bigr)=}{}
 +\bigl(k_i-k_i^{-1}\bigr)\nabla_i^e(x^y)+\bigl(u_i-u_i^{-1}\bigr)\nabla_i^o(x^y)
\label{formula}
\end{gather}
by rewriting the two quotients in~\eqref{spoint} in terms of the odd and even truncated difference operators.

We first consider the proof of~\eqref{formula} when $\alpha_i\smash{\bigl(Q^\vee\bigr)}=\mathbb{Z}_e$, i.e., when $\Phi_0$ is of type ${\rm C}_r$, $r\geq 1$, and $i=r$.
Then~\eqref{W0C} and $\alpha_r\smash{\bigl(Q^\vee\bigr)}=\mathbb{Z}_e$ imply that
\begin{equation*}
2\chi_-\bigl(w^{-1}\alpha_r\bigr)-\alpha_r(\lambda)=
\begin{cases}
2-\lfloor\alpha_r(y)\rfloor_e &\hbox{if }  w^{-1}\alpha_r\in\Phi_0^-\ \text{and}\ \alpha_r\bigl(wc^\mathcal{O}\bigr)=0,\\
-\lfloor\alpha_r(y)\rfloor_e &\textup{otherwise}.
\end{cases}
\end{equation*}
By~\eqref{W0C},
this can be reformulated as
\begin{equation*}
2\chi_-\bigl(w^{-1}\alpha_r\bigr)-\alpha_r(\lambda)=
\begin{cases}
2-\lfloor\alpha_r(y)\rfloor_e &\text{if}\ \chi_-\bigl(w^{-1}\alpha_r\bigr)\chi_e\bigl(\alpha_r\bigl(wc^\mathcal{O}\bigr)\bigr)=1,\\
-\lfloor\alpha_r(y)\rfloor_e &\text{if}\ \chi_-\bigl(w^{-1}\alpha_r\bigr)\chi_e\bigl(\alpha_r\bigl(wc^\mathcal{O}\bigr)\bigr)=0.
\end{cases}
\end{equation*}
This allows us to rewrite the second line of~\eqref{spoint} for $\Phi_0$ of type ${\rm C}_r$ and $i=r$ in terms of the even truncated divided difference operator,
\begin{equation}\label{claim1}
\Biggl(\frac{1-x^{(2\chi_-(w^{-1}\alpha_r)-\alpha_r(\lambda))\alpha_r^\vee}}{1-x^{2\alpha_r^\vee}}\Biggr)x^{y}=
\nabla_r^e(x^y)+
\chi_-\bigl(w^{-1}\alpha_r\bigr)\chi_e\bigl(\alpha_r\bigl(wc^\mathcal{O}\bigr)\bigr)x^{s_ry}
\end{equation}
\big(in case $\chi_-\bigl(w^{-1}\alpha_r\bigr)\chi_e\bigl(\alpha_r\bigl(wc^\mathcal{O}\bigr)\bigr)=1$ use the fact that $\lfloor\alpha_r(y)\rfloor_e=\alpha_r(y)$\big).
To rewrite the third line of~\eqref{spoint} for $\Phi_0$ of type ${\rm C}_r$ and $i=r$, note that
\begin{equation*}
1-\alpha_r(\lambda)=
\begin{cases}
2-\lfloor\alpha_r(y)\rfloor_o &\text{if}\ \alpha_r\bigl(wc^\mathcal{O}\bigr)=1,\\
-\lfloor\alpha_r(y)\rfloor_o &\text{otherwise},
\end{cases}
\end{equation*}
which can be reformulated as
\begin{equation*}
1-\alpha_r(\lambda)=
\begin{cases}
2-\lfloor\alpha_r(y)\rfloor_o &\text{if}\ \chi_+\bigl(w^{-1}\alpha_r\bigr)\chi_o\bigl(\alpha_r\bigl(wc^\mathcal{O}\bigr)\bigr)=1,\\
-\lfloor\alpha_r(y)\rfloor_o &\text{if}\ \chi_+\bigl(w^{-1}\alpha_r\bigr)\chi_o\bigl(\alpha_r\bigl(wc^\mathcal{O}\bigr)\bigr)=0
\end{cases}
\end{equation*}
in view of~\eqref{W0C}. It follows that
\begin{equation}\label{claim2}
\Biggl(\frac{x^{\alpha_r^\vee}-x^{(1-\alpha_r(\lambda))\alpha_r^\vee}}{1-x^{2\alpha_r^\vee}}\Biggr)x^y=
\nabla_r^o(x^y)+
\chi_+\bigl(w^{-1}\alpha_r\bigr)\chi_o\bigl(\alpha_r\bigl(wc^\mathcal{O}\bigr)\bigr)
x^{s_ry}
\end{equation}
(in case $\chi_+(w^{-1}\alpha_r)\chi_o(\alpha_r\bigl(wc^\mathcal{O}\bigr))=1$ use the fact that $\lfloor\alpha_r(y)\rfloor_o=\alpha_r(y)$).
Substituting~\eqref{claim1} and~\eqref{claim2} into formula~\eqref{spoint} for $i=r$, we obtain~\eqref{formula} for $\Phi_0$ of type ${\rm C}_r$ and for $i=r$.

We now prove~\eqref{formula} when $\alpha_i\smash{\bigl(Q^\vee\bigr)}=\mathbb{Z}$. Then $k_i=u_i$ and using that
\[
\bigl\{2\chi_-\bigl(w^{-1}\alpha_i\bigr)-\alpha_i(\lambda), 1-\alpha_i(\lambda)\bigr\}=\bigl\{\chi_-\bigl(w^{-1}\alpha_i\bigr)-\alpha_i(\lambda),1+\chi_-\bigl(w^{-1}\alpha_i\bigr)-\alpha_i(\lambda)\bigr\}
\]
as unordered $2$-sets,
formula~\eqref{spoint} simplifies to
\begin{equation}\label{spointred}
\psi^\mathcal{O}(T_ix^\lambda T_w)=\kappa_{s_iw}^\mathcal{O}x^{s_iy}
+\kappa_w^\mathcal{O}\bigl(k_i-k_i^{-1}\bigr)\Biggl(\frac{1-x^{(\chi_-\bigl(w^{-1}\alpha_i\bigr)-\alpha_i(\lambda))\alpha_i^\vee}}{1-x^{\alpha_i^\vee}}\Biggr)x^{y}.
\end{equation}
Using~\eqref{W0C}, we have for $w^{-1}\alpha_i\in\Phi_0^-$ that
\begin{equation*}
\chi_-\bigl(w^{-1}\alpha_i\bigr)-\alpha_i(\lambda)=
\begin{cases}
-\lfloor\alpha_i(y)\rfloor &\text{if}\ \chi_e\bigl(\alpha_i\bigl(wc^\mathcal{O}\bigr)\bigr)=0,\\
1-\lfloor\alpha_i(y)\rfloor &\text{if}\ \chi_e\bigl(\alpha_i\bigl(wc^\mathcal{O}\bigr)\bigr)=1
\end{cases}
\end{equation*}
and for $w^{-1}\alpha_i\in\Phi_0^+$ that
\begin{equation*}
\chi_-\bigl(w^{-1}\alpha_i\bigr)-\alpha_i(\lambda)=
\begin{cases}
-\lfloor\alpha_i(y)\rfloor &\text{if}\ \chi_o\bigl(\alpha_i\bigl(wc^\mathcal{O}\bigr)\bigr)=0,\\
1-\lfloor\alpha_i(y)\rfloor &\text{if}\ \chi_0\bigl(\alpha_i\bigl(wc^\mathcal{O}\bigr)\bigr)=1.
\end{cases}
\end{equation*}
This leads to the formula
\begin{gather*}
\Biggl(\frac{1-x^{(\chi_-\bigl(w^{-1}\alpha_i\bigr)-\alpha_i(\lambda))\alpha_i^\vee}}{1-x^{\alpha_i^\vee}}\Biggr)x^{y}\\
\qquad{}{}=
\nabla_i(x^y)+\bigl(\chi_-\bigl(w^{-1}\alpha_i\bigr)\chi_e\bigl(\alpha_i\bigl(wc^\mathcal{O}\bigr)\bigr)+\chi_+\bigl(w^{-1}\alpha_i\bigr)\chi_o\bigl(\alpha_i\bigl(wc^\mathcal{O}\bigr)\bigr)\bigr)x^{s_iy}.
\end{gather*}
Substituting into~\eqref{spointred} and using the formulas~\eqref{Aisum} and~\eqref{Aires}, we now also obtain~\eqref{formula} in case $\alpha_i\smash{\bigl(Q^\vee\bigr)}=\mathbb{Z}$.

The next step is choosing the normalization factors $\kappa_w^\mathcal{O}$, $w\in W_0$, in such a way that the coefficient of $x^{s_i\lambda}$ in~\eqref{formula}
only depends on $y=\lambda+wc^\mathcal{O}$. We will show that this is the case for the scalars
\begin{equation}\label{kappawy}
\kappa_w^\mathcal{O}:=\prod_{\alpha\in\Pi(w)}\mathbf{k}_\alpha^{\chi_e(\alpha(c^\mathcal{O}))}\mathbf{k}_{\alpha/2}^{-\chi_o(\alpha(c^\mathcal{O}))}
\end{equation}
(a small computation using~\eqref{W0C} shows that $\kappa_w^\mathcal{O}$ reduces for $\mathbf{k}\in\mathcal{K}^{{\rm res}}$ to the normalization factor~\cite[(4.11)]{SSV2}).
Since $\Pi(w)=\Phi_0^+\cap w^{-1}\Phi_0^-$ for $w\in W_0$, we have
\begin{equation*}
\Pi(s_iw)=
\begin{cases}
\Pi(w)\cup\big\{w^{-1}\alpha_i\big\} &\text{if}\ \chi_+\bigl(w^{-1}\alpha_i\bigr)=1,\\
\Pi(w)\setminus\big\{-w^{-1}\alpha_i\big\} &\text{if} \
\chi_-\bigl(w^{-1}\alpha_i\bigr)=1.
\end{cases}
\end{equation*}
Combined with the fact that $\chi_e,\chi_o\colon \mathbb{R}\rightarrow\{0,1\}$ are even functions, we obtain
\begin{equation*}
\frac{\kappa_{s_iw}^\mathcal{O}}{\kappa_w^\mathcal{O}}=
\begin{cases}
k_i^{\chi_e(\alpha_i (wc^\mathcal{O}))}u_i^{-\chi_o(\alpha_i (wc^\mathcal{O} ))} &\text{if}\ \chi_+\bigl(w^{-1}\alpha_i\bigr)=1,\\
k_i^{-\chi_e(\alpha_i (wc^\mathcal{O}))}u_i^{\chi_o(\alpha_i (wc^\mathcal{O} ))} &\text{if}\ \chi_-\bigl(w^{-1}\alpha_i\bigr)=1,
\end{cases}
\end{equation*}
and hence
\begin{gather*}
\frac{\kappa_{s_iw}^\mathcal{O}}{\kappa_w^\mathcal{O}}
+\bigl(k_i-k_i^{-1}\bigr)\chi_-\bigl(w^{-1}\alpha_i\bigr)\chi_e\bigl(\alpha_i\bigl(wc^\mathcal{O}\bigr)\bigr)+\bigl(u_i-u_i^{-1}\bigr)\chi_+\bigl(w^{-1}\alpha_i\bigr)\chi_o\bigl(\alpha_i\bigl(wc^\mathcal{O}\bigr)\bigr)\\
\qquad =k_i^{\chi_e(\alpha_i\bigl(wc^\mathcal{O}\bigr))}u_i^{\chi_o(\alpha_i\bigl(wc^\mathcal{O}\bigr))}
\end{gather*}
for $w\in W_0$ and $i=1,\dots,r$. So formula~\eqref{formula} reduces
for the specific choice~\eqref{kappawy} of $\kappa_w^\mathcal{O}$ to
\begin{align}
\frac{1}{\kappa_w^\mathcal{O}}\psi^\mathcal{O}\bigl(T_ix^\lambda T_w\bigr)
={}&k_i^{\chi_e(\alpha_i (wc^\mathcal{O} ))}u_i^{\chi_o(\alpha_i (wc^\mathcal{O} ))}x^{s_iy}\nonumber\\
&{}{+}\ \bigl(k_i-k_i^{-1}\bigr)\nabla_i^e(x^y)+\bigl(u_i-u_i^{-1}\bigr)\nabla_i^o(x^y).\label{formula2}
\end{align}
Now note that the map
\begin{equation}\label{functionQ}
E\rightarrow\mathbf{F}^\times,\qquad y\mapsto k_i^{\chi_e(\alpha_i(y))}u_i^{\chi_o(\alpha_i(y))}
\end{equation}
is $\tau\smash{\bigl(Q^\vee\bigr)}$-invariant. This is trivial when $\Phi_0$ is of type ${\rm C}_r$, $r\geq 1$, and $i=r$, since in this case $\alpha_r\smash{\bigl(Q^\vee\bigr)}=\mathbb{Z}_e$. In all other cases, $\alpha_i\smash{\bigl(Q^\vee\bigr)}=\mathbb{Z}$ hence $k_i=u_i$, in which case it follows from the fact that the function~\eqref{functionQ} reduces to $y\mapsto k_i^{\chi_{\mathbb{Z}}(\alpha_i(y))}$.
In~\eqref{formula2}, we may thus replace the coefficient
\[
k_i^{\chi_e(\alpha_i (wc^\mathcal{O} ))}u_i^{\chi_o(\alpha_i (wc^\mathcal{O} ))}
\] of $x^{s_iy}$ by
\smash{$k_i^{\chi_e(\alpha_i(y))}u_i^{\chi_o(\alpha_i(y))}$}.

In conclusion, for $\kappa_w^\mathcal{O}$ given by~\eqref{kappawy} the associated linear map $\psi^\mathcal{O}$~\eqref{psicmap} satisfies
\[
\psi^\mathcal{O}(T_ih)=\mathcal{D}_i^\mathcal{O}\bigl(\psi^\mathcal{O}(h)\bigr)\qquad \forall h\in H^X
\]
with $\mathcal{D}_i^\mathcal{O}\in\textup{End}(\mathbf{F}[\mathcal{O}])$ defined by
\begin{equation*}
\mathcal{D}_i^\mathcal{O}(x^y):=k_i^{\chi_e(\alpha_i(y))}u_i^{\chi_o(\alpha_i(y))}x^{s_iy}+\bigl(k_i-k_i^{-1}\bigr)\nabla_i^e(x^y)+\bigl(u_i-u_i^{-1}\bigr)\nabla_i^o(x^y)
\end{equation*}
for $y\in\mathcal{O}$.
The kernel of $\psi^\mathcal{O}\colon H^X\twoheadrightarrow\mathbf{F}[\mathcal{O}]$ thus is a left-ideal, and $\mathbf{F}[\mathcal{O}]$ inherits a $H^X$-action from $H^X/\textup{ker}\bigl(\psi^\mathcal{O}\bigr)$ with $T_i$ acting by $\mathcal{D}_i^\mathcal{O}=\pi^\mathcal{O}(T_i)$ and $x^\mu$ acting by
$\pi^\mathcal{O}(x^\mu)$. This completes the proof of the theorem.
\end{proof}

\begin{Remark}\label{H_remark}
The proof of Theorem~\ref{theoremH} involves a particular choice of normalisation factors $\kappa_w^{\mathcal{O}}\in\mathbf{F}^\times$ ($w\in W_0$). Any choice of $\kappa_w^\mathcal{O}$, $w\in W_0$, such that the coefficient of $x^{s_iy}$, ${y=\lambda+wc^{\mathcal{O}}}$, in~\eqref{formula} only depends on the coset $w\{v\in W_0 \mid vc^\mathcal{O}=c^\mathcal{O}\}$ for all ${y=\lambda+wc^\mathcal{O}\in\mathcal{O}}$ and ${i\in\{1,\dots,r\}}$
will lead to an explicit $H^X$-representation on $\mathbf{F}[\mathcal{O}]$ involving truncated Demazure--Lusztig type operators. The present choice~\eqref{kappawy} corresponds to a natural class of parabolically induced $H^X$-modules,
see Section~\ref{section_pi_mod} for details.
\end{Remark}
\begin{Corollary}\label{corpsi}\hfill
\begin{enumerate}\itemsep=0pt
\item[$(1)$] $\mathbf{F}[E]=\bigoplus_{\mathcal{O}\in E/W}\mathbf{F}[\mathcal{O}]$ is a decomposition of $\mathbf{F}[E]$ in $H^X$-submodules.
\item[$(2)$] Let $\mathcal{O}\in E/W$. Then
\begin{equation*}
\pi\bigl(x^\lambda T_w\bigr)x^{c^\mathcal{O}}=\kappa_w^\mathcal{O}x^{\lambda+wc^\mathcal{O}}
\end{equation*}
for $\lambda\in Q^\vee$ and $w\in W_0$,
with $\kappa_w^\mathcal{O}$ defined by~\eqref{kappawy}.
\end{enumerate}
\end{Corollary}

 \begin{proof}
 (1) This was remarked in the first paragraph of the proof of Theorem~\ref{theoremH}.

(2) By the last paragraph of the proof of Theorem~\ref{theoremH}, the epimorphism $\psi^\mathcal{O}\colon H^X\twoheadrightarrow\mathbf{F}[\mathcal{O}]$, mapping
$x^\lambda T_w$ to $\kappa_w^\mathcal{O}x^{\lambda+wc^\mathcal{O}}$ for $\lambda\in Q^\vee$ and $w\in W_0$, is $H^X$-linear for the special choice~\eqref{kappawy} of the scalars $\kappa_w^\mathcal{O}$. Hence
\[
\pi\bigl(x^\lambda T_w\bigr)x^{c^\mathcal{O}}=\pi\bigl(x^\lambda T_w\bigr)\psi^\mathcal{O}(1)=\psi^\mathcal{O}\bigl(x^\lambda T_w\bigr)=
\kappa_w^\mathcal{O}x^{\lambda+wc^\mathcal{O}}
\]
for $\lambda\in Q^\vee$ and $w\in W_0$.
 \end{proof}

 For the upgrade of Theorem~\ref{theoremH} to the double affine Hecke algebra $\mathbb{H}$ (see Section~\ref{wtsection}), it is useful to introduce the notation
 \begin{equation*}
 \pi^\mathcal{O}(\cdot):=\pi(\cdot)\vert_{\mathbf{F}[\mathcal{O}]}
 \end{equation*}
 for the representation map of the $H^X$-submodule $\mathbf{F}[\mathcal{O}]$ (as we in fact already have done in the proof of Theorem~\ref{theoremH}). The reason for this is that the extension of $\pi^{\mathcal{O}}\colon H^X\rightarrow\textup{End}(\mathbf{F}[\mathcal{O}])$ to a~$\mathbb{H}$-representation on $\mathbf{F}[\mathcal{O}]$ involves additional $\mathcal{O}$-dependent representation parameters.

\begin{Remark}\quad
\begin{enumerate}\itemsep=0pt
\item
By~\eqref{Aisum}, we have
\begin{equation*}
\pi(T_i)x^y=k_i^{\chi_{\mathbb{Z}}(\alpha_i(y))}x^{s_iy}+\bigl(k_i-k_i^{-1}\bigr)\nabla_i(x^y)\qquad
\text{if}\ k_i=u_i,
\end{equation*}
from which it follows that the $H^X$-representation $\pi^\mathcal{O}$ for $\mathbf{k}\in\mathcal{K}^{{\rm res}}$ is the restriction to $H^X$ of the quasi-polynomial $\mathbb{H}$-representation defined in~\cite[Theorem~1.1] {SSV2}.
\item Let $\Lambda\subset E$ be a $W_0$-invariant lattice containing $Q^\vee$. Then $\mathbf{F}[\Lambda]$ is a $\pi\bigl(H^X\bigr)$-submodule, and the action of $\pi(T_i)\vert_{\mathbf{F}[\Lambda]}$ can be written in terms of metaplectic Demazure--Lusztig type operators when $k_i=u_i$, see~\cite{SSV2}.
The resulting \smash{$H^X$}-representation $\pi(\cdot)\vert_{\mathbf{F}[\Lambda]}$ for $\mathbf{k}\in\mathcal{K}^{{\rm res}}$ is essentially the one introduced in~\cite[Theorem~3.7]{SSV1}. The proof of
 Theorem~\ref{theoremH} basically follows the same strategy as the proof of~\cite[Theorem~3.7]{SSV1}.
 \end{enumerate}
 \end{Remark}
\subsection{Parabolic data and parabolically induced modules}\label{section_pi_mod}
We show in this subsection that the representation $\pi^\mathcal{O}$ is a parabolically induced $H^X$-module when $\alpha_0\bigl(c^\mathcal{O}\bigr)\not=0$.

For a subset $I\subseteq\{1,\dots,r\}$, write
\begin{itemize}\itemsep=0pt
\item $W_{0,I}$ for the parabolic subgroup of $W_0$ generated by $s_i$, $i\in I$,
\item $W_0^I$ for the minimal coset representatives of $W_0/W_{0,I}$,
\item $H_{0,I}$ for the subalgebra of $H_0$ generated by $T_i$, $i\in I$.
\end{itemize}
The finite Hecke algebra $H_0$ is a free right $H_{0,I}$-module with basis $\{T_v\}_{v\in W_0^I}$, since
\begin{equation}\label{lengthadd}
\ell(vw)=\ell(v)+\ell(w)\qquad \text{for}\ v\in W_0^I \ \text{and}\ w\in W_{0,I}.
\end{equation}
Here is an immediate lift to the affine Hecke algebra $H^X$.
\begin{Lemma}\label{lemfree}\hfill
\begin{enumerate}\itemsep=0pt
\item[$(1)$]
$\{\tau(\lambda)v\}_{(\lambda,v)\in Q^\vee\times W_0^I}$ is a complete set of representatives of $W/W_{0,I}$.
\item[$(2)$] $H^X$ is a free right $H_{0,I}$-module with
basis \smash{$\bigl\{x^\lambda T_v\bigr\}_{(\lambda,v)\in Q^\vee\times W_0^I}$}.
\end{enumerate}
\end{Lemma}
\begin{proof}
The first part is a consequence of the fact that $(\lambda,w)\mapsto \tau(\lambda)w$ defines a bijection $Q^\vee\rtimes W_0\overset{\sim}{\longrightarrow} W$. For the second part, note that the multiplication map of $H^X$ restricts to
a~linear isomorphism $\mathbf{F}\smash{\bigl[Q^\vee\bigr]}\otimes H_0\overset{\sim}{\longrightarrow} H^X$.
\end{proof}

The more familiar parabolic structures on $W$ and on the associated affine Hecke algebra $H$ arise from their Coxeter type presentations. In this case it depends on a subset $J$ of $\{0,\dots,r\}$. We~write
\begin{itemize}\itemsep=0pt
\item $W_{J}$ for the parabolic subgroup of $W$ generated by $s_j$, $j\in J$,
\item $W^J$ for the minimal coset representatives of $W/W_{J}$,
\item $H_{J}$ for the subalgebra of $H$ generated by $T_j$, $j\in J$.
\end{itemize}
The length identity~\eqref{lengthadd} now also holds true for $v\in W^J$ and $w\in W_J$. As a consequence, $H$~is a free right $H_J$-module with basis
$\{T_g\}_{g\in W^J}$.

The closure $\overline{C_+}$ of the fundamental alcove $C_+$ splits in a disjoint union of facets
\begin{equation*}
\overline{C_+}=\bigsqcup_{J\subsetneq \{0,\dots,r\}}C_+^J,
\end{equation*}
with $C_+^J$ the set of vectors $y\in\overline{C_+}$ for which $\alpha_j(y)=0$ if and only if $j\in J$.
For $y\in E$ denote by $W_y\subset W$ the subgroup of $W$ fixing $y$. It is well known that
\[
W_c=W_J\qquad \text{for}\ c\in C_+^J.
\]

\begin{Definition}\label{IJ}
For a $W$-orbit $\mathcal{O}$ in $E$, we write
$J(\mathcal{O})$ for the subset of $\{0,\dots,r\}$ such that $c^\mathcal{O}\in C_+^{J(\mathcal{O})}$. We furthermore write
\begin{equation*}
I(\mathcal{O}):=J(\mathcal{O})\cap\{1,\dots,r\}.
\end{equation*}
\end{Definition}

Note that $I(\mathcal{O})=J(\mathcal{O})$ if and only if $\alpha_0\bigl(c^\mathcal{O}\bigr)\not=0$.
We will use the shorthand notations
\[
C_+^\mathcal{O},\ nW_\mathcal{O},\ W^\mathcal{O},\ H_\mathcal{O}, \ \dots
\]
for $C_+^{J(\mathcal{O})},W_{J(\mathcal{O})}, W^{J(\mathcal{O})}, H_{J(\mathcal{O})},\dots$ and
\[
W_{0,\mathcal{O}},\ W^{\mathcal{O}}_0,\ H_{0,\mathcal{O}},\ \dots
\]
for $W_{0,I(\mathcal{O})}, W_0^{I(\mathcal{O})}, H_{0,I(\mathcal{O})},\dots$.
The following lemma refines the decomposition~\eqref{FOpoldec} of the space of quasi-polynomials $\mathbf{F}[\mathcal{O}]$ as $\mathbf{F}\smash{\bigl[Q^\vee\bigr]}$-module when
$\alpha_0\bigl(c^\mathcal{O}\bigr)\not=0$.
\begin{Lemma}\label{Pcpoldec}
If $\mathcal{O}$ is a $W$-orbit such that $\alpha_0\bigl(c^\mathcal{O}\bigr)\not=0$, then
$\mathbf{F}[\mathcal{O}]$ decomposes as
\begin{equation}\label{FpoldecI}
\mathbf{F}[\mathcal{O}]=\bigoplus_{v\in W_0^\mathcal{O}}\mathbf{F}\smash{\bigl[Q^\vee\bigr]}x^{vc^\mathcal{O}}.
\end{equation}
\end{Lemma}
\begin{proof}
The assignment $gW_{0,\mathcal{O}}\mapsto gc^\mathcal{O}$ gives rise to a bijection
\[
W/W_{0,\mathcal{O}}\overset{\sim}{\longrightarrow}\mathcal{O},
\]
since $W_{c^\mathcal{O}}=W_{\mathcal{O}}=W_{0,\mathcal{O}}$ by the assumption on $\mathcal{O}$, and $W_0^{\mathcal{O}}$ is a complete set of representatives of the double coset space $\tau\smash{\bigl(Q^\vee\bigr)}\backslash W/W_{0,\mathcal{O}}$ by Lemma
\ref{lemfree}(1).
Hence $\{vc^\mathcal{O}\}_{v\in W_0^{\mathcal{O}}}=W_0c^\mathcal{O}$ is a~complete set of representatives of the $\tau\smash{\bigl(Q^\vee\bigr)}$-orbits in $\mathcal{O}$, and the lemma follows.
\end{proof}

For a $W$-orbit $\mathcal{O}$ in $E$, let $\mathbf{F}1^\mathcal{O}$
the trivial $H_\mathcal{O}$-module, defined by
\[
T_j1^\mathcal{O}=k_j1^\mathcal{O}\qquad \text{for}\ j\in J(\mathcal{O}).
\]
We will also view $\mathbf{F}1^\mathcal{O}$ as $H_{0,\mathcal{O}}$-module by restricting the action to $H_{0,\mathcal{O}}$.
Consider the induced $H^X$-module
\[
{\rm Ind}_{H_{0,\mathcal{O}}}^{H^X}\bigl(\mathbf{F}1^{\mathcal{O}}\bigr)=H^X\otimes_{H_{0,\mathcal{O}}}\mathbf{F}1^\mathcal{O}
\]
and write
\[
\mathbf{1}^{\mathcal{O}}=1\otimes_{H_{0,\mathcal{O}}}1^\mathcal{O}
\]
for its canonical cyclic vector.

\begin{Proposition}\label{indcor}
There exists a unique $H^X$-linear epimorphism
\begin{equation}\label{epiuse}
{\rm Ind}_{H_{0,\mathcal{O}}}^{H^X}\bigl(\mathbf{F}1^{\mathcal{O}}\bigr)\twoheadrightarrow
\bigl(\mathbf{F}[\mathcal{O}],\pi^\mathcal{O}\bigr)
\end{equation}
mapping $\mathbf{1}^{\mathcal{O}}$ to \smash{$x^{c^\mathcal{O}}$}. It is an isomorphism when $\alpha_0\bigl(c^\mathcal{O}\bigr)\not=0$.
\end{Proposition}

\begin{proof}
For the first statement, we need to show that the assignment $h\mathbf{1}^{\mathcal{O}}\mapsto\pi(h)x^{c^\mathcal{O}}$ for $h\in H^X$ is well defined. It suffices to note that
\[
\pi(T_i)x^{c^\mathcal{O}}=k_ix^{c^\mathcal{O}}\qquad \text{for}\ i\in I(\mathcal{O}).
\]
But for $i\in I(\mathcal{O})$, we have \smash{$s_ic^\mathcal{O}=c^{\mathcal{O}}$}, and hence
\[
\pi(T_i)x^{c^\mathcal{O}}=\kappa_{s_i}^\mathcal{O} x^{s_ic^\mathcal{O}}=k_ix^{c^\mathcal{O}}
\]
by Corollary~\ref{corpsi}\,(2) and~\eqref{kappawy}.

For the second statement, note first that
\begin{equation}\label{basis}
\bigl\{x^\lambda T_v\mathbf{1}^{\mathcal{O}}\mid (\lambda,v)\in Q^\vee\times W_0^\mathcal{O}\bigr\}
\end{equation}
is a basis of \smash{${\rm Ind}_{H_{0,\mathcal{O}}}^{H^X}\bigl(\mathbf{F}1^\mathcal{O}\bigr)$}. By Corollary~\ref{corpsi}\,(2), the basis element $x^\lambda T_v\mathbf{1}^{\mathcal{O}}$
is mapped by the epimorphism~\eqref{epiuse} to
\[
\pi(x^\lambda T_v)x^{c^\mathcal{O}}=\kappa_v^\mathcal{O}x^{\lambda+vc^\mathcal{O}}.
\]
By Lemma~\ref{Pcpoldec}, we conclude that the epimorphism~\eqref{epiuse} maps the basis~\eqref{basis} of
${\rm Ind}_{H_{0,\mathcal{O}}}^{H^X}\bigl(\mathbf{F}1^\mathcal{O}\bigr)$
to a basis of $\mathbf{F}[\mathcal{O}]$, which concludes the proof of the second statement.
\end{proof}

The natural analog of Proposition~\ref{indcor} for $W$-orbits $\mathcal{O}$ with $\alpha_0\bigl(c^\mathcal{O}\bigr)=0$ (i.e., with $0\in J(\mathcal{O})$) requires the extension of the $H^X$-action $\pi^\mathcal{O}$ on $\mathbf{F}[\mathcal{O}]$ to an action of the double affine Hecke algebra $\mathbb{H}$. This will be the subject of the next subsection.

\subsection[The quasi-polynomial representation of H]{The quasi-polynomial representation of $\boldsymbol{\mathbb{H}}$}\label{wtsection}
We now promote the quasi-polynomial representation $\pi^\mathcal{O}$ of the affine Hecke algebra $H^X$ to a~family of representations of the double affine Hecke algebra $\mathbb{H}$.
The number of additional parameters
depends on the facet \smash{$C^{J(\mathcal{O)}}_+$} containing $c^\mathcal{O}$.
For $\mathbf{k}\in\mathcal{K}^{{\rm res}}$ the extended representations are the quasi-polynomial $\mathbb{H}$-representations $\pi_{c^\mathcal{O},t}$ introduced in~\cite[Theorem~1.1]{SSV2}. For~${\mathbf{k}\in\mathcal{K}}$ and~$\Phi_0$ of type ${\rm C}_r$, $r\geq 1$, they will give quasi-polynomial extensions of the polynomial representation of the
type ${\rm C}^\vee{\rm C}_r$ double affine Hecke algebra $\mathbb{H}$~\cite{No,Sa}.

Recall the definition of the $\mathbf{F}$-torus $\mathbf{T}$ (see Definition~\ref{T_def}). For $J\subsetneq\{0,\dots,r\}$, consider its affine subtori
\begin{equation*}
\mathbf{T}_J:=\bigl\{t\in\mathbf{T}\mid t^{\alpha_j^\vee}=1\ \forall j\in J\bigr\}.
\end{equation*}
For a $W$-orbit $\mathcal{O}$ in $E$, write
\[\mathbf{T}_\mathcal{O}:=\mathbf{T}_{J(\mathcal{O})}.
\]
Note that
\begin{equation}\label{sat}
s_at=t\bigl(t^{-a^\vee}\bigr)^{\overline{a}}\qquad \text{for}\ a\in\Phi\ \text{and}\ t\in\mathbf{T},
\end{equation}
where \smash{$t\bigl(t^{-a^\vee}\bigr)^{\overline{a}}\in\mathbf{T}$} is viewed as character of $Q^\vee$ by \smash{$\lambda\mapsto t^\lambda\bigl(t^{-a^\vee}\bigr)^{\overline{a}(\lambda)}$} for
$\lambda\in Q^\vee$. By~\eqref{sat}, we~have
$s_at=t$ if $t^{a^\vee}=1$, and so
\begin{equation}\label{containedinv}
\mathbf{T}_J\subseteq\mathbf{T}^{W_J}:=\{t\in\mathbf{T} \mid gt=t\ \forall g\in W_J\}.
\end{equation}
Note that $\mathbf{T}_J$ and $\mathbf{T}^{W_J}$ are sub-tori of $\mathbf{T}$ when $0\not\in J$.

In~\cite[Lemma 4.2]{SSV2}, the $W_0$-action on $\mathbf{F}[\mathcal{O}]$ was extended to a family of $W$-actions on $\mathbf{F}[\mathcal{O}]$, parametrised by
$t\in\mathbf{T}_\mathcal{O}$. These actions are compatible with the $W$-action~\eqref{Wa} on $\mathbf{F}\smash{\bigl[Q^\vee\bigr]}$ by $q$-dilations and reflections.
The definition of this action requires the following definition.

\begin{Definition}\label{gy}
For $y\in E$, write $\mathbf{g}_y\in W$ for the unique element of shortest length in $W$ such that $\mathbf{g}_y^{-1}y\in\overline{C_+}$.
\end{Definition}
Note that if $y$ lies in the $W$-orbit $\mathcal{O}$ of $E$,
then $\mathbf{g}_y$ is the unique element in $W^\mathcal{O}$ such that
$y=\mathbf{g}_yc^\mathcal{O}$.
\begin{Lemma}\label{tdefaction}
Let $\mathcal{O}$ be a $W$-orbit in $E$. For $t\in\mathbf{T}_\mathcal{O}$,
 the formulas
\begin{gather}
w_tx^y:=w(x^y)=x^{wy},\qquad w\in W_0, \nonumber\\
\tau(\lambda)_tx^y:=(\mathbf{g}_yt)^{-\lambda}x^y,\qquad \lambda\in Q^\vee\label{taction}
\end{gather}
for $y\in\mathcal{O}$ define a linear left $W$-action on $\mathbf{F}[\mathcal{O}]$
 satisfying
\begin{equation}\label{wcomp}
g_t(pf)=(gp)(g_tf)
\end{equation}
for $g\in W$, $p\in\mathbf{F}\smash{\bigl[Q^\vee\bigr]}$ and $f\in\mathbf{F}[\mathcal{O}]$.
\end{Lemma}
It is instructive to recall the proof of~\eqref{wcomp} for $g=\tau(\lambda)$, $\lambda\in Q^\vee$. First note that by~\eqref{containedinv},
we may replace $\mathbf{g}_y$ in~\eqref{taction} by any other representative of the coset $\mathbf{g}_yW_\mathcal{O}$.
We then have for~${\mu\in Q^\vee}$,
\begin{align}
\tau(\lambda)_tx^{\mu+y}&=(\mathbf{g}_{\mu+y}t)^{-\lambda}x^{\mu+y} \nonumber\\
&=(\tau(\mu)\mathbf{g}_yt)^{-\lambda}x^{\mu+y} \nonumber\\
&=q^{-\langle\lambda,\mu\rangle}(\mathbf{g}_y t)^{-\lambda}x^{\mu+y} \nonumber\\
&=\bigl(q^{-\langle\lambda,\mu\rangle}x^\mu\bigr)\bigl((\mathbf{g}_yt)^{-\lambda}x^y\bigr),\label{tcomp}
\end{align}
where we used $\mathbf{g}_{y+\mu}W_\mathcal{O}=\tau(\mu)\mathbf{g}_yW_\mathcal{O}$ for the second equality, and~\eqref{Waction} for the third equality.
The last line in~\eqref{tcomp} clearly equals $(\tau(\lambda)x^\mu)(\tau(\lambda)_tx^y)$. Note in particular that the family of $W$-actions~\eqref{taction} depends on $q$ via the $W$-action~\eqref{Waction} on $\mathbf{T}$.

Note that by~\eqref{tcomp} we have the formula
\[
\tau(\lambda)_t\vert_{\mathbf{F}[Q^\vee]x^y}=(\mathbf{g}_yt)^{-\lambda}(x^y\circ\tau(\lambda)\circ x^{-y})\vert_{\mathbf{F}[Q^\vee]x^y},
\]
where $x^{\pm y}$ are regarded as multiplication operators on $\mathbf{F}[E]$.

The following result extends~\cite[Theorem~1.1]{SSV2} to multiplicity functions $\mathbf{k}\in\mathcal{K}$.
\begin{Theorem}\label{mthm}
Let $\mathcal{O}$ be a $W$-orbit in $E$ and $t\in\mathbf{T}_\mathcal{O}$. The formulas
\begin{gather}
\pi_{t}^\mathcal{O}\bigl(x^\lambda\bigr)x^y:=x^{y+\lambda}, \nonumber\\
\pi_{t}^\mathcal{O}(T_j)x^y:=k_j^{\chi_e(\overline{\alpha}_j(y))}u_j^{\chi_o(\overline{\alpha}_j(y))}s_{j,t}x^y+\bigl(k_j-k_j^{-1}\bigr)\nabla_j^e(x^y)+
\bigl(u_j-u_j^{-1}\bigr)\nabla_j^o(x^y)
\label{qprepHH}
\end{gather}
for $j=0,\dots,r$, $\lambda\in Q^\vee$ and $y\in\mathcal{O}$ define a representation $\pi_{t}^\mathcal{O}\colon \mathbb{H}\rightarrow\textup{End}(\mathbf{F}[\mathcal{O}])$.
\end{Theorem}
\begin{Remark}
Note that $\pi_{t}^\mathcal{O}\vert_{H^X}$ is the restriction $\pi^{\mathcal{O}}$ of the quasi-polynomial $H^X$-represen\-tation $\pi$ from Theorem~\ref{theoremH} to $\mathbf{F}[\mathcal{O}]$. Furthermore, for restricted multiplicity
parameters $\mathbf{k}\in\mathcal{K}^{{\rm res}}$ we have
\[
\pi_{t}^\mathcal{O}(T_j)=k_j^{\chi_{\mathbb{Z}}(\overline{\alpha}_j(y))}s_{j,t}x^y+\bigl(k_j-k_j^{-1}\bigr)\nabla_j(x^y)
\]
since $k_j=u_j$, hence $\pi_{t}^\mathcal{O}$ then coincides with the quasi-polynomial representation $\pi_{c^\mathcal{O},t}$ defined in~\cite[Theorem~1.1]{SSV2}.
\end{Remark}
\begin{proof}
Consider $\mathbb{H}$ as left regular $\mathbb{H}$-module.
Relative
to the $\mathbf{F}$-basis
\[
\bigl\{x^\lambda T_w Y^\mu\mid \lambda,\mu\in Q^\vee,\, w\in W_0\bigr\}
\]
of $\mathbb{H}$, the $H$-action is explicitly given by
\begin{gather}
T_jx^\lambda T_wY^\mu=s_j\bigl(x^\lambda\bigr)T_{s_{\overline{\alpha}_j}w}Y^{\mu-w^{-1}s_j(0)} \nonumber\\
\hphantom{T_jx^\lambda T_wY^\mu=}{}
+\big(k_j-k_j^{-1}\big)\Biggl(\frac{1-x^{(2\chi_-(w^{-1}\overline{\alpha}_j)-\overline{\alpha}_j(\lambda))\alpha_j^\vee}}{1-x^{2\alpha_j^\vee}}\Biggr)x^\lambda T_wY^\mu\nonumber\\
\hphantom{T_jx^\lambda T_wY^\mu=}{}
+\big(u_j-u_j^{-1}\big)\Biggl(\frac{x^{\alpha_j^\vee}-x^{(1-\overline{\alpha}_j(\lambda))\alpha_j^\vee}}{1-x^{2\alpha_j^\vee}}\Biggr)x^\lambda T_w Y^\mu\label{regularHH}
\end{gather}
for $w\in W_0$, $\lambda,\mu\in Q^\vee$ and $j=0,\dots,r$. For $j\in\{1,\dots,r\}$, formula~\eqref{regularHH} follows immediately from~\eqref{regularHX}.
For $j=0$, formula~\eqref{regularHH} follows by commuting $T_0$ and $x^\lambda$ using the cross relation~\eqref{crossX} and then applying
the identity
\begin{equation}\label{T0Tw}
T_0T_w=T_{s_{\overline{\alpha}_0}w}Y^{-w^{-1}s_0(0)}+\chi_-\bigl(w^{-1}\overline{\alpha}_0\bigr)\bigl(k_0-k_0^{-1}\bigr)T_w
\end{equation}
in $H$. For the proof of~\eqref{T0Tw}, first note that it is equivalent to the identity
\begin{equation}\label{Hid}
T_0^{\chi(w^{-1}\varphi)}T_{s_\varphi w}=T_wY^{w^{-1}\varphi^\vee}
\end{equation}
in $H$ since $s_0(0)=\varphi^\vee$, $\overline{\alpha}_0=-\varphi$ and $T_0^{-1}=T_0-k_0+k_0^{-1}$.
For a proof of~\eqref{Hid} see, for instance,~\cite[(3.3.6)]{Ma}.

Let $\mathcal{O}$ be a $W$-orbit in $E$.
Recall the linear epimorphism $\psi^\mathcal{O}\colon H^X\twoheadrightarrow\mathbf{F}[\mathcal{O}]$,
defined by~\eqref{psicmap}, which is an epimorphism of $H^X$-modules by Corollary~\ref{corpsi}\,(2). The goal is to find extensions of~$\psi^{\mathcal{O}}$ to linear epimorphisms $\mathbb{H}\twoheadrightarrow\mathbf{F}[\mathcal{O}]$ such that
\begin{enumerate}\itemsep=0pt
\item[(1)] their kernels are left ideals in $\mathbb{H}$,
\item[(2)]
right multiplication by $\mathbf{F}_Y\smash{\bigl[Q^\vee\bigr]}$ is turned into multiplication by a linear character of $\mathbf{F}_Y\smash{\bigl[Q^\vee\bigr]}$.
\end{enumerate}
This will upgrade the $H^X$-action on $\mathbf{F}[\mathcal{O}]$ to a family of $\mathbb{H}$-actions satisfying the additional property that the quasi-monomial \smash{$x^{c^\mathcal{O}}$} will be a simultaneous eigenvector
for the action of $\mathbf{F}_Y\smash{\bigl[Q^\vee\bigr]}$. The family of extended $\mathbb{H}$-actions on $\mathbf{F}[\mathcal{O}]$ will be natural parametrised by the associated linear characters of $\mathbf{F}_Y\smash{\bigl[Q^\vee\bigr]}$, which in turn can be described by an affine subtorus of the form $\mathfrak{s}\mathbf{T}_{\mathcal{O}}$ for a specific basepoint $\mathfrak{s}\in\mathbf{T}$, so be determined in due course.

So our starting point will be the desired property (2). Fix $t\in\mathbf{T}_{\mathcal{O}}$ and consider the extension of $\psi^\mathcal{O}$ to a ($t$-dependent) linear map $\psi_{t}^\mathcal{O}\colon \mathbb{H}\twoheadrightarrow\mathbf{F}[\mathcal{O}]$ by
\begin{equation}\label{psict}
\psi_{t}^\mathcal{O}\bigl(x^\lambda T_wY^\mu\bigr):=(\mathfrak{s}t)^{-\mu}\psi^\mathcal{O}\bigl(x^\lambda T_w\bigr)=\kappa_w^\mathcal{O}(\mathfrak{s}t)^{-\mu}x^{\lambda+wc^\mathcal{O}}
\end{equation}
for $w\in W_0$ and $\lambda,\mu\in Q^\vee$ \big(recall here that $\kappa_w^\mathcal{O}\in\mathbf{F}^\times$ is the explicit scalar defined by~\eqref{kappawy}\big).
By the first formula of~\eqref{psict} and Corollary~\ref{corpsi}\,(2), it follows that
\[
\psi_{t}^\mathcal{O}(hh^\prime)
=\pi_{t}^\mathcal{O}(h)\psi_{t}^\mathcal{O}(h^\prime)\qquad \text{for}\, h\in H^X\, \text{and}\, h^\prime\in\mathbb{H}.
\]
In particular, the kernel of $\psi_{t}^\mathcal{O}$ is a left $H^X$-submodule in $\mathbb{H}$.
To meet the first property (1), we~now fine-tune the choice of $\mathfrak{s}\in\mathbf{T}$ such that
\[
\psi_{t}^{\mathcal{O}}(T_0h^\prime)=\mathcal{D}_0\bigl(\psi_{t}^\mathcal{O}(h^\prime)\bigr) \qquad \textup{ for all }\, h^\prime\in\mathbb{H}
\]
for some $\mathcal{D}_0\in\textup{End}(\mathbf{F}[\mathcal{O}])$.

Note that for $\lambda\in Q^\vee$ and $w\in W_0$,
\begin{equation*}
s_{0,t}x^{\lambda+wc^\mathcal{O}}=s_0\bigl(x^\lambda\bigr)\bigl(s_{0,t}x^{wc^\mathcal{O}}\bigr)=t^{w^{-1}\varphi^\vee}s_0\bigl(x^\lambda\bigr)x^{s_\varphi wc^\mathcal{O}}
\end{equation*}
and hence, by~\eqref{regularHH},
\begin{gather}
\frac{(\mathfrak{s}t)^\mu}{\kappa_w^\mathcal{O}}\psi_{t}^\mathcal{O}(T_0x^\lambda T_wY^\mu)=\frac{\kappa_{s_\varphi w}^\mathcal{O}\,\mathfrak{s}^{w^{-1}\varphi^\vee}}{\kappa_w^\mathcal{O}}s_{0,t}x^{\lambda+wc^\mathcal{O}}\nonumber\\
\hphantom{\frac{(\mathfrak{s}t)^\mu}{\kappa_w^\mathcal{O}}\psi_{t}^\mathcal{O}(T_0x^\lambda T_wY^\mu)=}{}
+(k_0-k_0^{-1})\Biggl(\frac{1-x^{(2\chi_-(w^{-1}\overline{\alpha}_0)-\overline{\alpha}_0(\lambda))\alpha_0^\vee}}{1-x^{2\alpha_0^\vee}}\Biggr)x^{\lambda+wc^\mathcal{O}}\nonumber\\
\hphantom{\frac{(\mathfrak{s}t)^\mu}{\kappa_w^\mathcal{O}}\psi_{t}^\mathcal{O}(T_0x^\lambda T_wY^\mu)=}{}
+(u_0-u_0^{-1})\Biggl(\frac{x^{\alpha_0^\vee}-x^{(1-\overline{\alpha}_0(\lambda))\alpha_0^\vee}}{1-x^{2\alpha_0^\vee}}\Biggr)x^{\lambda+wc^\mathcal{O}}.
\label{spoint0}
\end{gather}
So we need to fine-tune $\mathfrak{s}\in\mathbf{T}$ such that the right-hand side of~\eqref{spoint0} can be written
$\mathcal{D}_0\bigl(x^{\lambda+wc}\bigr)$ for
some $\mathcal{D}_0\in\textup{End}(\mathbf{F}[\mathcal{O}])$.
We follow the proof of Theorem~\ref{theoremH}, but it requires some necessary additional computations due to the presence of the additional parameters $\mathfrak{s}$ and $t$ (compare also with~\cite[Section~5]{SSV2}, which deals with the case that $\mathbf{k}\in\mathcal{K}^{{\rm res}}$).

We first consider the case that $\Phi_0$ is of type ${\rm C}_r$, $r\geq 1$, so that $\alpha_0\smash{\bigl(Q^\vee\bigr)}=\mathbb{Z}_o$.
Fix $\lambda\in Q^\vee$ and $w\in W_0$ and set
\[
y:=\lambda+wc^\mathcal{O}.
\]
The second line of~\eqref{spoint0} can then be rewritten using the formula
\begin{equation}\label{claim3}
\Biggl(\frac{1-x^{(2\chi_-(w^{-1}\overline{\alpha}_0)-\overline{\alpha}_0(\lambda))\alpha_0^\vee}}{1-x^{2\alpha_0^\vee}}\Biggr)x^{y}=
\nabla_0^e(x^y)
+\chi_-\bigl(w^{-1}\overline{\alpha}_0\bigr)\chi_e\bigl(\overline{\alpha}_0\bigl(wc^\mathcal{O}\bigr)\bigr)s_{0,t}x^{y}.
\end{equation}
To prove~\eqref{claim3}, note that by~\eqref{W0C} we have
\begin{equation}
\label{sstep10}
2\chi_-\bigl(w^{-1}\overline{\alpha}_0\bigr)-\overline{\alpha}_0(\lambda)=
\begin{cases}
2-\lfloor\overline{\alpha}_0(y)\rfloor_e &\text{if}\ w^{-1}\overline{\alpha}_0\in\Phi_0^-\ \text{and}\ \overline{\alpha}_0\bigl(wc^\mathcal{O}\bigr)=0,\\
-\lfloor\overline{\alpha}_0(y)\rfloor_e &\text{otherwise},
\end{cases}
\end{equation}
and that $\overline{\alpha}_0\bigl(wc^\mathcal{O}\bigr)=0$ if and only if $\chi_e\bigl(\overline{\alpha}_0\bigl(wc^\mathcal{O}\bigr)\bigr)=1$.
This immediately implies~\eqref{claim3} unless
$w^{-1}\overline{\alpha}_0\in\Phi_0^-$ and $\overline{\alpha}_0\bigl(wc^\mathcal{O}\bigr)=0$.
So suppose now that $w^{-1}\overline{\alpha}_0\in\Phi_0^-$ and $\overline{\alpha}_0\bigl(wc^\mathcal{O}\bigr)=0$. Then $s_\varphi wc^\mathcal{O}=wc^\mathcal{O}$, and hence
\[
s_{0,t}x^y=s_0\bigl(x^\lambda\bigr)\bigl(s_{0,t}x^{wc^\mathcal{O}}\bigr)=\bigl(x^{\lambda-\overline{\alpha}_0(\lambda)\alpha_0^\vee}\bigr)\bigl(t^{w^{-1}\varphi^\vee}x^{wc^\mathcal{O}}\bigr)
=t^{w^{-1}\varphi^\vee}x^{\lambda-\overline{\alpha}_0(\lambda)\alpha_0^\vee+wc^\mathcal{O}}.
\]
But \smash{$\bigl(w^{-1}\varphi\bigr)\bigl(c^\mathcal{O}\bigr)=0$}, hence
\[
w^{-1}\varphi^\vee\in\Phi_0^\vee\cap\bigoplus_{i\in I(\mathcal{O})}\mathbb{Z}\alpha_i^\vee.
\]
 Since $t\in\mathbf{T}_\mathcal{O}$,
we conclude that ${t^{w^{-1}\varphi^\vee}=1}$, so
\begin{equation}\label{s0simple}
s_{0,t}x^y=x^{\lambda-\overline{\alpha}_0(\lambda)\alpha_0^\vee+wc^\mathcal{O}}.
\end{equation}
Then~\eqref{claim3} follows by combining the first case of~\eqref{sstep10} and~\eqref{s0simple}.

Similarly, we rewrite the third line of~\eqref{spoint0} using the formula
\begin{equation}\label{claim4}
\Biggl(\frac{x^{\alpha_0^\vee}-x^{(1-\overline{\alpha}_0(\lambda))\alpha_0^\vee}}{1-x^{2\alpha_0^\vee}}\Biggr)x^y=
\nabla_0^o(x^y)+
\chi_+\bigl(w^{-1}\overline{\alpha}_0\bigr)\chi_o\bigl(\overline{\alpha}_0\bigl(wc^\mathcal{O}\bigr)\bigr)
s_{0,t}x^y.
\end{equation}
For the proof of~\eqref{claim4}, we now use that
\begin{equation}\label{sstep20}
1-\overline{\alpha}_0(\lambda)=
\begin{cases}
2-\lfloor\overline{\alpha}_0(y)\rfloor_o &\text{if}\ \overline{\alpha}_0\bigl(wc^\mathcal{O}\bigr)=1,\\
-\lfloor\overline{\alpha}_0(y)\rfloor_o &\text{otherwise},
\end{cases}
\end{equation}
and the observation that $\overline{\alpha}_0\bigl(wc^\mathcal{O}\bigr)=1$ is equivalent to
$\chi_+\bigl(w^{-1}\overline{\alpha}_0\bigr)\chi_o\bigl(\overline{\alpha}_0\bigl(wc^\mathcal{O}\bigr)\bigr)\!=\!1$. Then~\eqref{claim4} is immediate if $\overline{\alpha}_0\bigl(wc^\mathcal{O}\bigr)\not=1$.
So suppose now that $\overline{\alpha}_0\bigl(wc^\mathcal{O}\bigr)=1$. Then
\begin{align}
s_{0,t}x^y
&=t^{w^{-1}\varphi^\vee}x^{\lambda-\overline{\alpha}_0(\lambda)\alpha_0^\vee+s_\varphi wc^\mathcal{O}}\nonumber\\
&=q_\varphi t^{w^{-1}\varphi^\vee}x^{\lambda-(1+\overline{\alpha}_0(\lambda))\alpha_0^\vee+wc^\mathcal{O}}=
x^{\lambda-(1+\overline{\alpha}_0(\lambda))\alpha_0^\vee+wc^\mathcal{O}}.\label{s0help}
\end{align}
The last equality follows from the fact that the affine root $a:=\bigl(w^{-1}\varphi,1\bigr)$ satisfies $a\bigl(c^\mathcal{O}\bigr)=0$, so \smash{$a^\vee\in\Phi^+\cap\bigoplus_{j\in J(\mathcal{O})}\mathbb{Z}\alpha_j^\vee$} and hence
\smash{$q_\varphi t^{w^{-1}\varphi^\vee}=t^{a^\vee}=1$} (cf.~\cite[Lemma 3.4]{SSV2}). Then~\eqref{claim4} follows from the first case of~\eqref{sstep20} and~\eqref{s0help}.

Returning now to~\eqref{spoint0},
we conclude from~\eqref{claim3} and~\eqref{claim4} that
\begin{equation}\label{formula0}
\frac{(\mathfrak{s}t)^\mu}{\kappa_w^\mathcal{O}}\psi_{t}^\mathcal{O}\bigl(T_0x^\lambda T_wY^\mu\bigr)
={\rm coeff}_w^\mathcal{O}s_{0,t}x^y+\bigl(k_0-k_0^{-1}\bigr)\nabla_0^e(x^y)+\bigl(u_0-u_0^{-1}\bigr)\nabla_0^o(x^y)
\end{equation}
with
\begin{gather}
{\rm coeff}_w^\mathcal{O}:=\frac{\kappa_{s_\varphi w}^\mathcal{O}\mathfrak{s}^{w^{-1}\varphi^\vee}}{\kappa_w^\mathcal{O}}+\bigl(k_0-k_0^{-1}\bigr)\chi_-(w^{-1}\overline{\alpha}_0)\chi_e \bigl(\overline{\alpha}_0\bigl(wc^\mathcal{O}\bigr)\bigr)\nonumber\\
\hphantom{{\rm coeff}_w^\mathcal{O}:=}{}
+\bigl(u_0-u_0^{-1}\bigr)\chi_+\bigl(w^{-1}\overline{\alpha}_0\bigr)\chi_o\bigl(\overline{\alpha}_0\bigl(wc^\mathcal{O}\bigr)\bigr).
\label{coeff}
\end{gather}
Formula~\eqref{formula0} also holds true when $\Phi_0$ is not of type ${\rm C}_r$. In fact, if $\Phi_0$ is not of type ${\rm C}_r$ then $\mathcal{K}=\mathcal{K}^{{\rm res}}$ and $k_0=u_0$,
and the formula can be recovered from~\cite[Section~5]{SSV2}. It can also be derived directly, similarly as the proof of~\eqref{formula} when $\alpha_i\smash{\bigl(Q^\vee\bigr)}=\mathbb{Z}$.
So we will now continue the proof of the theorem for $\Phi_0$ of any type, taking~\eqref{formula0} as the starting point.

Properties of the normalization factors $\kappa_w^\mathcal{O}$~\eqref{kappawy}
were derived in~\cite{SSV2} for restricted parameters $\mathbf{k}\in\mathcal{K}^{{\rm res}}$, which led to an explicit expression of the quotient $\kappa_{s_\varphi w}^\mathcal{O}/\kappa_w^\mathcal{O}$ (see~\cite[Lemma 5.6\,(1)]{SSV2}, as well as case 2 of the proof of~\cite[Lemma 5.9]{SSV2}).
This can be easily extended to $\kappa_w^\mathcal{O}$ for arbitrary parameters $\mathbf{k}\in\mathcal{K}$. It leads to the formula
\begin{gather*}
\frac{\kappa_{s_\varphi w}^\mathcal{O}}{\kappa_w^\mathcal{O}}=\mathbf{k}_\varphi^{-\chi(w^{-1}\varphi)\chi_e(\varphi(wc^\mathcal{O}))}
\mathbf{k}_{\frac{\varphi}{2}}^{\chi(w^{-1}\varphi)\chi_o(\varphi(wc^\mathcal{O}))}
\prod_{\alpha\in\Phi_0^+}\mathbf{k}_\alpha^{\chi_e(\alpha(c^\mathcal{O}))\alpha(w^{-1}\varphi^\vee)}\mathbf{k}_{\frac{\alpha}{2}}^{-\chi_o(\alpha(c^\mathcal{O}))\alpha(w^{-1}\varphi^\vee)},
\end{gather*}
where $\chi$ is given by~\eqref{chi}. Hence
\begin{gather*}
{\rm coeff}_w^\mathcal{O}=k_r^{-\chi(w^{-1}\varphi)\chi_e(\varphi(wc^\mathcal{O}))}u_r^{\chi(w^{-1}\varphi)\chi_o(\varphi(wc^\mathcal{O}))}
\biggl(\mathfrak{s}\prod_{\alpha\in\Phi_0^+}\mathbf{k}_\alpha^{\chi_e(\alpha(c^\mathcal{O}))\alpha} \mathbf{k}_{\frac{\alpha}{2}}^{-\chi_o(\alpha(c^\mathcal{O}))\alpha}\biggr)^{w^{-1}\varphi^\vee}\\
\hphantom{{\rm coeff}_w^\mathcal{O}=}{}
+\bigl(k_0-k_0^{-1}\bigr)\chi_+\bigl(w^{-1}\varphi\bigr)\chi_e\bigl(\varphi\bigl(wc^\mathcal{O}\bigr)\bigr)+
\bigl(u_0-u_0^{-1}\bigr)\chi_-\bigl(w^{-1}\varphi\bigr)\chi_o\bigl(\varphi\bigl(wc^\mathcal{O}\bigr)\bigr),
\end{gather*}
where the factor in big brackets in the first line is considered as element in $\mathbf{T}$ with value at~${\lambda\in Q^\vee}$ given by
\[
\mathfrak{s}^\lambda\prod_{\alpha\in\Phi_0^+}\mathbf{k}_\alpha^{\chi_e(\alpha(c^\mathcal{O}))\alpha(\lambda)}\,\mathbf{k}_{\frac{\alpha}{2}}^{-\chi_o(\alpha(c^\mathcal{O}))\alpha(\lambda)}.
\]
Now we pick $\mathfrak{s}\in\mathbf{T}$ to be
\begin{equation}\label{s}
\mathfrak{s}_\mathcal{O}:=\prod_{\alpha\in\Phi_0^+}\bigl(\mathbf{k}_{\alpha}\mathbf{k}_{(\alpha,1)}\bigr)^{-\frac{\chi_e(\alpha(c^\mathcal{O}))}{2}\alpha}
\bigl(\mathbf{k}_{\frac{\alpha}{2}}\mathbf{k}_{(\frac{\alpha}{2},\frac{1}{2})}\bigr)^{\frac{\chi_o(\alpha(c^\mathcal{O}))}{2}\alpha},
\end{equation}
whose value at $\lambda\in Q^\vee$ is
\[
\mathfrak{s}^\lambda_\mathcal{O}:=\prod_{\alpha\in\Phi_0^+}\bigl(\mathbf{k}_{\alpha}\mathbf{k}_{(\alpha,1)}\bigr)^{-\frac{\chi_e(\alpha(c^\mathcal{O}))\alpha(\lambda)}{2}}
\bigl(\mathbf{k}_{\frac{\alpha}{2}}\mathbf{k}_{(\frac{\alpha}{2},\frac{1}{2})}\bigr)^{\frac{\chi_o(\alpha(c^\mathcal{O}))\alpha(\lambda)}{2}}.
\]
When \smash{$\alpha\bigl(Q^\vee\bigr)=\mathbb{Z}$}, the factor in this product should be read as
\[
\mathbf{k}_\alpha^{-\chi_e(\alpha(c^\mathcal{O}))\alpha(\lambda)}\,\mathbf{k}_{\frac{\alpha}{2}}^{\chi_o(\alpha(c^\mathcal{O}))\alpha(\lambda)}=
\mathbf{k}_{\alpha}^{(\chi_o(\alpha(c^\mathcal{O}))-\chi_e(\alpha(c^\mathcal{O})))\alpha(\lambda)}
\]
\big(recall that \smash{$\mathbf{k}_\alpha=\mathbf{k}_{(\alpha,1)}=\mathbf{k}_{\frac{\alpha}{2}}=
\mathbf{k}_{(\frac{\alpha}{2},\frac{1}{2})}$} when $\alpha\smash{\bigl(Q^\vee\bigr)}=\mathbb{Z}$\big).

For the remainder of the proof, we set $\mathfrak{s}=\mathfrak{s}_\mathcal{O}$, hence the linear map $\psi_t^\mathcal{O}\colon \mathbb{H}\twoheadrightarrow\mathbf{F}[\mathcal{O}]$ is now given by
\begin{equation}\label{psitO}
\psi_t^\mathcal{O}\bigl(x^\lambda T_wY^\mu\bigr)=\kappa_w^\mathcal{O}(\mathfrak{s}_\mathcal{O}t)^{-\mu}x^{\lambda+wc^\mathcal{O}}
\end{equation}
for $\lambda,\mu\in Q^\vee$ and $w\in W_0$.
We get
\begin{gather}
{\rm coeff}_w^\mathcal{O}=k_r^{-\chi(w^{-1}\varphi)\chi_e(\varphi(wc^\mathcal{O}))}
u_r^{\chi(w^{-1}\varphi)\chi_o(\varphi(wc^\mathcal{O}))}\nonumber\\
\hphantom{{\rm coeff}_w^\mathcal{O}=}\quad{}
 \times\prod_{\alpha\in\Phi_0^+: \alpha\smash{(Q^\vee)}=\mathbb{Z}_e}
\bigl(\mathbf{k}_\alpha \mathbf{k}_{(\alpha,1)}^{-1}\bigr)^{\frac{\chi_e(\alpha(c^\mathcal{O}))\alpha(w^{-1}\varphi^\vee)}{2}} \bigl(\mathbf{k}_{\frac{\alpha}{2}}^{-1}\mathbf{k}_{(\frac{\alpha}{2},\frac{1}{2})}\bigr)^{\frac{\chi_o(\alpha(c^\mathcal{O}))\alpha(w^{-1}\varphi^\vee)}{2}}\nonumber\\
\hphantom{{\rm coeff}_w^\mathcal{O}=}{}
+\bigl(k_0-k_0^{-1}\bigr)\chi_+\bigl(w^{-1}\varphi\bigr)\chi_e\bigl(\varphi\bigl(wc^\mathcal{O}\bigr)\bigr)+
\bigl(u_0-u_0^{-1}\bigr)\chi_-\bigl(w^{-1}\varphi\bigr)\chi_o\bigl(\varphi\bigl(wc^\mathcal{O}\bigr)\bigr).
\label{coeffmid}
\end{gather}
Recall that there only exist roots $\alpha\in\Phi_0$ with $\alpha\smash{\bigl(Q^\vee\bigr)}=\mathbb{Z}_e$ when $\Phi_0$ is of type ${\rm C}_r$. In this case
\smash{$\bigl\{\alpha\in\Phi_0\mid \alpha\bigl(Q^\vee\bigr)=\mathbb{Z}_e\bigr\}$} is the set $\Phi_{0,\ell}$ of long roots in $\Phi_0$, and
\[
\bigl\{\alpha\in\Phi_{0,\ell} \mid \alpha\bigl(w^{-1}\varphi^\vee\bigr)\not=0\bigr\}=\big\{w^{-1}\varphi,-w^{-1}\varphi\big\}.
\]
We thus have
\begin{gather*}
\prod_{\alpha\in\Phi_0^+\colon \alpha\smash{(Q^\vee)}=\mathbb{Z}_e}
\bigl(\mathbf{k}_\alpha \mathbf{k}_{(\alpha,1)}^{-1}\bigr)^{\frac{\chi_e(\alpha(c^\mathcal{O}))\alpha(w^{-1}\varphi^\vee)}{2}} \bigl(\mathbf{k}_{\frac{\alpha}{2}}^{-1}\mathbf{k}_{(\frac{\alpha}{2},\frac{1}{2})}\bigr)^{\frac{\chi_o(\alpha(c^\mathcal{O}))\alpha(w^{-1}\varphi^\vee)}{2}}\\
\qquad{}=
\bigl(k_0^{-1}k_r\bigr)^{\chi(w^{-1}\varphi)\chi_e(\varphi(wc^\mathcal{O}))}
\bigl(u_0u_r^{-1}\bigr)^{\chi(w^{-1}\varphi)\chi_o(\varphi(wc^\mathcal{O}))}.\label{hhhelp}
\end{gather*}
Note that this formula is correct for $\Phi_0$ of arbitrary type. Indeed, if $\Phi_0$ is not of type $C_r$ then the product on the left-hand side is an empty product, and the right-hand side also reduces to~$1$.
Substituting~\eqref{hhhelp} into
\eqref{coeffmid}, we see that
the dependence on $k_r$ and $u_r$ drops out, and we end up with the formula
\begin{gather*}
{\rm coeff}_w^\mathcal{O}=k_0^{-\chi(w^{-1}\varphi)\chi_e(\varphi(wc^\mathcal{O}))}
u_0^{\chi(w^{-1}\varphi)\chi_o(\varphi(wc^\mathcal{O}))}\\
\hphantom{{\rm coeff}_w^\mathcal{O}=}{}
+\bigl(k_0-k_0^{-1}\bigr)\chi_+\bigl(w^{-1}\varphi\bigr)\chi_e\bigl(\varphi\bigl(wc^\mathcal{O}\bigr)\bigr)+
\bigl(u_0-u_0^{-1}\bigr)\chi_-\bigl(w^{-1}\varphi\bigr)\chi_o\bigl(\varphi\bigl(wc^\mathcal{O}\bigr)\bigr).
\end{gather*}
Now note that
\begin{gather*}
k_0^{-\chi(w^{-1}\varphi)\chi_e(\varphi(wc^\mathcal{O}))}
u_0^{\chi(w^{-1}\varphi)\chi_o(\varphi(wc^\mathcal{O}))}=\bigl(k_0\chi_-\bigl(w^{-1}\varphi\bigr)+k_0^{-1}\chi_+\bigl(w^{-1}\varphi\bigr)\bigr)\chi_e \bigl(\varphi\bigl(wc^\mathcal{O}\bigr)\bigr)\\
\hphantom{k_0^{-\chi(w^{-1}\varphi)\chi_e(\varphi(wc^\mathcal{O}))}
u_0^{\chi(w^{-1}\varphi)\chi_o(\varphi(wc^\mathcal{O}))}=}{}
+\bigl(u_0\chi_+\bigl(w^{-1}\varphi\bigr)+u_0^{-1}\chi_-\bigl(w^{-1}\varphi\bigr)\bigr)\chi_o\bigl(\varphi\bigl(wc^\mathcal{O}\bigr)\bigr),
\end{gather*}
and hence
\begin{equation}\label{coefffin}
{\rm coeff}_w^\mathcal{O}=k_0\chi_e\bigl(\varphi\bigl(wc^\mathcal{O}\bigr)\bigr)+u_0\chi_o\bigl(\varphi\bigl(wc^\mathcal{O}\bigr)\bigr)=
k_0^{\chi_e(\varphi(wc^\mathcal{O}))}u_0^{\chi_o(\varphi(wc^\mathcal{O}))}.
\end{equation}
Note that \smash{$y\mapsto k_0^{\chi_e(\varphi(y))}u_0^{\chi_o(\varphi(y))}$} is $\tau\smash{\bigl(Q^\vee\bigr)}$-invariant (compare with the proof of Theorem~\ref{theoremH}),
hence~\eqref{formula0} and~\eqref{coefffin} lead to the formula
\begin{equation*}
\frac{(\mathfrak{s}_\mathcal{O}t)^\mu}{\kappa_w^\mathcal{O}}\psi_{t}^\mathcal{O}\bigl(T_0x^\lambda T_wY^\mu\bigr)=\mathcal{D}_0(x^y)
\end{equation*}
for $y=\lambda+wc^\mathcal{O}\in\mathcal{O}$ and $\mu\in Q^\vee$, where $\mathcal{D}_0$ is the linear operator on $\mathbf{F}[\mathcal{O}]$ defined by
\[
\mathcal{D}_0(x^y):=
k_0^{\chi_e(\overline{\alpha}_0(y))}u_0^{\chi_o(\overline{\alpha}_0(y))}s_{0,t}x^y+\bigl(k_0-k_0^{-1}\bigr)\nabla_0^e(x^y)+\bigl(u_0-u_0^{-1}\bigr)\nabla_0^o(x^y)
\]
for $y\in\mathcal{O}$.

In conclusion, the kernel of $\psi_t^\mathcal{O}\colon \mathbb{H}\twoheadrightarrow\mathbf{F}[\mathcal{O}]$ (see~\eqref{psitO})
is a left $\mathbb{H}$-module, and the resulting isomorphism \smash{$\mathbb{H}/\textup{ker}\bigl(\psi_t^\mathcal{O}\bigr)\overset{\sim}{\longrightarrow}\mathbf{F}[\mathcal{O}]$}
extends the quasi-polynomial $H^X$-action $\pi$ on $\mathbf{F}[\mathcal{O}]$ to an action of $\mathbb{H}$ with $T_0\in H\subset\mathbb{H}$ acting by
$\mathcal{D}_0$, hence the corresponding representation map is \smash{$\pi_t^\mathcal{O}$}. This concludes the proof.
\end{proof}

\begin{Remark}\label{CCpolHH}
For the $W$-orbit $\mathcal{O}=Q^\vee$ containing the origin, we have $\mathbf{T}_{Q^\vee}=\{1_{\mathbf{T}}\}$ since $J\smash{\bigl(Q^\vee\bigr)}=\{1,\dots,r\}$.
By Lemma~\ref{polredNabla} and by the fact that $g_{1_{\mathbf{T}}}f=g(f)$ for $g\in W$ and $f\in\mathbf{F}\smash{\bigl[Q^\vee\bigr]}$ (see Lemma~\ref{tdefaction}),
we conclude that
\[
\pi_{1_{\mathbf{T}}}^{Q^\vee}(T_j)x^\mu=k_js_j(x^\mu)+\bigl(k_j-k_j^{-1}+\bigl(u_j-u_j^{-1}\bigr)x^{\alpha_j^\vee}\bigr)\biggl(\frac{x^\mu-s_j(x^\mu)}{1-x^{2\alpha_j^\vee}}\biggr)
\]
for $j=0,\dots,r$ and $\mu\in Q^\vee$. Hence
\smash{$\pi_{1_{\mathbf{T}}}^{Q^\vee}\colon \mathbb{H}\rightarrow\textup{End}\bigl(\mathbf{F}\smash{\bigl[Q^\vee\bigr]}\bigr)$} is the basic representation of $\mathbb{H}$, due to Cherednik~\cite{ChKZ2} for $\mathbf{k}\in\mathcal{K}^{{\rm res}}$ and due to Noumi~\cite{No} and Sahi~\cite{Sa} when $\mathbf{k}\in\mathcal{K}$ and $\Phi_0$ is of type ${\rm C}_r$.
\end{Remark}
\begin{Corollary}\label{corindstep}
We have
\[
\pi_t^\mathcal{O}\bigl(x^\lambda T_wY^\mu)x^{c^\mathcal{O}}=\kappa_w^\mathcal{O}(\mathfrak{s}_\mathcal{O}t)^{-\mu}x^{\lambda+wc^\mathcal{O}}
\]
for $\lambda,\mu\in Q^\vee$ and $w\in W_0$, with $\mathfrak{s}_\mathcal{O}\in\mathbf{T}$ given by~\eqref{s} and $\kappa_w^\mathcal{O}$ given by~\eqref{kappawy}.
\end{Corollary}
\begin{proof}
In the proof of Theorem~\ref{mthm}, we showed that the epimorphism
$\psi_t^\mathcal{O}\colon \mathbb{H}\twoheadrightarrow\mathbf{F}[\mathcal{O}]$, defined by~\eqref{psitO}, is $\mathbb{H}$-linear.
Since \smash{$\psi_t^\mathcal{O}(1)=x^{c^\mathcal{O}}$}, the result now follows immediately from formula~\eqref{psitO}, cf.\ the proof of Corollary~\ref{corpsi}\,(2).
\end{proof}

\begin{Remark}
A similar remark as for the proof of Theorem~\ref{theoremH} (see Remark~\ref{H_remark}) can be made for the proof of Theorem~\ref{mthm}. With the $\kappa_w$ chosen to be~\eqref{kappawy},
the proof of Theorem~\ref{mthm} involves choosing some $\mathfrak{s}\in\mathbf{T}$ such that the coefficients ${\rm coeff}_w^\mathcal{O}$ (see~\eqref{coeff}) only depends on the coset $wW_{0,\mathcal{O}}$ for all $w\in W_0$.
For any such choice, one gets explicit realisations of quotients of cyclic $Y$-parabolically induced $\mathbb{H}$-modules with the associated induction datum given by $\mathfrak{s}t$ (cf.~Corollary~\ref{corindstep}). This forces $\mathfrak{s}$ to lie in suitable affine subtori of $\mathbf{T}$.
With the present choice $\mathfrak{s}=\mathfrak{s}_\mathcal{O}$ (see~\eqref{s}) one obtains the explicit realisations of all cyclic $Y$-parabolically induced $\mathbb{H}$-modules,
see Section~\ref{InducedModuleSection} for details.
\end{Remark}

\section[Quasi-polynomial analogs of the nonsymmetric Macdonald--Koornwinder polynomials]{Quasi-polynomial analogs of the nonsymmetric\\ Macdonald--Koornwinder polynomials}\label{qKsection}

We fix a $W$-orbit $\mathcal{O}$ in $E$ and $t\in\mathbf{T}_{\mathcal{O}}$ throughout this section.

In the first part of this section, we show that the commuting operators $\pi_{t}^\mathcal{O}(Y^\mu)$, $\mu\in Q^\vee$, on~$\mathbf{F}[\mathcal{O}]$ are triangular with respect to an appropriate partial order on the basis $\{x^y\}_{y\in\mathcal{O}}$ of quasi-monomials. This will lead to the definition of the quasi-polynomial analogs of the nonsymmetric Macdonald--Koorn\-win\-der polynomials as the simultaneous eigenfunctions of the operators~$\pi_{t}^\mathcal{O}(Y^\mu)$, $\mu\in Q^\vee$. The techniques in this section again closely follow the paper~\cite{SSV2}, in which these results are derived for $\mathbf{k}\in\mathcal{K}^{{\rm res}}$.

We first establish triangularity for a family $G_{t}^\mathcal{O}(a)$, $a\in\Phi$, of operators
closely related to the~$\pi_{t}^\mathcal{O}(T_j)$.
The linear operator $G_{t}^\mathcal{O}(a)$ on $\mathbf{F}[\mathcal{O}]$ is defined by the formula
\begin{gather*}
G_{t}^\mathcal{O}(a)x^y:=\mathbf{k}_a^{\chi_e(\overline{a}(y))}\mathbf{k}_{\frac{a}{2}}^{\chi_o(\overline{a}(y))}x^y +\bigl(\mathbf{k}_a-\mathbf{k}_a^{-1}\bigr)\Biggl(\frac{1-x^{\lfloor\overline{a}(y)\rfloor_e a^\vee}}{1-x^{-2a^\vee}}\Biggr)s_{a,t}x^y\\
\hphantom{G_{t}^\mathcal{O}(a)x^y:=}{}
+\bigl(\mathbf{k}_{\frac{a}{2}}-\mathbf{k}_{\frac{a}{2}}^{-1}\bigr)\Biggl(\frac{x^{-a^\vee}-x^{\lfloor \overline{a}(y)\rfloor_o a^\vee}}{1-x^{-2a^\vee}}\Biggr)
s_{a,t}x^y
\end{gather*}
for $y\in\mathcal{O}$.
\begin{Lemma}\label{elementaryG}
We have
\begin{enumerate}\itemsep=0pt
\item[$(1)$] $G_{t}^\mathcal{O}(\alpha_j)=s_{j,t}\pi_{t}^\mathcal{O}(T_j)$ for $j=0,\dots,r$.
\item[$(2)$] $G_{t}^\mathcal{O}(ga)=g_tG_t^\mathcal{O}(a)g_t^{-1}$ for $g\in W$.
\end{enumerate}
\end{Lemma}

\begin{proof}This follows from a direct computation using~\eqref{wcomp} and~\eqref{awcomp}, cf.~\cite[Section~5.5]{SSV2}.
\end{proof}

Recall the definition of $\mathbf{g}_y\in W^{\mathcal{O}}$ from Definition~\ref{gy}. Denote by $\leq$ the partial order on $E$ defined by
\[
y\leq z\quad\Leftrightarrow\quad y\in Wz \ \text{and}\ \mathbf{g}_y\leq_B\mathbf{g}_z
\]
with $\leq_B$ the Bruhat order of $(W,\{s_0,\dots,s_r\})$. Note that for each $z\in E$,
\[
\{y\in E \mid y\leq z\}
\]
is a finite set contained in the $W$-orbit $Wz$ of $z$. Various other properties of this partial order are obtained in~\cite[Section~5.4]{SSV2}.

\begin{Definition}
For $f\in\mathbf{F}[E]$, we write
\[
f=dx^y+ {\rm l.o.t.}
\]
if $f\in dx^y+\bigoplus_{z<y}\mathbf{F}x^z$ with $d\in\mathbf{F}^\times$. We then say that $f$ is of degree $y$ with leading coefficient~$d$.
\end{Definition}

Define $\eta_e,\eta_o\colon \mathbb{R}\rightarrow\{-1,0,1\}$ by
\begin{equation*}
\eta_e=\chi_{2\mathbb{Z}_{\geq 1}}-\chi_{2\mathbb{Z}_{\leq 0}},\qquad
\eta_o=\chi_{1+2\mathbb{Z}_{\geq 0}}-\chi_{1+2\mathbb{Z}_{<0}}.
\end{equation*}
We also set
\begin{equation*}
\eta:=\eta_e+\eta_o,
\end{equation*}
which is equal to $\chi_{\mathbb{Z}_{>0}}-\chi_{\mathbb{Z}_{\leq 0}}$. We have the following extension of
\cite[Lemma~5.27]{SSV2}.

\begin{Lemma}\label{lemG}
For $a\in\Phi_0^+\times\mathbb{Z}$, we have
\[
G_{t}^\mathcal{O}(a)x^y=\mathbf{k}_a^{-\eta_e(\overline{a}(y))}\mathbf{k}_{\frac{a}{2}}^{-\eta_o(\overline{a}(y))}x^y+ {\rm l.o.t.}
\]
for all $y\in\mathcal{O}$.
\end{Lemma}

\begin{proof}
This is covered by~\cite[Lemma~5.27]{SSV2} unless $\Phi_0$ is of type ${\rm C}_r$. For $\Phi_0$ of type ${\rm C}_r$, one checks using~\cite[Lemma~4.3\,(1)]{SSV2},~\cite[Proposition~5.20]{SSV2}
and~\cite[Lemma~5.24]{SSV2} that for $a\in\Phi_0^+\times\mathbb{Z}$,
\begin{gather*}
\Biggl(\frac{1-x^{\lfloor\overline{a}(y)\rfloor_e a^\vee}}{1-x^{-2a^\vee}}\Biggl)s_{a,t}x^y=-\chi_{2\mathbb{Z}_{\geq 1}}(\overline{a}(y))x^y+ {\rm l.o.t.},\\
\Biggl(\frac{x^{-a^\vee}-x^{\lfloor\overline{a}(y)\rfloor_o a^\vee}}{1-x^{-2a^\vee}}\Biggl)s_{a,t}x^y=-\chi_{1+2\mathbb{Z}_{\geq 0}}(\overline{a}(y))x^y+ {\rm l.o.t.},
\end{gather*}
and hence
\begin{gather*}
G_{t}^\mathcal{O}(a)x^y\\
\qquad{}=\bigl(\mathbf{k}_a^{\chi_{e}(\overline{a}(y))}\mathbf{k}_{\frac{a}{2}}^{\chi_{o}(\overline{a}(y))}
+\bigl(\mathbf{k}_a^{-1}-\mathbf{k}_a\bigr)\chi_{2\mathbb{Z}_{\geq 1}}(\overline{a}(y))+
\bigl(\mathbf{k}_{\frac{a}{2}}^{-1}-\mathbf{k}_{\frac{a}{2}}\bigr)\chi_{1+2\mathbb{Z}_{\geq 0}}(\overline{a}(y))\bigr)x^y+ {\rm l.o.t.}\\
\qquad{}=\mathbf{k}_a^{-\eta_e(\overline{a}(y))}\mathbf{k}_{\frac{a}{2}}^{-\eta_o(\overline{a}(y))}x^y+ {\rm l.o.t.},
\end{gather*}
as desired.
\end{proof}

\begin{Definition}
For $y\in E$, define $\mathfrak{s}_y\in\mathbf{T}$ by
\[
\mathfrak{s}_y:=\prod_{\alpha\in\Phi_0^+}\bigl(\mathbf{k}_\alpha\mathbf{k}_{(1,\alpha)}\bigr)^{\frac{\eta_e(\alpha(y))}{2}\alpha}
\bigl(\mathbf{k}_{\frac{\alpha}{2}}\mathbf{k}_{(\frac{1}{2},\frac{\alpha}{2})}\bigr)^{\frac{\eta_o(\alpha(y))}{2}\alpha}.
\]
\end{Definition}

In other words, the value of $\mathfrak{s}_y$ at $\lambda\in Q^\vee$ is
\begin{equation*}
\mathfrak{s}_y^\lambda:=\prod_{\alpha\in\Phi_0^+}\bigl(\mathbf{k}_\alpha\mathbf{k}_{(\alpha,1)}\bigr)^{\frac{\eta_e(\alpha(y))\alpha(\lambda)}{2}}
\bigl(\mathbf{k}_{\frac{\alpha}{2}}\mathbf{k}_{(\frac{\alpha}{2},\frac{1}{2})}\bigr)^{\frac{\eta_o(\alpha(y))\alpha(\lambda)}{2}}.
\end{equation*}
If $\alpha\smash{\bigl(Q^\vee\bigr)}=\mathbb{Z}_e$, then the factors in this product are clearly well defined.
If $\alpha\smash{\bigl(Q^\vee\bigr)}=\mathbb{Z}$, then \smash{$\mathbf{k}_\alpha=\mathbf{k}_{(\alpha,1)}=\mathbf{k}_{\frac{\alpha}{2}}=\mathbf{k}_{(\frac{\alpha}{2},\frac{1}{2})}$}, and the corresponding factor in the
product should be read~as $\mathbf{k}_\alpha^{\eta(\alpha(y))\alpha(\lambda)}$. In particular,
\[
\mathfrak{s}_y=\prod_{\alpha\in\Phi_0^+}\mathbf{k}_\alpha^{\eta(\alpha(y))\alpha}\qquad \textup{for}\ \mathbf{k}\in\mathcal{K}^{{\rm res}},
\]
which is the base-point considered in~\cite{SSV2} (see~\cite[Definition~5.1]{SSV2}). By~\cite[Lemma~2.5]{SSV2}, the function $E\rightarrow \mathbf{T}$, $y\mapsto\mathfrak{s}_y$ is constant on the faces of the affine root hyperplane arrangement.

\begin{Remark}
If $c\in\overline{C_+}$, then~\eqref{W0C} implies that
\[
\mathfrak{s}_c=\prod_{\alpha\in\Phi_0^+}\bigl(\mathbf{k}_{\alpha}\mathbf{k}_{(\alpha,1)}\bigr)^{-\frac{\chi_e(\alpha(c))}{2}\alpha}
\bigl(\mathbf{k}_{\frac{\alpha}{2}}\mathbf{k}_{(\frac{\alpha}{2},\frac{1}{2})}\bigr)^{\frac{\chi_o(\alpha(c))}{2}\alpha}.
\]
In particular,
\[
\mathfrak{s}_{c^\mathcal{O}}=\mathfrak{s}_\mathcal{O},
\]
with $\mathfrak{s}_\mathcal{O}\in\mathbf{T}$ the torus element appearing before in Corollary~\ref{corindstep} (see~\eqref{s}).
\end{Remark}

\begin{Corollary}\label{Ytriang}
For all $\mu\in Q^\vee$ and $y\in\mathcal{O}$, we have
\[
\pi_{t}^\mathcal{O}(Y^\mu)x^{y}= (\mathfrak{s}_y\mathbf{g}_yt)^{-\mu}x^y+ {\rm l.o.t.}
\]
\end{Corollary}

\begin{proof}
This is~\cite[Proposition~5.28]{SSV2} when $\mathbf{k}\in\mathcal{K}^{{\rm res}}$. Using Lemma~\ref{lemG} as replacement of~\cite[Lemma~5.27]{SSV2},
the proof of~\cite[Proposition~5.28]{SSV2} extends to the case $\mathbf{k}\in\mathcal{K}$.
\end{proof}

We now first derive some further properties of the map $E\rightarrow\mathbf{T}$, $y\mapsto\mathfrak{s}_y$.
The following lemma extends~\cite[Lemma~5.3]{SSV2}.

\begin{Lemma}\label{lemsy}
Let $W_y:=\{w\in W \mid wy=y\}$ be the subgroup of $W$ fixing $y\in E$, and let ${j\in\{0,\dots,r\}}$.
\begin{enumerate}\itemsep=0pt
\item[$(1)$] If $s_j\in W_y$, then
\[
\mathfrak{s}_y^{\overline{\alpha}_j^\vee}=\widetilde{\mathbf{k}}_{\alpha_j}^{-1}\widetilde{\mathbf{k}}_{\frac{\alpha_j}{2}}^{-1}\qquad \text{and}
\qquad s_{\overline{\alpha}_j}\mathfrak{s}_y=\bigl(\widetilde{\mathbf{k}}_{\alpha_j}\widetilde{\mathbf{k}}_{\frac{\alpha_j}{2}}\bigr)^{\overline{\alpha}_j}\mathfrak{s}_y.
\]
\item[$(2)$] If $s_j\not\in W_y$, then $s_{\overline{\alpha}_j}\mathfrak{s}_y=\mathfrak{s}_{s_jy}$.
\end{enumerate}
\end{Lemma}

\begin{proof}
We give here the required adjustments to the proof of~\cite[Lemma~5.3]{SSV2}.

For $1\leq i\leq r$, we have $\overline{\alpha_i}=\alpha_i$, $s_{\overline{\alpha}_i}=s_i$ and $\Pi(s_i)=\{\alpha_i\}$, hence
\[
s_i\mathfrak{s}_y=\bigl(\mathbf{k}_{\alpha_i},\mathbf{k}_{(\alpha_i,1)}\bigr)^{-(\eta_e(\alpha_i(y))+\eta_e(-\alpha_i(y)))\frac{\alpha_i}{2}}
\bigl(\mathbf{k}_{\frac{\alpha_i}{2}}\mathbf{k}_{(\frac{\alpha_i}{2},\frac{1}{2})}\bigr)^{-(\eta_o(\alpha_i(y))+\eta_o(-\alpha_i(y)))\frac{\alpha_i}{2}}
\mathfrak{s}_y,
\]
with the obvious interpretation of the right-hand side when $\alpha_i\smash{\bigl(Q^\vee\bigr)}=\mathbb{Z}$. Then (1) and (2) follow from the fact that
\[
\eta_e(z)+\eta_e(-z)=-2\chi_{\{0\}}(z),\qquad \eta_o(z)+\eta_o(-z)=0
\]
for $z\in\mathbb{R}$.

We now prove the lemma for $s_0$ when $\Phi_0$ is of type ${\rm C}_r$ (the other types are covered by
\cite[Lemma~5.3]{SSV2}). Write $\Phi_{0,\ell}^{\pm}$ \big(resp.\ $\Phi_{0,s}^{\pm}$\big) for the positive and negative long (resp.\ short) roots in~$\Phi_0$. Clearly,
\[
\Pi(s_\varphi)=\Pi_\ell(s_\varphi)\sqcup\Pi_s(s_\varphi)
\]
with $\Pi_\ell(w):=\Phi_{0,\ell}^+\cap w^{-1}\Phi_{0,\ell}^-$ and $\Pi_s(w):=\Phi_{0,s}^+\cap w^{-1}\Phi_{0,s}^-$ for $w\in W_0$. Furthermore,
\[
\Pi_\ell(s_\varphi)=\{\varphi\},\qquad \Pi_s(s_\varphi)=\bigl\{\alpha\in\Phi_{0,s}^+\mid \alpha\bigl(\varphi^\vee\bigr)=1\bigr\},
\]
where we used for the first equality that there are no long positive roots $\alpha$ with \smash{$\alpha\bigl(\varphi^\vee\bigr)=1$}, because
$\Phi_0$ is of type ${\rm C}_r$. Following the proof of~\cite[Lemma~5.3]{SSV2} and using that
\[
\mathbf{k}_\beta=\mathbf{k}_{(\beta,1)}=\mathbf{k}_{\frac{\beta}{2}}=\mathbf{k}_{(\frac{\beta}{2},\frac{1}{2})}
\]
for $\beta\in\Phi_{0,s}^+$, we get
\begin{gather}
s_\varphi\mathfrak{s}_y=\mathfrak{s}_{s_0y}(k_0k_r)^{-(\eta_e(\varphi(y))+\eta_e(2-\varphi(y)))\frac{\varphi}{2}}
(u_0u_r)^{-(\eta_o(\varphi(y))+\eta_o(2-\varphi(y)))\frac{\varphi}{2}} \nonumber\\
\hphantom{s_\varphi\mathfrak{s}_y=}{}
\times\prod_{\beta\in\Phi_{0,s}^+\colon \beta(\varphi^\vee)=1}
\mathbf{k}_\beta^{-(\eta(-\beta(s_\varphi y))+\eta(1+\beta(s_\varphi y)))\beta}.\label{sy0}
\end{gather}
Now the product in the second line of~\eqref{sy0} is $1_{\mathbf{T}}$ since $\eta(z)+\eta(1-z)=0$ for $z\in\mathbb{R}$. Applying the elementary formulas
\begin{equation}\label{els0}
\eta_e(z)+\eta_e(2-z)=0,\qquad \eta_o(z)+\eta_o(2-z)=2\chi_{\{1\}}(z)
\end{equation}
for $z\in\mathbb{R}$ to the first line of~\eqref{sy0}, the identity~\eqref{sy0} reduces to
\[
s_\varphi\mathfrak{s}_y=\mathfrak{s}_{s_0y}(u_0u_r)^{-\chi_{\{1\}}(\varphi(y))\varphi},
\]
from which the lemma for $j=0$ follows immediately.
\end{proof}

We denote by $\overline{w}\in W_0$ the image
of $w\in W$ under the group homomorphism $W\twoheadrightarrow W_0$, ${v\tau(\lambda)\mapsto v}$, $v\in W_0$, $\lambda\in Q^\vee$.
Note that $\overline{s_{a}}=s_{\overline{a}}$ for $a\in\Phi$.

\begin{Corollary}\label{sycor}
For $y\in\mathcal{O}$, we have
\[
\mathfrak{s}_y=\overline{\mathbf{g}_y}\mathfrak{s}_\mathcal{O}\qquad \text{and}\qquad \mathfrak{s}_y\mathbf{g}_yt=\mathbf{g}_y(\mathfrak{s}_{\mathcal{O}}t).
\]
Furthermore, \smash{$\mathfrak{s}_\mathcal{O}^{\overline{\alpha}_j^\vee}=\widetilde{\mathbf{k}}_{\alpha_j}^{-1}\widetilde{\mathbf{k}}_{\frac{\alpha_j}{2}}^{-1}$} for all $j\in J(\mathcal{O})$.
\end{Corollary}

\begin{proof}
Similar to the proof of~\cite[Proposition~5.4]{SSV2} and~\cite[Corollary~5.5]{SSV2}.
\end{proof}

By Corollary~\ref{sycor}, we have $\mathfrak{s}_\mathcal{O}\mathbf{T}_\mathcal{O}=\mathbf{L}_\mathcal{O}$
with
\begin{equation}\label{LO}
\mathbf{L}_\mathcal{O}:=
\bigl\{\gamma\in\mathbf{T}\mid \gamma^{\alpha_j^\vee}=\widetilde{\mathbf{k}}_{\alpha_j}^{-1}\widetilde{\mathbf{k}}_{\frac{\alpha_j}{2}}^{-1}\quad \forall j\in
J(\mathcal{O})\bigr\}.
\end{equation}
Note here that \smash{$\widetilde{\mathbf{k}}_{\alpha_0}\widetilde{\mathbf{k}}_{\frac{\alpha_0}{2}}=u_0u_r$}, \smash{$\widetilde{\mathbf{k}}_{\alpha_i}\widetilde{\mathbf{k}}_{\frac{\alpha_i}{2}}=k_i^2=k^2$} for $1\leq i<r$, and \smash{$\widetilde{\mathbf{k}}_{\alpha_r}\widetilde{\mathbf{k}}_{\frac{\alpha_r}{2}}=
k_0k_r$}.

Write
\[
\mathbf{T}_\mathcal{O}^\prime:=\bigl\{t\in\mathbf{T}_\mathcal{O} \mid \textup{the map}\ W^\mathcal{O}\rightarrow\mathbf{T},\ g\mapsto g(\mathfrak{s}_\mathcal{O}t)\ \text{is injective}\bigr\}.
\]
Then $\mathbf{T}_{\mathcal{O}}^\prime\not=\varnothing$ for generic $q\in\mathbf{F}^\times$ and $\mathbf{k}\in\mathcal{K}$.
We are now in the position to extend the definition of the quasi-polynomial analogs of the monic nonsymmetric Macdonald polynomials, introduced in~\cite[Theorem~6.2]{SSV2}
for multiplicity functions $\mathbf{k}\in\mathcal{K}^{{\rm res}}$, to multiplicity functions
$\mathbf{k}\in\mathcal{K}$.

\begin{Theorem}\label{EJycor}
For $t\in\mathbf{T}_\mathcal{O}^\prime$ and $y\in\mathcal{O}$, there exists a unique quasi-polynomial
\[
E^\mathcal{O}_y(x)=E_y^\mathcal{O}(x;\mathbf{k},t;q)\in\mathbf{F}[\mathcal{O}]
\]
satisfying the following two properties:
\begin{enumerate}\itemsep=0pt
\item[$(a)$] $E_y^\mathcal{O}(x)=x^y+ {\rm l.o.t.}$,
\item[$(b)$] $\pi_{t}^\mathcal{O}(Y^\mu)E_y^\mathcal{O}(x)=(\mathbf{g}_y(\mathfrak{s}_\mathcal{O}t))^{-\mu}E_y^\mathcal{O}(x)$ for all $\mu\in Q^\vee$.
\end{enumerate}
\end{Theorem}

\begin{proof}
This is a direct consequence of Corollaries~\ref{Ytriang} and~\ref{sycor}.
\end{proof}

Only the Koornwinder/${\rm C}^\vee{\rm C}_r$-case of Theorem~\ref{EJycor} is new compared to~\cite[Theorem~6.2]{SSV2}. In this case, $\Phi_0$ is of type ${\rm C}_r$, $r\geq 1$, and $E_y^\mathcal{O}(x)$ depends on five multiplicity parameters ${k_0,u_0,k_r,u_r,k\in\mathbf{F}^\times}$
(four in case of $r=1$), on the dilation parameter $q$, and on the representation parameter $t\in\mathbf{T}_\mathcal{O}$.

To see how Sahi's~\cite{Sa} monic nonsymmetric Koornwinder polynomials fit into this picture, consider the special case that $\mathcal{O}=Q^\vee$. Then
$J\smash{\bigl(Q^\vee\bigr)}=\{1,\dots,r\}$ and $\mathbf{T}_{Q^\vee}=\{1_{\mathbf{T}}\}$. Then Theorem~\ref{EJycor} requires that \smash{$1_{\mathbf{T}}\in \mathbf{T}_{\{1,\dots,r\}}^\prime$}, which amounts to generic conditions on $q\in\mathbf{F}^\times$ and~${\mathbf{k}\in\mathcal{K}}$ (including, typically, the condition that $q$ is not a root of unity). By Remark~\ref{CCpolHH}, the resulting Laurent polynomial
\[
E_\mu^{Q^\vee}(x;\mathbf{k},1_{\mathbf{T}};q)\in\mathbf{F}\smash{\bigl[Q^\vee\bigr]}
\]
is Sahi's~\cite[Theorem~6.2]{Sa} monic nonsymmetric Koornwinder polynomial $E_\mu$ of degree $\mu\in Q^\vee$, with $n$ and the multiplicity parameters $t_0$, $u_0$, $t_n$, $u_n$, $t_i$, $i\not=0,n$, in~\cite{Sa} corresponding to $r$ and~$k_0$, $u_0$, $k_r$, $u_r$, $k$.

Various properties of the quasi-polynomial generalisations of the Macdonald polynomials obtained in~\cite{SSV2} have direct analogs in the Koornwinder case, such as
the face limit transitions~\cite[Proposition~6.15]{SSV2}, the creation formulas~\cite[Theorem~6.12]{SSV2} in terms of double affine Hecke algebra $Y$-intertwiners, the orthogonality relations~\cite[Theorem~6.42]{SSV2}, and (anti)symmetrisation~\cite[Section~6.6]{SSV2}. We do not give the details here. The quasi-polynomial generalisations of the symmetric Macdonald--Koornwinder polynomials will be the topic of an upcoming paper.

\section[The quasi-polynomial representation as Y-parabolically induced module]{The quasi-polynomial representation\\ as $\boldsymbol{Y}$-parabolically induced module}\label{InducedModuleSection}

In this section, $\mathcal{O}$ is a $W$-orbit in $E$ and $t\in\mathbf{T}_\mathcal{O}$. Then $\mathfrak{s}_\mathcal{O}t\in\mathfrak{s}_\mathcal{O}\mathbf{T}_{\mathcal{O}}=\mathbf{L}_\mathcal{O}$ with
$\mathbf{L}_\mathcal{O}$ the affine subtorus of $\mathbf{T}$ defined by~\eqref{LO} and $\mathfrak{s}_\mathcal{O}$ given by~\eqref{s}. Recall the definition of the subset $I(\mathcal{O})\subseteq\{1,\dots,r\}$ from Definition~\ref{IJ}.

\begin{Lemma}\label{cyclicaction}
The $\mathbb{H}$-module $\bigl(\mathbf{F}[\mathcal{O}],\pi_{t}^\mathcal{O}\bigr)$ is cyclic with
cyclic vector \smash{$x^{c^\mathcal{O}}$}. Furthermore,
\begin{gather}
\pi_{t}^\mathcal{O}(T_i)x^{c^\mathcal{O}}=k_ix^{c^\mathcal{O}},\qquad i\in I(\mathcal{O}),\nonumber\\
\pi_{t}^\mathcal{O}(Y^\mu)x^{c^\mathcal{O}}=(\mathfrak{s}_\mathcal{O}t)^{-\mu}x^{c^\mathcal{O}},\qquad \mu\in Q^\vee.\label{genfor0}
\end{gather}
\end{Lemma}
\begin{proof}
This follows immediately from Corollary~\ref{corindstep} and the fact that $\kappa_{s_i}^\mathcal{O}=k_i$ for $i\in I(\mathcal{O})$.
\end{proof}

Lemma~\ref{cyclicaction} prompts the following definition.

\begin{Definition}
We write $H_{0,\mathcal{O}}[Y]$ for the subalgebra of the
affine Hecke algebra $H=H(\mathbf{k})$ generated by $H_{0,\mathcal{O}}$ and $\mathbf{F}_Y\smash{\bigl[Q^\vee\bigr]}$.
\end{Definition}
By Lemma
\ref{cyclicaction}, \smash{$\mathbf{F}x^{c^\mathcal{O}}$} is a one-dimensional $H_{0,\mathcal{O}}[Y]$-submodule of
\smash{$\bigl(\mathbf{F}[\mathcal{O}],\pi_t^\mathcal{O}\vert_{H_{0,\mathcal{O}}[Y]}\bigr)$}, with the action defined by
the unique linear character $\zeta_{t}^\mathcal{O}\colon H_{0,\mathcal{O}}[Y]\rightarrow\mathbf{F}$ satisfying
\begin{gather}
\zeta_{t}^\mathcal{O}(T_i)=k_{i}\qquad \textup{for}\ i\in I(\mathcal{O}),\nonumber\\
\zeta_{t}^\mathcal{O}(Y^\mu)=(\mathfrak{s}_\mathcal{O}t)^{-\mu}\qquad\textup{for}\ \mu\in Q^\vee.\label{zzz}
\end{gather}
The existence of the linear character $\zeta_t^\mathcal{O}\colon H_{0,\mathcal{O}}[Y]\rightarrow\mathbf{F}$ can also be established without reference to Lemma~\ref{cyclicaction} using the Bernstein presentation of $H_{0,\mathcal{O}}[Y]\subseteq H$ in terms of the algebraic generators $T_i$, $i\in I(\mathcal{O})$, and $Y^\mu$, $\mu\in Q^\vee$. From now on, we write $\mathbf{F}1_{t}^\mathcal{O}$ for the one-dimensional $H_{0,\mathcal{O}}[Y]$-representation with representation map
$\zeta_t^\mathcal{O}$.

\begin{Proposition}\label{indform}
We have a unique surjection of $\mathbb{H}$-modules
\begin{equation*}
{\rm Ind}_{H_{0,\mathcal{O}}[Y]}^{\mathbb{H}}\bigl(\mathbf{F}1_t^\mathcal{O}\bigr)\twoheadrightarrow \bigl(\mathbf{F}[\mathcal{O}],\pi_t^\mathcal{O}\bigr)
\end{equation*}
mapping \smash{$1\otimes_{H_{0,\mathcal{O}}[Y]}1_t^\mathcal{O}$} to $x^{c^{\mathcal{O}}}$. It is an isomorphism when $\alpha_0\bigl(c^\mathcal{O}\bigr)\not=0$ $($i.e., when $I(\mathcal{O})=J(\mathcal{O}))$.
\end{Proposition}
\begin{proof}
The first statement is immediate from Lemma~\ref{cyclicaction}.

In view of the PBW theorem for $\mathbb{H}$,
\[
\bigl\{x^\lambda T_v\otimes_{H_{0,\mathcal{O}}[Y]}1_t^\mathcal{O}
\mid \lambda\in Q^\vee,\, v\in W_0^\mathcal{O}\bigr\}
\]
is a $\mathbf{F}$-basis of \smash{${\rm Ind}_{H_{0,\mathcal{O}}[Y]}^{\mathbb{H}}\bigl(\mathbf{F}1_t^\mathcal{O}\bigr)$}. This
is mapped to
\[
\bigl\{\kappa_v^\mathcal{O}x^{\lambda+vc^\mathcal{O}} \mid \lambda\in Q^\vee,\, v\in W_0^\mathcal{O}\bigr\}
\]
by Corollary~\ref{corindstep}, which is
a basis of $\mathbf{F}[\mathcal{O}]$ when $\alpha_0\bigl(c^\mathcal{O}\bigr)\not=0$
by Corollary~\ref{Pcpoldec}.
\end{proof}

\begin{Remark}
For an associative $\mathbf{F}$-algebra $A$, denote by ${\rm mod}_A$ the category of left $A$-modules. The image of
\smash{${\rm Ind}_{H_{0,\mathcal{O}}[Y]}^{\mathbb{H}}\bigl(\mathbf{F}1_t^\mathcal{O}\bigr)$} under the
restriction functor
\[
\textup{Res}_{H^X}^{\mathbb{H}}\colon\ {\rm mod}_{\mathbb{H}}\rightarrow{\rm mod}_{H^X}
\]
is isomorphic to \smash{${\rm Ind}_{H_{0,\mathcal{O}}}^{H^X}\bigl(\mathbf{F}1^\mathcal{O}\bigr)$} because \smash{${\rm Ind}_{H_{0,\mathcal{O}}[Y]}^{\mathbb{H}}\bigl(\mathbf{F}1_t^\mathcal{O}\bigr)$} is already generated by $\mathbf{1}_t^\mathcal{O}$ as
a $H^X$-module and
 $\zeta_t^\mathcal{O}\vert_{H_{0,\mathcal{O}}}$ is the trivial linear character of $H_{0,\mathcal{O}}$. Hence Proposition~\ref{indcor} follows from
 Proposition~\ref{indform} by applying the restriction functor $\textup{Res}_{H^X}^{\mathbb{H}}$.
\end{Remark}

We finish this section by realising \smash{$\bigl(\mathbf{F}[\mathcal{O}],\pi_t^\mathcal{O}\bigr)$} as a $Y$-parabolically induced $\mathbb{H}$-module when $\alpha_0\bigl(c^\mathcal{O}\bigr)=0$.
For $y\in E$ and $w\in W$, set
\[
k_w(y):=\prod_{\alpha\in\Phi_0^+}\mathbf{k}_\alpha^{\frac{\eta_e(\alpha(wy))-\eta_e(\alpha(y))}{2}}
\mathbf{k}_{\frac{\alpha}{2}}^{\frac{\eta_o(\alpha(wy))-\eta_o(\alpha(y))}{2}},
\]
which is well defined since the product involves integer powers of the multiplicity parameters. Indeed, this follows from the observation that
\begin{equation}\label{kcocycle}
k_{ww^\prime}(y)=k_w(w^\prime y)k_{w^\prime}(y)
\end{equation}
for $w,w^\prime\in W$ and the formulas
\begin{equation}\label{kcocycle1}
k_{s_i}(y)=
\begin{cases}
\mathbf{k}_{\alpha_i}^{-\eta_e(\alpha_i(y))}\mathbf{k}_{\frac{\alpha_i}{2}}^{-\eta_o(\alpha_i(y))} &\text{if}\ \alpha_i(y)\not=0,\\
1 &\text{if}\ \alpha_i(y)=0
\end{cases}
\end{equation}
for $i=1,\dots,r$ and
\begin{equation}\label{kcocycle2}
k_{s_0}(y)=
\begin{cases}
\prod_{\alpha\in\Pi(s_\varphi)}\mathbf{k}_{\alpha}^{-\eta_e(\alpha(y))}\mathbf{k}_{\frac{\alpha}{2}}^{-\eta_o(\alpha(y))} &\text{if}\ \alpha_0(y)\not=0,\\[1.1mm]
\prod_{\alpha\in\Pi(s_\varphi)\setminus\{\varphi\}}\mathbf{k}_{\alpha}^{-\eta_e(\alpha(y))}\mathbf{k}_{\frac{\alpha}{2}}^{-\eta_o(\alpha(y))} &\text{if}\ \alpha_0(y)=0,
\end{cases}
\end{equation}
which in turn follow by a computation in the spirit of Lemma~\ref{lemsy}.

The following result extends Lemma~\ref{cyclicaction} (see~\cite[Lemma~5.11]{SSV2} and~\cite[Proposition~5.29]{SSV2}). Recall the definition of the duality anti-algebra isomorphism
$\delta=\delta_{\mathbf{k}}\colon \mathbb{H}\rightarrow\widetilde{\mathbb{H}}$ with inverse $\widetilde{\delta}=\delta_{\widetilde{\mathbf{k}}}$ from Section~\ref{2.5}.

\begin{Proposition}\label{triangdelta}
We have
\begin{gather}
\pi_{t}^\mathcal{O}\bigl(\widetilde{\delta}(T_j)\bigr)x^{c^\mathcal{O}}=\widetilde{\mathbf{k}}_{\alpha_j}x^{c^\mathcal{O}} \qquad \textup{for}\ j\in J(\mathcal{O}),\nonumber\\
\pi_{t}^\mathcal{O}\bigl(\widetilde{\delta}(T_{w^{-1}})\bigr)x^{c^\mathcal{O}}=k_w\bigl(c^\mathcal{O}\bigr)x^{wc^\mathcal{O}}+ {\rm l.o.t.}\qquad \textup{for}\ w\in W^\mathcal{O}.
\label{genfor}
\end{gather}
\end{Proposition}

Note that for $j\in I(\mathcal{O})$ we have \smash{$\widetilde{\delta}(T_j)=T_j$}
and \smash{$\widetilde{\mathbf{k}}_{\alpha_j}=k_j$}, so the first line of~\eqref{genfor} is consistent with the first line of~\eqref{genfor0}.
Before proving Proposition~\ref{triangdelta}, let me explain how it leads to the interpretation of $\bigl(\mathbf{F}[\mathcal{O}],\pi_t^\mathcal{O}\bigr)$ as a $Y$-parabolically induced $\mathbb{H}$-module.

Consider the following algebras:
\begin{itemize}\itemsep=0pt
\item the subalgebra $\widetilde{H}_\mathcal{O}$ of $\widetilde{H}=H(\widetilde{\mathbf{k}})$, generated by $T_j$, $j\in J(\mathcal{O})$,
\item the subalgebra $H_\mathcal{O}^\delta:=\widetilde{\delta}(\widetilde{H}_\mathcal{O})$ of $H^X$,
\item the subalgebra $\widetilde{H}_\mathcal{O}[X]$ of $\widetilde{\mathbb{H}}$, generated by $\widetilde{H}_{\mathcal{O}}$ and $\mathbf{F}\smash{\bigl[Q^\vee\bigr]}$,
\item the subalgebra $H_\mathcal{O}^\delta[Y]:=\widetilde{\delta}\bigl(\widetilde{H}_\mathcal{O}[X]\bigr)$ of $\mathbb{H}$.
\end{itemize}
Note that
\[
H_{\mathcal{O}}^\delta=H_{0,\mathcal{O}}\qquad \textup{and}\qquad H_{\mathcal{O}}^\delta[Y]=H_{0,\mathcal{O}}[Y]\qquad \textup{when}\ \alpha_0\bigl(c^\mathcal{O}\bigr)\not=0.
\]
On the other hand, if $\alpha_0\bigl(c^\mathcal{O}\bigr)=0$,
then $H_{\mathcal{O}}^\delta[Y]$ is generated as algebra by $H_{0,\mathcal{O}}[Y]$ and
\begin{equation}\label{deltaT00}
\widetilde{\delta}(T_0)=Y^{-\varphi^\vee}T_0x^{-\varphi^\vee}
\end{equation}
\big(the equality in~\eqref{deltaT00} follows from the fact that $Y^{\varphi^\vee}=T_0T_{s_\varphi}$\big). In this case, $H_\mathcal{O}^\delta[Y]$ no longer is a subalgebra of $H$.

Recall the linear character $\zeta_t^\mathcal{O}\colon
H_{0,\mathcal{O}}[Y]\rightarrow\mathbf{F}$, defined by~\eqref{zzz}. By Lemma~\ref{cyclicaction} and Proposition~\ref{triangdelta}, we can define the following extension of
$\zeta_t^\mathcal{O}$ to a linear character
of $H_\mathcal{O}^\delta[Y]$, which we again denote by $\zeta_t^\mathcal{O}$.

\begin{Definition}We write $\zeta_t^\mathcal{O}\colon H_{\mathcal{O}}^\delta[Y]\rightarrow\mathbf{F}$ for the unique linear character of $H_{\mathcal{O}}^\delta[Y]$
satisfying
\begin{gather*}
\zeta_t^\mathcal{O}\bigl(\widetilde{\delta}(T_j)\bigr)=\widetilde{k}_{\alpha_j}\qquad \textup{for}\ j\in J(\mathcal{O}),\\
\zeta_t^\mathcal{O} (Y^\mu )= (\mathfrak{s}_\mathcal{O}t )^{-\mu}\qquad \textup{for}\ \mu\in Q^\vee.
\end{gather*}
\end{Definition}

So if $\alpha_0\bigl(c^\mathcal{O}\bigr)=0$, then $\zeta_t^\mathcal{O}$ is characterised by~\eqref{zzz} and the formula
\[
\zeta_{t}^\mathcal{O}\bigl(\widetilde{\delta}(T_0)\bigr)=\widetilde{\mathbf{k}}_{\alpha_0}=u_r.
\]
The existence of the linear character \smash{$\zeta_t^\mathcal{O}\colon H_\mathcal{O}^\delta[Y]\!\rightarrow\!\mathbf{F}$} can be proven without referring to the~quasi-polynomial representation, but by using instead that \smash{$H_\mathcal{O}^\delta[Y]=\widetilde{\delta}\bigl(\widetilde{H}_\mathcal{O}[X]\bigr)$} and the Bernstein-type presentation of \smash{$\widetilde{H}_\mathcal{O}[X]$} involving the cross relations~\eqref{crossX} for $j\in J(\mathcal{O})$ and $\mu\in Q^\vee$ with~$\mathbf{k}$ replaced by $\widetilde{\mathbf{k}}$.
 It is discussed in detail in~\cite[Section~3.2]{SSV2} when $\mathbf{k}\in\mathcal{K}^{{\rm res}}$, in which case the duality anti-algebra involution does not affect the multiplicity parameters.

\begin{Definition}
We write $\mathbf{F}1_{t}^\mathcal{O}$ for the one-dimensional $H_\mathcal{O}^\delta[Y]$-module with
representation map~$\zeta_{t}^\mathcal{O}$, and
\[
\mathbb{M}_{t}^\mathcal{O}:={\rm Ind}_{H_\mathcal{O}^\delta[Y]}^{\mathbb{H}}\bigl(\mathbf{F}1_{t}^\mathcal{O}\bigr)
\]
for the resulting induced $\mathbb{H}$-module.
\end{Definition}
Note that $\mathbb{M}_t^\mathcal{O}={\rm Ind}_{H_{0,\mathcal{O}}[Y]}^{\mathbb{H}}\bigl(\mathbf{F}1_t^\mathcal{O}\bigr)$ when $\alpha_0\bigl(c^\mathcal{O}\bigr)\not=0$, which is the induced $\mathbb{H}$-module appearing in Proposition~\ref{indform}.

We denote the canonical cyclic vector of $\mathbb{M}_t^\mathcal{O}$ by
\[
\mathbbm{1}_{t}^\mathcal{O}:=1\otimes_{H_\mathcal{O}^\delta[Y]}1_{t}^\mathcal{O}.
\]
The following theorem extends \cite[Theorem~4.5\,(2)]{SSV2} to multiplicity parameters in $\mathcal{K}$.

\begin{Theorem}
With the above notations and conventions, we have a unique isomorphism
\[
\mathbb{M}_{t}^\mathcal{O}\overset{\sim}{\longrightarrow} \bigl(\mathbf{F}[\mathcal{O}],\pi_{t}^\mathcal{O}\bigr)
\]
of $\mathbb{H}$-modules mapping
$\mathbbm{1}_{t}^\mathcal{O}$ to $x^{c^\mathcal{O}}$.
\end{Theorem}

\begin{proof}By Lemma~\ref{cyclicaction} and the first line of~\eqref{genfor}, we have a unique epimorphism
\begin{equation}\label{epi2}
\mathbb{M}_{t}^\mathcal{O}\twoheadrightarrow \bigl(\mathbf{F}[\mathcal{O}],\pi_{t}^\mathcal{O}\bigr)
\end{equation}
of $\mathbb{H}$-modules mapping $\mathbbm{1}_{t}^\mathcal{O}$ to $x^{c^\mathcal{O}}$. By the PBW theorem for $\mathbb{H}$,
\[
\bigl\{x^\mu T_uT_{w^{-1}}\mid \mu\in Q^\vee,\, u\in W_\mathcal{O},\, w\in W^\mathcal{O} \bigr\}
\]
is a basis of $\mathbb{H}$, and hence \smash{$\bigl\{\widetilde{\delta}(T_{w^{-1}})H_\mathcal{O}^\delta[Y]\mid w\in W^\mathcal{O}\bigr\}$} is a $\mathbf{F}$-basis of $\mathbb{H}/H_\mathcal{O}^\delta[Y]$. The resulting $\mathbf{F}$-basis
\smash{$\bigl\{\widetilde{\delta}(T_{w^{-1}})\mathbbm{1}_{t}^\mathcal{O}\mid w\in W^\mathcal{O}\bigr\}$} of $\mathbb{M}_{t}^\mathcal{O}$ is mapped by the epimorphism~\eqref{epi2} to
\[
\bigl\{\pi_{t}^\mathcal{O}\bigl(\widetilde{\delta}(T_{w^{-1}})\bigr)x^{c^\mathcal{O}}\mid w\in W^\mathcal{O}\bigr\},
\]
which is a basis of $\mathbf{F}[\mathcal{O}]$ due to the second line of~\eqref{genfor}. Hence the map~\eqref{epi2} is an isomorphism.\looseness=1
\end{proof}

\begin{proof}[Proof of Proposition~\ref{triangdelta}]
{\it First line of~\eqref{genfor}.} By Lemma~\ref{cyclicaction}, it suffices to check it for $j=0$ when $\alpha_0\bigl(c^\mathcal{O}\bigr)=0$,
which we assume from now on.

By~\eqref{deltaT00} and Lemma~\ref{cyclicaction}, we have
\begin{equation}\label{start0}
\pi_{t}^\mathcal{O}\bigl(\widetilde{\delta}\bigl(T_0^{-1}\bigr)\bigr)x^{c^\mathcal{O}} =(\mathfrak{s}_\mathcal{O}t)^{-\varphi^\vee}x^{\varphi^\vee}\pi_{t}^\mathcal{O}\bigl(T_0^{-1}\bigr)x^{c^\mathcal{O}}.
\end{equation}
By the explicit expression~\eqref{qprepHH} of $\pi_{t}^\mathcal{O}(T_0)$, we have, since $\overline{\alpha}_0\bigl(c^\mathcal{O}\bigr)=-1$,
\[
\pi_{t}^\mathcal{O}(T_0)x^{c^\mathcal{O}}=u_0s_{0,t}x^{c^\mathcal{O}}+\bigl(k_0-k_0^{-1}\bigr)x^{c^\mathcal{O}},
\]
and hence
\begin{equation*}
\pi_{t}^\mathcal{O}\bigl(T_0^{-1}\bigr)x^{c^\mathcal{O}}=\pi_{t}^\mathcal{O}\bigl(T_0-k_0+k_0^{-1}\bigr)x^{c^\mathcal{O}}=u_0s_{0,t}x^{c^\mathcal{O}}.
\end{equation*}
By~\eqref{ssplit} and~\eqref{taction}, we then have
\[
\pi_{t}^\mathcal{O}\bigl(T_0^{-1}\bigr)x^{c^\mathcal{O}}=u_0t^{\varphi^\vee}x^{s_\varphi c^\mathcal{O}}=u_0t^{\varphi^\vee}x^{c^\mathcal{O}-\varphi^\vee},
\]
where we used that $c^\mathcal{O}=s_0c^\mathcal{O}=s_\varphi c^\mathcal{O}+\varphi^\vee$ for the second equation.
Returning to~\eqref{start0}, we conclude that
\[
\pi_{t}^\mathcal{O}\bigl(\widetilde{\delta}\bigl(T_0^{-1}\bigr)\bigr)x^{c^\mathcal{O}}=
\mathfrak{s}_\mathcal{O}^{-\varphi^\vee}u_0x^{c^\mathcal{O}}.
\]
But $0\in J(\mathcal{O})$ by the assumption that $\alpha_0\bigl(c^\mathcal{O}\bigr)=0$, so
\[
\mathfrak{s}_\mathcal{O}^{-\varphi^\vee}=\mathfrak{s}_\mathcal{O}^{\overline{\alpha}_0}=u_0^{-1}u_r^{-1}
\]
by Corollary~\ref{sycor}, and we conclude that
\[
\pi_{t}^\mathcal{O}\bigl(\widetilde{\delta}(T_0)\bigr)x^{c^\mathcal{O}}=u_rx^{c^\mathcal{O}}=\widetilde{\mathbf{k}}_{\alpha_0}x^{c^\mathcal{O}},
\]
as desired.

{\it Second line of~\eqref{genfor}.} The proof uses the following lemma.

\begin{Lemma}\label{twtinvlemma}
For $j\in\{0,\dots,r\}$ and $y\in\mathcal{O}$ with $\alpha_j(y)>0$, we have
\begin{equation*}
\pi_{t}^\mathcal{O}\bigl(\widetilde{\delta}(T_j)\bigr)x^y=k_{s_j}(y)x^{s_jy}+ {\rm l.o.t.}
\end{equation*}
\end{Lemma}

\begin{proof}
The proof we give here deviates from the proof of~\cite[Proposition~5.29]{SSV2}. We will make use of Lemma~\ref{lemG}, which simplifies the computations.

(1) Consider first the case that $j=i\in\{1,\dots,r\}$. By~\eqref{kcocycle1}, we then have to show that
 \begin{equation}\label{twtinv}
\pi_{t}^\mathcal{O}(T_i)x^y=\mathbf{k}_{\alpha_i}^{-\eta_e(\alpha_i(y))}\mathbf{k}_{\frac{\alpha_i}{2}}^{-\eta_0(\alpha_i(y))}x^{s_iy}+ {\rm l.o.t.}
\end{equation}
for $y\in\mathcal{O}$ satisfying $\alpha_i(y)>0$.

For any $y\in\mathcal{O}$, we have
\begin{equation}\label{twt}
\pi_{t}^\mathcal{O}\bigl(T_i^{-1}\bigr)x^y=G_{t}^\mathcal{O}(\alpha_i)^{-1}x^{s_iy}=\mathbf{k}_{\alpha_i}^{\eta_e(-\alpha_i(y))}\mathbf{k}_{\frac{\alpha_i}{2}}^{\eta_0(-\alpha_i(y))}x^{s_iy}+ {\rm l.o.t.}
\end{equation}
by Lemmas~\ref{elementaryG}\,(1) and~\ref{lemG}. If in addition $\alpha_i(y)>0$,
then $y<s_iy$ by~\cite[Proposition~5.21]{SSV2} and a direct computation shows that
$\eta_e(-\alpha_i(y))=-\eta_e(\alpha_i(y))$ and $\eta_o(-\alpha_i(y))=-\eta_o(\alpha_i(y))$.
Formula~\eqref{twtinv} for $y\in\mathcal{O}$ satisfying $\alpha_i(y)>0$ then follows from~\eqref{twt} and the fact that
$T_i=T_i^{-1}+k_i-k_i^{-1}$.

(2) Consider now the case that $j=0$. By~\eqref{kcocycle2}, we then have to show that
\begin{equation}\label{twt2inv}
\pi_{t}^\mathcal{O}\bigl(\widetilde{\delta}(T_0)\bigr)x^y=
\biggl(\prod_{\alpha\in\Pi(s_\varphi)}\mathbf{k}_\alpha^{-\eta_e(\alpha(y))}\mathbf{k}_{\frac{\alpha}{2}}^{-\eta_0(\alpha(y))}\biggr)x^{s_0y}+ {\rm l.o.t.}
\end{equation}
for $y\in\mathcal{O}$ satisfying $\alpha_0(y)>0$.

Consider the element
\[
U_0:=x^{-\alpha_0^\vee}T_0^{-1}=q_\varphi^{-1}x^{\varphi^\vee}T_{s_\varphi}Y^{-\varphi^\vee}\in\mathbb{H}
\]
(the second equality follows from the fact that \smash{$x^{\alpha_0^\vee}=q_{\varphi}x^{-\varphi^\vee}$} and from formula~\eqref{Hid} for~${w=1}$). For type ${\rm C}_r$, the element $U_0$ was introduced by Sahi~\cite{Sa}, who in particular showed that
$U_0$ satisfies the Hecke relation $(U_0-u_0)\bigl(U_0+u_0^{-1}\bigr)=0$ (but we are not going to need this here). By formula~\eqref{Hid} with $w=1$, we have
\begin{equation}\label{deltaT0}
\widetilde{\delta}(T_0)=T_{s_\varphi}^{-1}x^{-\varphi^\vee}=q_\varphi^{-1}Y^{-\varphi^\vee}U_0^{-1}.
\end{equation}
So it suffices to focus on the quasi-monomial expansion of $\pi_{t}^\mathcal{O}\bigl(U_0^{-1}\bigr)x^y$ and then use Corollary~\ref{Ytriang}.

For the moment, suppose that $y\in\mathcal{O}$ is arbitrary. We compute, using Lemma~\ref{elementaryG},
\begin{equation}\label{sss1}
\pi_{t}^\mathcal{O}\bigl(U_0^{-1}\bigr)x^y= \pi_{t}^\mathcal{O}(T_0)\bigl(x^{y+\alpha_0^\vee}\bigr)=
G_{t}^\mathcal{O}(-\alpha_0)s_{0,t}\bigl(x^{y+\alpha_0^\vee}\bigr).
\end{equation}
By~\eqref{taction} and~\eqref{wcomp}, we have
\[
s_{0,t}\bigl(x^{y+\alpha_0^\vee}\bigr)=s_0\bigl(x^{\alpha_0^\vee}\bigr)s_{0,t}(x^y)=x^{-\alpha_0^\vee}(\mathbf{g}_yt)^{\varphi^\vee}x^{s_\varphi y}=
q_\varphi^{-1}(\mathbf{g}_yt)^{\varphi^\vee}x^{s_0y}.
\]
Substituting in~\eqref{sss1} then gives
\[
\pi_{t}^\mathcal{O}\bigl(U_0^{-1}\bigr)x^y=q_\varphi^{-1}(\mathbf{g}_yt)^{\varphi^\vee}G_{t}^\mathcal{O}(-\alpha_0)x^{s_0y}.
\]
Now $-\alpha_0\in\Phi_0^+\times\mathbb{Z}$, so by Lemma~\ref{lemG},
\begin{equation*}
\pi_{t}^\mathcal{O}\bigl(U_0^{-1}\bigr)x^y=q_\varphi^{-1}(\mathbf{g}_yt)^{\varphi^\vee}k_0^{-\eta_e(\varphi(s_0y))}u_0^{-\eta_0(\varphi(s_0y))}x^{s_0y}+ {\rm l.o.t.}
\end{equation*}
Combined with~\eqref{deltaT0} and Corollary~\ref{Ytriang}, we conclude that
\begin{equation}\label{almostthere}
\pi_{t}^\mathcal{O}\bigl(\widetilde{\delta}(T_0)\bigr)x^y=q_\varphi^{-2}\mathfrak{s}_{s_0y}^{\varphi^\vee}(\mathbf{g}_{s_0y}t)^{\varphi^\vee}(\mathbf{g}_yt)^{\varphi^\vee}
k_0^{-\eta_e(\varphi(s_0y))}u_0^{-\eta_0(\varphi(s_0y))}x^{s_0y}+ {\rm l.o.t.}
\end{equation}

From now on, we assume that $\alpha_0(y)>0$. Then $s_\varphi\mathfrak{s}_y=\mathfrak{s}_{s_0y}$ by Lemma~\ref{lemsy}\,(2), hence
\smash{$\mathfrak{s}_{s_0y}^{\varphi^\vee}=\mathfrak{s}_y^{-\varphi^\vee}$}. Furthermore,
\[
\mathbf{g}_{s_0y}t=s_0\mathbf{g}_yt=q^{\varphi^\vee}s_\varphi\mathbf{g}_yt
\]
in $\mathbf{T}$, where we used in the first equality that $t\in\mathbf{T}_{\mathcal{O}}\subseteq\mathbf{T}^{W_{\mathcal{O}}}$, hence we may replace $\mathbf{g}_{s_0y}$
with any other affine Weyl group element mapping $c^\mathcal{O}$ to $s_0y$. Hence \smash{$(\mathbf{g}_{s_0y}t)^{\varphi^\vee}=q_\varphi^2(\mathbf{g}_yt)^{-\varphi^\vee}$}. So the leading coefficient in~\eqref{almostthere} reduces to
\[
\mathfrak{s}_y^{-\varphi^\vee}k_0^{-\eta_e(\varphi(s_0y))}u_0^{-\eta_0(\varphi(s_0y))}.
\]
Note that $\varphi(s_0y)=2-\varphi(y)\not=1$ since $\alpha_0(y)>0$, so by~\eqref{els0} the
leading coefficient in~\eqref{almostthere} reduces further to
\[
\mathfrak{s}_y^{-\varphi^\vee}k_0^{\eta_e(\varphi(y))}u_0^{\eta_0(\varphi(y))}.
\]
To complete the proof of~\eqref{twt2inv}, it thus suffices to show that
\begin{equation}\label{ldformula}
\mathfrak{s}_y^{\varphi^\vee}=k_0^{\eta_e(\varphi(y))}u_0^{\eta_0(\varphi(y))}\prod_{\alpha\in\Pi(s_\varphi)}\mathbf{k}_\alpha^{\eta_e(\alpha(y))}\mathbf{k}_{\frac{\alpha}{2}}^{\eta_0(\alpha(y))}.
\end{equation}

By the definition of $\mathfrak{s}_y$, we have
\[
\mathfrak{s}_y^{\varphi^\vee}=\prod_{\alpha\in\Phi_0^+}\bigl(\mathbf{k}_\alpha\mathbf{k}_{(\alpha,1)}\bigr)^{\frac{\eta_e(\alpha(y))\alpha(\varphi^\vee)}{2}}
\bigl(\mathbf{k}_{\frac{\alpha}{2}}\mathbf{k}_{(\frac{\alpha}{2},\frac{1}{2})}\bigr)^{\frac{\eta_o(\alpha(y))\alpha(\varphi^\vee)}{2}}.
\]
Consider the decomposition of $\Phi_0^+$ as the disjoint union of the subsets
\[
\Phi_0^+[m]:=\bigl\{\alpha\in\Phi_0^+ \mid \alpha\bigl(\varphi^\vee\bigr)=m\bigr\},\qquad m\in\mathbb{Z}.
\]
We have $\Phi_0^+[m]=\varnothing$ unless $m=0,1,2$, and
\[
\Phi_0^+[0]=\Phi_0^+\setminus\Pi(s_\varphi),\qquad \Phi_0^+[1]=\Pi(s_\varphi)\setminus\{\varphi\},\qquad \Phi_0^+[2]=\{\varphi\}.
\]
Hence
\[
\mathfrak{s}_y^{\varphi^\vee}=(k_0k_r)^{\frac{\eta_e(\varphi(y))}{2}}(u_0u_r)^{\frac{\eta_o(\varphi(y))}{2}}
\prod_{\alpha\in\Pi(s_\varphi)}\bigl(\mathbf{k}_\alpha\mathbf{k}_{(\alpha,1)}\bigr)^{\frac{\eta_e(\alpha(y))}{2}}
\bigl(\mathbf{k}_{\frac{\alpha}{2}}\mathbf{k}_{(\frac{\alpha}{2},\frac{1}{2})}\bigr)^{\frac{\eta_o(\alpha(y))}{2}}.
\]
By~\eqref{equalparameter}, formula~\eqref{ldformula} immediately follows if $\Phi_0$ is not of type ${\rm C}_r$, $r\geq 1$. If
$\Phi_0$ is of type ${\rm C}_r$, $r\geq 1$, then
\[
\Phi_0^+[1]=\Pi(s_\varphi)\setminus\{\varphi\}=\Pi_s(s_\varphi)
\]
with $\Pi_s(s_\varphi)$ the positive {\it short} roots in $\Phi_0$ mapped to
negative roots by $s_\varphi$. Hence \smash{$\mathbf{k}_\alpha=\mathbf{k}_{\frac{\alpha}{2}}$} and \smash{$\mathbf{k}_{\frac{\alpha}{2}}=\mathbf{k}_{(\frac{\alpha}{2},\frac{1}{2})}$}
for $\alpha\in\Phi_0^+[1]=\Pi_s(s_\varphi)$, and we conclude that
\begin{align*}
\mathfrak{s}_y^{\varphi^\vee}&=(k_0k_r)^{\eta_e(\varphi(y))}(u_0u_r)^{\eta_o(\varphi(y))}
\prod_{\alpha\in\Pi_s(s_\varphi)}\bigl(\mathbf{k}_\alpha\mathbf{k}_{(\alpha,1)}\bigr)^{\frac{\eta_e(\alpha(y))}{2}}
\bigl(\mathbf{k}_{\frac{\alpha}{2}}\mathbf{k}_{(\frac{\alpha}{2},\frac{1}{2})}\bigr)^{\frac{\eta_o(\alpha(y))}{2}}\\
&=(k_0k_r)^{\eta_e(\varphi(y))}(u_0u_r)^{\eta_o(\varphi(y))}
\prod_{\alpha\in\Pi_s(s_\varphi)}\mathbf{k}_\alpha^{\eta_e(\alpha(y))}\mathbf{k}_{\frac{\alpha}{2}}^{\eta_o(\alpha(y))}\\
&=k_0^{\eta_e(\varphi(y))}u_0^{\eta_o(\varphi(y))}
\prod_{\alpha\in\Pi(s_\varphi)}\mathbf{k}_\alpha^{\eta_e(\alpha(y))}\mathbf{k}_{\frac{\alpha}{2}}^{\eta_o(\alpha(y))},
\end{align*}
as desired.
\end{proof}

We can now complete the proof of the second line of~\eqref{genfor} (and hence of Proposition~\ref{triangdelta}) as in~\cite[Proposition~5.29]{SSV2}:
let $w\in W^\mathcal{O}$ and fix a reduced expression $w=s_{j_1}s_{j_2}\cdots s_{j_\ell}$. Then $\Pi(w)=\{b_1,\dots,b_\ell\}$ with
\[
b_i:=s_{j_\ell}\cdots s_{j_{i+1}}\alpha_{j_i}
\]
(for $i=\ell$ this should be read as $b_\ell=\alpha_{j_\ell}$). Since $w\in W^\mathcal{O}$, we have $w\Phi_\mathcal{O}^+\subseteq\Phi^+$ with $\Phi_{\mathcal{O}}^+$
defined~by
\[
\Phi_\mathcal{O}^+:=
\Phi^+\cap\biggl(\bigoplus_{j\in J(\mathcal{O})}\mathbb{Z}\alpha_j\biggr),
\]
and hence $b_i\in\Phi^+\setminus\Phi_\mathcal{O}^+$ for $i=1,\dots,\ell$. Since $c^\mathcal{O}\in C_+^\mathcal{O}$, it follows that
\[
\alpha_{j_i}\bigl(s_{j_{i+1}}\cdots s_{j_\ell}c^\mathcal{O}\bigr)=b_i\bigl(c^\mathcal{O}\bigr)>0
\]
for $i=1,\dots,\ell$.
Hence
\[
\pi_{t}^\mathcal{O}\bigl(\widetilde{\delta}(T_{w^{-1}})\bigr)x^{c^\mathcal{O}} =\pi_{t}^\mathcal{O}\bigl(\widetilde{\delta}(T_{j_1})\bigr)\cdots\pi_{t}^\mathcal{O}\bigl(\widetilde{\delta}(T_{j_\ell})\bigr)x^{c^\mathcal{O}}=
k_w\bigl(c^\mathcal{O}\bigr)x^{wc^\mathcal{O}}+ {\rm l.o.t.}
\]
by Lemma~\ref{twtinvlemma} and~\eqref{kcocycle}. This completes the proof of the second line of~\eqref{genfor} (and hence of Proposition~\ref{triangdelta}).
\end{proof}

\subsection*{Acknowledgements}
I thank Siddhartha Sahi and Vidya Venkateswaran for several interesting discussions about the topic of the paper. I also thank the referees for their thoughtful comments.

\pdfbookmark[1]{References}{ref}
\LastPageEnding

\end{document}